\documentstyle{article}

\newcommand{\N}{{\rm I\kern-.5ex N}}
\newcommand{\Z}{{\sf \vrule height 1.55ex depth-1.2ex width.03em\kern-.11em Z \kern-.9ex Z\kern-.11em\vrule height 0.3ex depth0ex width.03em}}
\newcommand{\Q}{{\rm\kern.2ex\vrule height1.55ex depth-.05ex width.03em\kern-.7ex Q}}
\newcommand{\R}{{\rm I\kern-.5ex R}}
\newcommand{\Rvar}{{\rm I\kern-.5ex R}}
\newcommand{\C}{{\rm\kern.3ex\vrule height1.55ex depth-.05ex width.03em\kern-.7ex C}}
\newcommand{\Cvar}{{\, \rm\kern.3ex\vrule height1.1ex depth-.05ex width.03em\kern-.7ex C}}

\newcommand{\golf}{\tilde{\rule{0ex}{1.7ex}}\,}
\newcommand{\streep}{\overline{\rule{0ex}{1.7ex}\hspace{1ex}\rule{0ex}{1.7ex}}}

\newcommand{\spat}{\hspace{4ex}}

\newcommand{\bld}{\text{\rm Ran}\,}
\newcommand{\lan}{\langle}
\newcommand{\ran}{\rangle}

\newcommand{\cK}{{\cal K}}
\newcommand{\cB}{{\cal B}}
\newcommand{\cR}{{\cal R}}
\newcommand{\cL}{{\cal L}}

\newcommand{\od}{\odot}
\newcommand{\ot}{\otimes}

\newcommand{\om}{\omega}
\newcommand{\io}{\iota}
\newcommand{\vfi}{\varphi}
\newcommand{\vep}{\varepsilon}

\newcommand{\al}{\alpha}
\newcommand{\be}{\beta}

\newcommand{\th}{\theta}

\newcommand{\si}{\sigma}

\newcommand{\comp}{\,\rule[.5ex]{.2ex}{.2ex}\,}
\newcommand{\dplus}{\dot{+}}

\newcommand{\text}[1]{\mbox{#1}}

\newcommand{\cst}{\text{C}$^*$}
\newcommand{\wst}{\text{W}$^*$}

\newcommand{\qed}{\ \hfill \rule{2mm}{2mm}}
\newenvironment{demo}{\medskip\noindent\bf Proof :\ \  \rm}{\qed\bigskip\par }

\newtheorem{definition}{Definition}[section]
\newtheorem{proposition}[definition]{Proposition}
\newtheorem{lemma}[definition]{Lemma}
\newtheorem{corollary}[definition]{Corollary}
\newtheorem{remark}[definition]{Remark}
\newtheorem{theorem}[definition]{Theorem}

\newtheorem{notation}[definition]{Notation}
\newtheorem{result}[definition]{Result}

\newtheorem{terminology}[definition]{Terminology}
\newtheorem{question}[definition]{Question}

\begin{document}
\begin{center}
\LARGE\bf The functional calculus of regular operators on Hilbert \cst-modules revisited

\end{center}

\bigskip

\begin{center}
\rm J. Kustermans  \footnote{Research Assistant of the
National Fund for Scientific Research (Belgium)}

Institut for Matematik og Datalogi

Odense Universitet

Campusvej 55

5230 Odense M

Denmark

\bigskip

e-mail : johank@imada.ou.dk

\bigskip\medskip

\bf June 1997 \rm
\end{center}

\bigskip

\subsection*{Abstract}
In \cite{Wor6}, Woronowicz introduced a functional calculus for normal regular operators in a Hilbert \cst-module. In this paper,
we have translated several concepts known in Hilbert space theory to the Hilbert \cst-module framework. We looked for instance into 
the functional calculus for strictly positive elements, the Fuglede Putnam theorem in Hilbert \cst-modules, commuting normal regular operators, ... 

It appears that if we are a little bit careful, most of the Hilbert space results have their analogues in the Hilbert \cst-module case. 
An exception to this rule is of course the polar decomposition of a regular operator.

\section*{Introduction}

Regular operators on Hilbert \cst-modules were studied in \cite{Baa1}. A nice extensive overview concerning Hilbert \cst-modules is given 
in \cite{Lan}. Another standard reference for refular operators is \cite{Wor6}. 

\medskip

There is a lot of interest for regular operators (or elements affiliated with a \cst-algebra) in the \cst-algebraic quantum group scene. 
There are several reasons for this :
\begin{itemize}
\item An interesting object assoiated to a locally compact group is the modular function which connects the left and right Haar measure.
 This modular function is a continuous group homomorphism from the group in the complex numbers.

In the quantum group case, the analogue of this modular function still exists, but it is now a strictly positive element affiliated 
with the \cst-algebra.
\item Some important quantum groups are defined by generators and relations. In the compact case, this generators belong to the \cst-algebra.
 But in the non compact case, these generators can (and will be mostly) elements affiliated to the \cst-algebra (see e.g \cite{Wor6}, \cite{Wor5}).
\end{itemize}

\medskip

Another application of regular operators can be found in the theory of so-called \cst-valued weights (see \cite{Kus1}).

\bigskip

The main aim of this paper is the translation of Hilbert space results to the Hilbert \cst-module framework in order to make the Hilbert 
\cst-module theory a little bit more flexible. Another part of this effort can be found in \cite{Kus}.

As a consequence, you will not find  highly original results but only a collection of the most basic and useful properties. Most of 
the results are proven within the Hilbert \cst-module framework and do not make use of Hilbert space techniques. 

The proofs of the results depend only on the results proven in \cite{Wor6} and \cite{Lan}.

\bigskip

The following topics are included in this paper : 
\begin{enumerate}
\item Regular operators in Hilbert \cst-modules (preliminaries)
\item Left, right and middle multipliers
\item The functional calculus of normal elements  (preliminaries)
\item A further development of the functional calculus
\item Invertible regular operators
\item The functional calculus revisited
\item The Fuglede-Putnam theorem for Hilbert \cst-modules
\item Natural left and right multipliers of a regular operator
\item \cst-algebras of adjointable operators on Hilbert \cst-modules
\item Representing Hilbert \cst-modules on Hilbert spaces
\item Commuting normal operators
\item Tensor products of Hilbert \cst-modules
\end{enumerate}

\bigskip

We end this section with some conventions and notations.

\medskip

If $G$ is a set, we will denote the identity mapping on $G$ by $\io_G$. If it is clear which set is under consideration, we will drop
the subscript $G$. For any function $f$ and any set $G$, we denote the restriction of $f$ to $G$ by the symbol $f_G$.

For any subset $K$ of $\C$, we define $K_0 = K \setminus \{0\}$. 

\medskip

If $E$ is a normed space, we denote the set of bounded operators by $\cB(E)$.

Consider normed spaces $E$, $F$, $G$ and linear mappings $S : E \rightarrow F$ and
$T : F \rightarrow G$. If the composition $T S$ is closable, we will always denote the closure by $S \comp T$. We will however use this 
notation only in well controlled circumstances.

\medskip

All Hilbert \cst-modules  in these paper are right modules over the \cst-algebra in question. The inner valued product will always be 
linear in the first and adjoined linear in the second factor.

Consider Hilbert \cst-module $E$,$F$. Then we define $\cL(E,F)$  as the set of adjointable operators from $E$ into $F$. We will also
use the notation $\cL(E) = \cL(E,E)$.

\bigskip

\section{Regular operators on Hilbert \cst-modules}

Regular operators on Hilbert \cst-modules were studied in \cite{Baa1}. Other standard 
references for regular operators are \cite{Wor6} and \cite{Lan}. The results gathered in this first section come from \cite{Wor6} and \cite{Lan}.

\bigskip

Let us start of with the definition of the adjoint of a densely defined linear operator between
Hilbert \cst-modules. The definition is completely same as the one for operators in Hilbert spaces but we have to use the \cst-algebra 
valued inner product.

\begin{definition}
Consider  Hilbert \cst-modules $E$,$F$ over a \cst-algebra $A$ and let $T$ be a densely
defined $A$-linear mapping  from within $E$ into $F$. Then we define the mapping $T^*$
from within $F$ into $E$ such that the domain of $D(T^*)$ is equal to
$$\{ \, v \in F \mid \text{There exists } w \in E \text{ such that }
\lan T(u) , v \ran = \lan u , w \ran \text{ for every } u \in D(T) \,\}$$
and such that $\lan T(u) , v \ran = \lan u , T^*(v) \ran$ for every $v \in D(T^*)$ and $u
\in D(T)$.
\end{definition}

It is then easy to check that $T^*$ is a closed $A$-linear mapping.

\bigskip\medskip

We can now introduce the definition of a regular operator (see page 96 of \cite{Lan}).

\begin{definition}
Consider  Hilbert \cst-modules $E$,$F$ over a \cst-algebra $A$ and let $T$ be  a densely
defined $A$-linear mapping from within $E$ into $F$ such that $T$ is closed, $T^*$ is densely defined and
$1+T^*T$ has dense range. Then we call $T$ a regular operator from within $E$ into $F$.
\end{definition}

\begin{notation}
Consider  Hilbert \cst-modules $E$,$F$ over a \cst-algebra $A$. Then we denote the set of regular operators from within $E$ into $F$ by $\cR(E,F)$.
\end{notation}

\begin{remark}
Consider a Hilbert \cst-module $E$ over a \cst-algebra $A$. Then regular operators from within $E$ into $E$ are called regular operators in $E$. 
We will also use the notation $\cR(E) = \cR(E,E)$.
\end{remark}

\bigskip

There is an alternative definition of a regular operator between Hilbert \cst-modules (see
definition 1.1 of \cite{Wor6}. It is proven in \cite{Lan} that both are equivalent.

\begin{theorem}
Consider  Hilbert \cst-modules $E$,$F$ over a \cst-algebra $A$ and let $T$ be a densely
defined mapping from within $E$ into $F$. Then $T$ is regular $\Leftrightarrow$ There
exists an element $z \in \cL(E,F)$ with $\|z\| \leq 1$ such that $D(T) = (1 - z^*
z)^\frac{1}{2} \, E$ and $T\bigl((1-z^* z)^\frac{1}{2} \, v) = z \, v$ for every $v \in
E$.
\end{theorem}

\medskip

\begin{remark} \rm
Consider  Hilbert \cst-modules $E$,$F$ over a \cst-algebra $A$.
Then we have the following properties concerning the element $z$ in the definition above.
\begin{itemize}
\item Let $T$ be an element in $\cR(E,F)$. Then there exists a unique element $z \in \cL(E,F)$
with $\|z\| \leq 1$ such that $D(T) = (1 - z^* z)^\frac{1}{2} \, E$ and $T\bigl((1-z^*
z)^\frac{1}{2} \, v) = z \, v$ for every $v \in E$. We use the notation $z_T = z$ and call
$z_T$ the $z$-transform of $T$.
\item Consider two regular operators $S, T \in \cR(E,F)$. Then $S = T$ if and
only if $z_S = z_T$.
\item Let $z$ be an element in $\cL(E,F)$ such that $\|z\| \leq 1$ and such that
$(1-z^* z)^\frac{1}{2}\, E$ is dense in $E$. Then there exists a unique element $T \in
\cR(E,F)$ such that $z_T = z$.
\end{itemize}
This $z$-transform turns out to be very useful in many proofs concerning regular
operators. It allows to transfer problems concerning an unbounded regular operator to a
bounded adjointable operator.
\end{remark}

\medskip

A first result in this respect connects the boundedness of $T$ to a certain property of $z_T$ : 

\begin{result}
Consider  Hilbert \cst-modules $E$,$F$ over a \cst-algebra $A$. Then the following holds :
\begin{itemize}
\item Consider $T \in \cL(E,F)$. Then $T$ is regular and $\|z_T\| < 1$.
\item Consider $T \in \cR(E,F)$. Then $T$ belongs to $\cL(E,F)$ $\Leftrightarrow$ $D(T) = E$
$\Leftrightarrow$ $T$ is bounded $\Leftrightarrow$ $\|z_T\| < 1$.
\end{itemize}
\end{result}

\bigskip

\begin{remark} \rm
An important class of regular operators arises from \cst-algebras. In fact, Woronowicz
looks in \cite{Wor6} only  at this kind of regular operators but most of  his proofs can be
immediately copied to the case of regular operators between Hilbert \cst-modules.

Consider a \cst-algebra $A$ and define $E$ to be the Hilbert \cst-module over $A$ such 
that $E = A$ as a right $A$-module and $\lan a , b \ran = b^* a$ for every $a,b \in A$. 
Then the elements of $\cR(E)$ are called elements affiliated with $A$. We write also 
$T \eta A$ instead of $T \in \cR(E)$.
\end{remark}

\medskip

An important result, proven by Woronowicz in \cite{Wor6}, states that a non-degenerate $^*$-homomorphism can
be extended to the set of affiliated elements :

\begin{theorem} \label{prel1.thm1}
Consider a Hilbert \cst-module $E$ over a \cst-algebra $A$. Let $B$ be a \cst-algebra and
$\pi$ be a non-degenerate $^*$-homomorphism from $B$ into $\cL(E)$.

Consider an element $T$ affiliated with $B$. Then there exists a unique element $S \in
\cR(E)$ such that $z_s = \pi(z_T)$ and we define $S = \pi(T)$.
We have moreover that $\pi(D(T))\,E$ is a core for $\pi(T)$ and $\pi(T)(\pi(b)v)
= \pi(T(b)) \, v$ for every $b \in D(T)$ and $v \in E$.
\end{theorem}

The last part of this theorem implies that $\pi(D) \, K$ is a core for $\pi(T)$ if $D$ is
a core for $T$ and $K$ is a dense subspace of $E$.

\medskip

\begin{remark} \rm
Suppose moreover that $\pi$ is injective. Then the canonical extension of $\pi$ to $M(B)$ is also injective. 
Using the $z$-transform, this implies immediately the following result.

Let $S$ and $T$ be two elements affiliated with $B$. Then $S = T$ if and only if $\pi(S) =
\pi(T)$.
\end{remark}

\medskip

The following result can also be found in \cite{Wor6} (theorem 1.2). It follows easily
using the $z$-transform.

\begin{proposition} \label{prel1.prop2}
Consider a Hilbert \cst-module $E$ over a \cst-algebra $A$. Let $B$,$C$ be two
\cst-algebras. Consider a non-degenerate $^*$-homomorphism $\pi$ from $B$ into $M(C)$ and
a non-degenerate $^*$-homomorphism $\th$ from $C$ into $\cL(E)$. Then $(\th \pi)(T) =
\th(\pi(T))$ for every $T \, \eta \, B$.
\end{proposition}

\bigskip

Concerning the adjoint, we have the following key result (see theorem 1.4 of \cite{Wor6}).

\begin{proposition}
Consider  Hilbert \cst-modules $E$,$F$ over a \cst-algebra $A$ and let $T$ be an element in $\cR(E,F)$. 
Then $T^*$ is a regular operator from within $F$ into $E$. We have moreover that $z_{T^*} = (z_T)^*$  and $T^{**}=T$.
\end{proposition}

\bigskip

We have the usual definitions of selfadjointness, normality and positivity. Normality and selfadjointness were 
already considered in \cite{Wor6}. Positivity seems to be the logical definition.

\begin{definition}
Consider a Hilbert \cst-module $E$ over a \cst-algebra $A$ and let $T$ be an element
in $\cR(E)$. Then we have the following definitions.
\begin{itemize}
\item We call $T$ normal $\Leftrightarrow$ $D(T) = D(T^*)$ and $\lan T v, T v \ran
= \lan T^* v , T^* v \ran$ for every $v \in D(T)$.
\item We call $T$ selfadjoint $\Leftrightarrow$ $T^* = T$
\item We call $T$ positive $\Leftrightarrow$ $T$ is normal and
$\lan T v , v \ran \geq 0$ for every $v \in D(T)$.
\end{itemize}
\end{definition}

\medskip

These definitions behave rather well with respect to the $z$-thransforms (see equivalence 1.15 of \cite{Wor6}) :

\begin{result}
Consider a Hilbert \cst-module $E$ over a \cst-algebra $A$ and let $T$ be an element
in $\cR(E)$. Then we have the following equivalencies.
\begin{itemize}
\item $T$ is normal $\Leftrightarrow$ $z_T$ is normal
\item $T$ is selfadjoint $\Leftrightarrow$ $z_T$ is selfadjoint
\item $T$ is positive $\Leftrightarrow$ $z_T$ is positive
\end{itemize}
\end{result}
\begin{demo}
The result concerning normality is proven in equivalence 1.15 of \cite{Wor6}. Using the fact that $z_{T^*} =
(z_T)^*$, the result concerning the selfadjointness is trivial. So we turn to the
positivity.
\begin{trivlist}
\item[$\ \,\Rightarrow$] Suppose that $T$ is positive. Because $T$ is normal, we have
that $z_T$ is normal.

Choose $v \in E$. We know that $(1-z_T^* z_T)^\frac{1}{2} v$ belongs to $D(T)$ and $T(
(1-z_T^* z_T)^\frac{1}{2} v ) = z_T v$. Therefore, the positivity of $T$ implies that the element  
$\lan z_T v , (1-z_T^* z_T)^\frac{1}{2} v \ran $ is positive.

\medskip

Take $w \in E$, and replace $v$ in the previous equality by $(1-z_T^* z_T)^\frac{1}{4} w$.
So we get that the element \newline $\lan z_T (1-z_T^* z_T)^\frac{1}{4} w , (1-z_T^* z_T)^\frac{3}{4} w \ran$ is positive.

By the normality of $z_T$, we have that
\begin{eqnarray*}
\lan z_T (1 - z_T^* z_T)^\frac{1}{2} w , (1 - z_T^* z_T)^\frac{1}{2} w \rangle
& = & \lan (1 - z_T^* z_T)^\frac{1}{4} z_T (1 - z_T^* z_T)^\frac{1}{4} w ,
(1 - z_t^* z_T)^\frac{1}{2} w \rangle \\
& = & \lan z_T (1 - z_T^* z_T)^\frac{1}{4} w , (1 - z_T^* z_T)^\frac{3}{4} w \rangle
\end{eqnarray*}
So we see that the element $\lan z_T (1 - z_T^* z_T)^\frac{1}{2} w , (1 - z_t^*
z_T)^\frac{1}{2} w \rangle$ is positive.

\medskip

Because $(1 - z_T^* z_T)^\frac{1}{2} E$ is dense in $E$, this implies that $\lan z_T u , u
\ran$ is positive for every $u \in E$. By lemma 4.1 of \cite{Lan}, we get that $z_T$ is positive. 
\item[$\ \,\Leftarrow$] Suppose that $z_T$ is positive. Because $z_T$ is normal,
$T$ will be normal.

Choose $v \in D(T)$. Then there exists $w \in E$ such that $v = (1-z_T^2)^\frac{1}{2} w$.
This implies that
$$ \lan T v , v \ran  =
\lan z_T w , (1 - z_T^2)^\frac{1}{2} w \ran
= \lan (1 - z_T^2)^\frac{1}{4} z_T w , (1 - z_T^2)^\frac{1}{4} w \ran
= \lan z_T  (1 - z_T^2)^\frac{1}{4} w , (1 - z_T^2)^\frac{1}{4} w \ran$$
which is positive by the positivity of $z_T$.
\end{trivlist}
\end{demo}

\medskip

\begin{remark} \rm
This result implies immediately that a positive regular operator is selfadjoint.
\end{remark}

\medskip

We will also need the following notion of strict positivity .

\begin{definition}
Consider a Hilbert \cst-module $E$ over a \cst-algebra $A$ and let $T$ be an element in $\cR(E)$. 
Then we call $T$ strictly positive $\Leftrightarrow$ $T$ is positive and
$T$ has dense range.
\end{definition}

A strictly positive element is automatically injective. The reverse is not true. A counterexample can be 
easily constructed in the commutative case by taking a positive continuous function which is 0 in one point.

Strict positivity is intimately connected with invertibility of a regular operator which will be considered in section \ref{art2}.

\medskip

We have also the following results. The only difficult part of the proof of the first one is resolved by  
proposition 9.9 of \cite{Lan}. The second one is the content of lemma 9.2 of \cite{Lan}.

\begin{proposition}
Consider Hilbert \cst-modules $E$,$F$ over a \cst-algebra $A$ and let $T$ be an element in $\cR(E,F)$. 
Then $T^* T$ and $T T^*$ are positive regular operators in $E$.
\end{proposition}

\begin{proposition} \label{prel1.prop1}
Consider Hilbert \cst-modules $E$,$F$ over a \cst-algebra $A$ and let $T$ be an element in $\cR(E,F)$. 
Then $D(T^* T)$ is a core for $T$.
\end{proposition}

\bigskip

\begin{result} \label{prel1.res1}
Consider a Hilbert \cst-module over a \cst-algebra $A$. Let $B$ be a \cst-algebra and
$\pi$ a non-degenerate $^*$-homomorphism from $B$ into $\cL(E)$. Then we have for every
element $T$ affiliated with $B$ the equalities 
$$\pi(T)^* = \pi(T)^*  \hspace{2cm} \pi(T^* T) = \pi(T)^* \pi(T) \hspace{2cm} \pi(T T^*) = \pi(T) \pi(T)^*$$
\end{result}
\begin{demo}
We have that $z_{\pi(T^*)} = \pi(z_{T^*})= \pi(z_T^*) = \pi(z_T)^* =
 z_{\pi(T)}^* = z_{\pi(T)^*}$. Hence $\pi(T^*) = \pi(T)^*$.

\medskip

Take  $x \in D(T^* T)$ and $v \in E$. Then $\pi(x) v$ belongs to $D(\pi(T^* T))$ and
$\pi(T^* T) (\pi(x) v) = \pi((T^* T)(x)) \, v$.

Because $x$ belongs to $D(T)$, we have that $\pi(x) v$ belongs to $D(\pi(T))$ and that
$\pi(T) (\pi(x) v) = \pi(T(x))\,v$. 
Because $T(x)$ belongs to $D(T^*)$, this implies that $\pi(T) (\pi(x) v)$ belongs to
$D(\pi(T^*))$ and that 
$$\pi(T^*)\bigl(\pi(T) (\pi(x) \, v)\bigr) = \pi(T^*)(\pi(T(x)) v) =
\pi\bigl(T^*(T(x))\bigr) \, v = \pi((T^* T)(x))\,v = \pi(T^* T)\,(\pi(x) v) \ . $$
This implies that $\pi(x) v$ belongs to $D(\pi(T)^* \pi(T))$ and that $(\pi(T)^* \pi(T)) (\pi(x)
v) = \pi(T^* T) (\pi(x) v)$.

We know that $\pi(D(T^* T))\,E$ is a core for $\pi(T^* T)$. So we get that $\pi(T^* T)
\subseteq \pi(T)^* \pi(T)$ by the closedness of $\pi(T)^* \pi(T)$.

Because both are selfadjoint (remember the first equality of this result), we get that $\pi(T^* T) = \pi(T)^* \pi(T)$. 
The other equality follows by symmetry.
\end{demo}

\medskip

Using the $z$-transform, the first three statements of the following proposition are easy
to prove. The last one is then not so difficult to prove.

\begin{proposition} \label{prel1.prop3}
Consider a Hilbert \cst-module $E$ over a \cst-algebra $A$. Let $B$ be a \cst-algebra and
$\pi$ a non-degenerate $^*$-homomorphism from $B$ into $\cL(E)$. Suppose that $T$ is an element affiliated with $B$. 
Then we have the following properties.
\begin{enumerate}
\item If $T$ is normal, then $\pi(T)$ is normal.
\item If $T$ is selfadjoint, then $\pi(T)$ is selfadjoint.
\item If $T$ is positive, then $\pi(T)$ is positive.
\item If $T$ is strictly positive, then $\pi(T)$ is strictly positive.
\end{enumerate}
\end{proposition}

If $\pi$ is moreover injective, then the first three implications become equivalencies (use the $z$-transform). 
This is however not true for the last one!

\bigskip

\section{Left, right and middle multipliers}

This section serves merely to introduce some conventions, notations and basic results concerning
certain compositions of adjointable and regular operators between Hilbert \cst-modules. 

\medskip

Let us start of with some definitions.

\begin{terminology}
Consider  Hilbert \cst-modules $E,F,G$ over a \cst-algebra $A$. 
\begin{enumerate}
\item Consider $T \in \cR(F,G)$ and $x \in \cL(E,F)$. We call $x$ a right multiplier of $T$ $\Leftrightarrow$ 
$T \, x$ is bounded and $\overline{T\,x}$ belongs to $\cL(E,G)$.
\item Consider $T \in \cR(E,F)$ and $x \in \cL(F,G)$. We call $x$ a left multiplier of $T$ $\Leftrightarrow$ 
$x \, T$ is bounded and  $\overline{x\,T}$ belongs to $\cL(E,G)$.
\end{enumerate}
\end{terminology}

\medskip

We will combine this with the following notations.

\begin{notation}
Consider  Hilbert \cst-modules $E, F, G$ over a \cst-algebra $A$.
\begin{enumerate}
\item Consider $T \in \cR(F,G)$ and $x \in \cL(E,F)$ such that $x$ is a right multiplier of $T$. Then we define
$T \comp x = \overline{T\,x}$. So $T \comp  x$ belongs to $\cL(E,G)$.
\item Consider $T \in \cR(E,F)$ and $x \in \cL(F,G)$ such that $x$ is a left multiplier of $T$. Then we define 
$x \comp T = \overline{x\,T}$. So $x \comp T$ belongs to $\cL(E,G)$.
\end{enumerate}
\end{notation}

\medskip

Concerning right multipliers, we have the following useful characterization.

\begin{result}  \label{prel2.res1}
Consider  Hilbert \cst-modules $E, F, G$ over a \cst-algebra $A$. Let $T$ be an element in $\cR(F,G)$ and 
$x$ an element in $\cL(E,F)$. Then
\begin{enumerate}
\item $x$ is a right multiplier of $T$ $\Leftrightarrow$ $x \, E \subseteq  D(T)$.
\item If $x$ is a right multiplier of $T$, then $T \comp x = T \, x$.
\end{enumerate}
\end{result}
\begin{demo}
\begin{itemize}
\item Suppose that $x$ is a right multiplier of $T$. Because $D(\,\overline{T\,x}\,) = E$,
we have that $D(T\,x)$ is dense in $E$.

Choose $v \in E$. Then there exists a sequence $(v_n)_{n=1}^\infty$ in
$D(T\,x)$ such that $(v_n)_{n=1}^\infty$ converges to $v$.

So $(x \, v_n)_{n=1}^\infty$ converges clearly to $x \, v$.

We have for every $n \in \N$ that $x\,v_n$ belongs to $D(T)$ and  that $T(x \, v_n) = (T \comp x) \, v_n$.  
So we get that $(T(x \, v_n))_{n=1}^\infty$ converges to $(T \comp x) \, v$. Hence, the closedness of $T$ 
implies that $x \, v$ belongs to $D(T)$ and $T(x \, v) = (T \comp x)\,v$.

So we see that $x \, E \subseteq D(T)$ and that $T\,x = T \comp x$.

\item Suppose that $x \, E \subseteq D(T)$. Then $T\,x$ is a closed linear operator from
$E$ into $E$ so it must be bounded by the closed graph theorem. Then we have of course that
$\overline{T\,x} = T\,x$.

Because $T\,x$ is bounded, the operator $(T\,x)^*$ is also bounded. We have also that
$(T\,x)^*$ is closed. It is not difficult to check that $x^*  T^* \subseteq
(T\,x)^*$, so $(T\,x)^*$ has a dense domain. Combining these three facts, we get that $D((T\,x)^*) = E$. But
this implies that $T\,x$ belongs to $\cL(E,F)$. So $x$ is by definition a right multiplier of $T$.
\end{itemize}
\end{demo}

\medskip

Another useful terminology is the following one :

\begin{terminology}
Consider Hilbert \cst-modules $E, F, G, H$ over a \cst-algebra $A$.
Let $S$ be an element in $\cR(E,F)$ and $T$ an element in $\cR(G,H)$. 
Consider $x \in \cL(F,G)$. Then we call $x$ a middle multiplier of \newline $S$,
$T$ if and only if $x(\text{Ran}\,T) \subseteq D(S)$, $S\,x\,T$ is bounded and
$\overline{S\,x\,T}$ belongs to $\cL(E,H)$.
\end{terminology}

\begin{notation}
Consider Hilbert \cst-modules $E, F, G, H$ over a \cst-algebra $A$.
Let $S$ be an element in $\cR(E,F)$ and $T$ an element in $\cR(G,H)$. 
Consider $x \in \cL(F,G)$ such that $x$ is a middle multiplier of $S$,$T$. 
Then we define $S \comp x \comp T =  \overline{S\,x\,T}$. So $S \comp x \comp T$
belongs to $\cL(E,H)$.
\end{notation}

\medskip

\begin{remark} \rm \label{prel2.rem1}
Consider  Hilbert \cst-modules $E$,$F$,$G$ over a \cst-algebra $A$ and let $S$ be an element in $\cR(F,G)$.
Consider an element $x$ in $\cL(E,F)$. Then we have clearly the following
properties :
\begin{itemize}
\item $x$ is a right multiplier of $S$ $\Leftrightarrow$ $x$ is a middle multiplier of $S$, $1$.
\item If $x$ is a right multiplier of $S$, then $S \comp x = S \comp x \comp 1$.
\end{itemize}
A similar remark applies of course also for left multipliers.
\end{remark}

\medskip

Concerning the adjoint operation, we have the following expected result.

\begin{result}
Consider Hilbert \cst-modules $E, F, G, H$ over a \cst-algebra $A$.
Let $S$ be an element in $\cR(E,F)$ and $T$ an element in $\cR(G,H)$. Consider $x \in \cL(F,G)$. Then we have the 
following properties.
\begin{itemize}
\item $x$ is a middle multiplier of $S$,$T$ $\Leftrightarrow$ $x^*$ is a middle muliplier
of $T^*$,$S^*$
\item If $x$ is a middle multiplier of $S$,$T$, then $(S \comp x \comp T)^* = T^* \comp x^* \comp S^*$.
\end{itemize}
\end{result}
\begin{demo}
Suppose that $x$ is a middle multiplier of $S$,$T$.

Take $v \in D(S^*)$. Then we have for every $w \in D(T)$ that
$$\lan x^* S^*(v) , T(w) \ran = \lan S^*(v) , x T(w) \ran = \lan v , S(x T(w)) \ran
= \lan v , (S \comp x \comp T) w \ran = \lan (S \comp  x \comp T)^* v , w \ran $$
which implies that $x^* S^*(v)$ belongs to $D(T^*)$ and that $T^*( x^* S^*(v)) = (S \comp x \comp T)^* v$.

So we get by definition that $x^*$ is a middle multiplier of $T^*$,$S^*$ and $T^* \comp  x^* \comp  S^* = (S \comp x \comp T)^*$.

If $x^*$ is a middle multiplier of $T^*$,$S^*$, we get in a similar way that $x$ is a
middle multiplier of $S$,$T$.
\end{demo}

\medskip

A special case hereof is the following corollary (see also remark \ref{prel2.rem1}).

\begin{corollary} \label{prel2.cor1}
Consider  Hilbert \cst-modules $E, F, G$ over a \cst-algebra $A$ and let $T$ be an element in $\cR(E,F)$
and $x$  an element in $\cL(F,G)$. Then we have the following properties.
\begin{itemize}
\item $x$ is a left multiplier of $T$ $\Leftrightarrow$ $x^*$ is a right multiplier of $T^*$
\item If $x$ is a left multiplier of $T$, then $(x \comp T)^* = T^* \comp x^*$.
\end{itemize}
\end{corollary}

\medskip

\begin{remark} \rm
We have all kinds of multiplication properties (each of which is easy to prove). For instance,  we have the following one.

Consider Hilbert \cst-modules $E$,$F$,$G$,$H$ over a \cst-algebra $A$. Let $T \in \cR(E,F)$, $y \in \cL(F,G)$ and 
$x \in \cL(G,H)$ such that $y$ is a left multiplier of $T$.  
Then $x y$ is a left multiplier of $T$ and $(x y) \comp T = x \, (y \comp T)$.
\end{remark}

\bigskip

\begin{proposition}
Consider a Hilbert \cst-module $E$ over a \cst-algebra $A$. Let $B$ be a \cst-algebra and
$\pi$ a non-degenerate $^*$-homomorphism from $B$ into $\cL(E)$. Consider $S,T \in \cR(B)$ and $x \in \cL(B)$ such that 
$x$ is a middle multiplier of $S$,$T$. 
Then $\pi(x)$ is a middle multiplier of $S$,$T$ and
$\pi(S \comp x \comp T) = \pi(S) \comp \pi(x) \comp \pi(T)$.
\end{proposition}
\begin{demo}
Choose $a \in D(T)$ and $v \in E$. Then $\pi(a) v \in D(\pi(T))$ and
$\pi(T)(\pi(a) v) = \pi(T(a)) \, v$.

Because $x$ is a middle multiplier of $S$,$T$ , we have that $x\,T(a)$ belongs to $D(S)$ and
$S(x \, T(a)) = (S \comp x \comp T)\,a$.
This implies that $\pi(x \, T(a)) v \in D(\pi(S))$ and
$$\pi(S)(\pi(x \, T(a))v) = \pi(S(x \, T(a))) \, v = \pi((S \comp x \comp T)\,a)\, v
= \pi(S \comp x \comp T)\,\pi(a) v$$
So we get that $\pi(x)(\pi(T) \, (\pi(a) v))$ belongs to $D(\pi(S))$ and
$\pi(S)\bigl(\pi(x)(\pi(T) \, (\pi(a) v))\bigr) =  \pi(S \comp x \comp T)\,\pi(a) v$.

We know that $\pi(D(T)) \, E$ is a core for $\pi(T)$. Combining this with the closedness of $\pi(S)$, the last result 
implies easily for every $w \in D(\pi(T))$ that
$\pi(x) ( \pi(T) w) \in D(\pi(S))$ and  $\pi(S)\bigl(\pi(x) ( \pi(T) w)\bigr) = 
\pi(S \comp x \comp T) \, w$.

So we have also that $\pi(S) \pi(x) \pi(T)$ is bounded and
$\overline{\pi(S) \pi(x) \pi(T)} = \pi(S \comp x \comp T)$ which belongs to $\cL(E)$.

Consequently is $\pi(x)$ a middle multiplier of $\pi(S)$,$\pi(T)$ and
$\pi(S) \comp \pi(x) \comp \pi(T) = \pi(S \comp x \comp T)$.
\end{demo}

\medskip

\begin{corollary}
Consider a Hilbert \cst-module $E$ over a \cst-algebra $A$. Let $B$ be a \cst-algebra and
$\pi$ a non-degenerate $^*$-homomorphism from $B$ into $\cL(E)$. Consider moreover $T \in \cR(B)$ and $x \in \cL(E)$.
\begin{itemize}
\item If $x$ is a right multiplier of $T$, then $\pi(x)$ is a right multiplier of $\pi(T)$
and $\pi(T) \comp \pi(x) = \pi(T \comp x)$.
\item If $x$ is a left multiplier of $T$, then $\pi(x)$ is a left multiplier of $\pi(T)$
and $\pi(x) \comp \pi(T) = \pi(x \comp T)$.
\end{itemize}
\end{corollary}

\bigskip

\section{The functional calculus of normal elements}

In \cite{Wor6}, Woronowicz introduced the functional calculus of normal regular operators. We will give an overview of 
the most important results.

\medskip

The main result is the content of the next theorem (see theorem 1.5 of \cite{Wor6}).

\begin{theorem}
Consider a Hilbert \cst-module  $E$ over a \cst-algebra $A$ and let $T$ be a normal element in $\cR(E)$. Then there 
exists a unique non-degenerate $^*$-homomorphism $\vfi_T$ from
$\text{C}_0(\C)$ into $\cL(E)$ such that $\vfi_T(\io_\Cvar) = T$.
\end{theorem}

\medskip

If $T$ belongs to $\cL(E)$, we have already a notion of a functional calculus and spectrum
for $T$ (looked upon as an element in the \cst-algebra $\cL(E)$).

\begin{lemma} \label{prel3.lem1}
Consider a Hilbert \cst-module  $E$ over a \cst-algebra $A$. Let $T$ be a normal element in
$\cL(E)$. Then we have the following properties.
\begin{itemize}
\item We have for every $f \in \text{C}(\C)$ that $\vfi_T(f) = f(T)$.
\item The spectrum $\si(T)$ is equal to the set $\{ \, c \in \C \mid f(c) = 0 \text{ for all }
f \in \ker \vfi_T \, \}$
\end{itemize}
\end{lemma}
\begin{demo}
\begin{itemize}
\item Define the mapping $\pi$ from $\text{C}_0(\C)$ into $\cL(E)$ such
that $\pi(f) = f_{\si(T)}(T)$ for every $f \in \text{C}_0(\C)$. As a
composition of non-degenerate $^*$-homomorphisms, $\pi$ is non-degenerate.

Using proposition \ref{prel1.prop2}, it is not so difficult to check that $\pi(f) =
f_{\si(T)}(T)$ for every $f \in \text{C}(\si(T))$.

We have in particular that $\pi(\io_\Cvar) = \io_{\si(T)}(T) = T$. So $\vfi_T = \pi$.
\item The inclusion $\subseteq$ follows immediately from the previous property.
The inclusion $\supseteq$ requires the lemma of Urysohn ($\si(T)$ is compact!).
\end{itemize}
\end{demo}

\medskip

Woronowicz introduced the following consistent definition.

\begin{definition}
Consider a Hilbert \cst-module  $E$ over a \cst-algebra $A$. Let $T$ be a normal 
element in $\cR(E)$. Then the spectrum of $T$ is defined to be set
$$\si(T) = \{ \, c \in \C \mid f(c) = 0 \text{ for all } f \in \ker \vfi_T \, \} \ .$$
\end{definition}

It is clear that $\si(T)$ is a closed subset of $\,\C$. Because $\vfi_T$ is non-degenerate,
$\vfi_T \neq 0$ so $\ker \vfi_T \neq \text{C}_0(\C)$. We have also that $\text{C}_0(\si(T)) \cong 
\text{C}_0(\C) / \ker \vfi_T$, which implies that $\si(T) \neq \emptyset$.

\bigskip

Using this spectrum, we have the following theorem (theorem 1.6 of \cite{Wor6}).

\begin{theorem}
Consider a Hilbert \cst-module  $E$ over a \cst-algebra $A$. Let $T$ be a normal 
element in $\cR(E)$. Then there exists a unique non-degenerate injective  $^*$-homomorphism
$\psi_T$ from  $\text{\rm C}_0(\si(T))$ into $\cL(E)$ such that $\psi_T(\io_{\si(T)}) = T$.
\end{theorem}

\medskip

If $T$ belongs to $\cL(E)$, we have immediately that $\psi_T(f) = f(T)$ for every $f \in
\text{C}(\si(T))$. So we get in this case  for every continuous function $f$ from within
$\,\C$ into $\,\C$ which is defined on $\si(T)$ that \newline $f(T) = f_{\si(T)}(T) =
\psi_T(f_{\si(T)})$. This remark justifies the following notation.

\medskip

\begin{notation}
Consider a Hilbert \cst-module  $E$ over a \cst-algebra $A$. Let $T$ be a normal 
element in $\cR(E)$. For every continuous function $f$ from within $\,\C$ into $\,\C$ which is
defined on $\si(T)$, we define the element $f(T) = \psi_T(f_{\si(T)})$ , so
$f(T)$ is a normal element in $\cR(E)$.
\end{notation}

If $f$ is bounded on $\si(T)$, we see that $f(T)$ belongs to $\cL(E)$.

\medskip

We have for every $f \in \text{C}(\si(T))$, that $f(T) = \psi_T(f)$.
This implies immediately the usual functional calculus rules for
elements in $\text{C}_b(\si(T))$.

\medskip

\begin{lemma}
Consider a Hilbert \cst-module  $E$ over a \cst-algebra $A$. Let $T$ be a normal 
element in $\cR(E)$. We have for every $f \in \text{C}(\C)$ that $\vfi_T(f) = f(T)$.
\end{lemma}
\begin{demo}
Define the non-degenerate $^*$-homomorphism $\pi$ from $\text{C}_0(\C)$ into
$\text{C}_b(\si(T))$ such that $\pi(f) =$   $f_{\si(T)}$ for every $f \in
\text{C}_0(\C)$. Then $\pi(f) = f_{\si(T)}$ for every $f \in \text{C}(\C)$.

We have that  $\psi_T \pi$ is a non-degenerate $^*$-homomorphism from $\text{C}_0(\C)$
into $\cL(E)$ such that $(\psi_T \pi)(\io_\Cvar) = \psi_T(\io_{\si(T)}) = T$. But this
implies that $\psi_T \pi = \vfi_T$, so we have for every $f \in \text{C}(\C)$ that
$$\vfi_T(f) = \psi_T(\pi(f)) = \psi_T(f_{\si(T)}) = f_{\si(T)}(T)
= f(T) \ . $$
\end{demo}

\medskip

The spectrum of $T$ and $z_T$ are intimately connected (See the remarks at the bottom of 
page 408 of \cite{Wor6}).

\begin{result} \label{prel3.res1}
Consider a Hilbert \cst-module  $E$ over a \cst-algebra $A$. Let $T$ be a normal 
element in $\cR(E)$. Then
\begin{enumerate}
\item $\si(T) = \{\,\frac{c}{(1-|c|^2)^\frac{1}{2}} \mid c \in \si(z_T) \text{ with }
|c| < 1 \,\}$.
\item $\si(z_T) = \text{ the closure of } \{ \frac{c}{\,(1+|c|^2)^\frac{1}{2}} \mid
c \in \si(T) \, \} $.
\end{enumerate}
\end{result}
\begin{demo}
The proof of this result depends on the proof of theorem 1.5 of \cite{Wor6}. Call $D$
the closed unit disk in the complex plane, $D^0$ its interior and $\partial D$ its boundary. Define the function
$J$ from $\,\C$ into $\,\C$ such that $J(c) = \frac{c}{(1 + |c|^2)^\frac{1}{2}}$ for every $c
\in \C$. Then $J$ is a homeomorphism from $\,\C$ to $D^0$ such that $J^{-1}(c) =
\frac{c}{(1-|c|^2)^\frac{1}{2}}$ for every $c$ in $D^0$.

Some reflection allows us  to conclude that
$$\text{C}_0(\C) = \{ \, f \circ J \mid f \in \text{C}(D) \text{ such that  } f=0
\text{ on  } \partial D \, \} \ . \text{\ \ \ \ \ \ \ (a)} $$

Looking at the proof of theorem 1.5 of \cite{Wor6} (equality 1.19), we know that
$\vfi_T(f\!\circ\!J) = f(z_T)$ for every $f \in \text{C}(D)$  \ \ \ \ \ \ \ (b)

Because $\si(z_T) \subseteq D$, we have, similar to the equation in lemma \ref{prel3.lem1},
that 
$$\si(z_T) = \{ \, d \in D \mid f(d) = 0 \text{ for all } f \in \text{C}(D)
\text{ such that } f(z_T)=0 \, \} \text{\ \ \ \ \ \ \ (c)}$$ which implies that
the set $\si(z_T) \cap D^0$ is equal to the set
$$\{\, d \in D^0 \mid f(d)= 0 \text{ for all } f \in \text{C}(D)
\text{ such that } f(z_T)=0 \text{ and } 
  f = 0 \text{ on  } \partial D\, \} \ .$$
For inclusion $\supseteq$, use the fact that for every $d \in D^0$ there exists an
element $g \in \text{C}(D)$ such that $g(d)=1$ and $g=0$ on $\partial D$.

Combining this equality, the definition of $\si(T)$ and equalities (a) and  (b), it is now
easy to prove that  $J(\si(T)) = \si(z_T) \cap D^0$.

\medskip

The previous equation implies that $J(\si(T)) \subseteq \si(z_T)$, so $\overline{J(\si(T))} \subseteq \si(z_T)$.

Choose $d \in \si(z_T)$. Take $f \in \text{C}(D)$ such that $f=0$ on $J(\si(T))$. Then
$f \circ J = 0$ on $\si(T)$, which implies that $\vfi_T(f \circ J) = 0$ by the lemma
before this result. So $f(z_T) = \vfi_T(f \circ J) = 0$. This implies that $f(d)=0$ by
(c). So the lemma of Urysohn implies that $d$ belongs to $\overline{J(\si(T))}$.

Consequently we get that $\si(z_T) = \overline{J(\si(T))}$.
\end{demo}

\medskip

\begin{corollary}
Consider a Hilbert \cst-module  $E$ over a \cst-algebra $A$. Let $T$ be a normal 
element in $\cR(E)$. Then
\begin{enumerate}
\item $T$ is selfadjoint $\Leftrightarrow$ $\si(T) \subseteq \R$.
\item $T$ is positive $\Leftrightarrow$ $\si(T) \subset \R^+$.
\end{enumerate}
\end{corollary}
\begin{demo}
If $T$ is selfadjoint, we know that $z_T$ is selfadjoint which implies that $\si(z_T)
\subseteq \R$, therefore $\si(T) \subseteq \R$. If $\si(T) \subseteq \R$, then
$\io_{\si(T)}$ is real valued which implies that $T = \psi_T(\io_{\si(T)})$ is
selfadjoint.

The result about positivity is proven in a similar way.
\end{demo}

So we see that the usual results about selfadjointness and positivity
remain true.

\bigskip

Also, the following result holds (remarks concerning equality 1.21 of \cite{Wor6})

\begin{proposition} \label{prel3.prop1}
Consider a Hilbert \cst-module  $E$ over a \cst-algebra $A$. Let $B$ be a \cst-algebra and $\pi$ a 
non-degenerate $^*$-homomorphism from $B$ into $\cL(E)$. Let $T$ be a normal element
affiliated with $B$. Then $\pi(T)$ is a normal element in $\cR(E)$  such that
$\si(\pi(T)) \subseteq \si(T)$. We have moreover that $f(\pi(T)) = \pi(f(T))$ for every
continuous function $f$ from within $\,\C$ into $\,\C$ which is defined on $\si(T)$.
\end{proposition}

If $\pi$ is injective, then we have even that $\si(\pi(T)) = \si(T)$.
This follows easily from result \ref{prel3.res1} because in this case $\overline{\pi}$ is injective which 
implies that $\si(z_{\pi(T)}) = \si(\pi(z_T)) = \si(z_T)$.

\medskip

The spectrum of a continuous function on a locally compact space is the closure of its
image. Combining this with the previous remark gives

\begin{result} 
Consider a Hilbert \cst-module  $E$ over a \cst-algebra $A$. Let $T$ be a normal 
element in $\cR(E)$. Let $f$ be a continuous function from within $\,\C$ into $\,\C$ which is
defined on $\si(T)$. Then $\si(f(T)) = \overline{f(\si(T))}$.
\end{result}

The following result is a special case of proposition \ref{prel3.prop1} with $\pi$ equal to 
$\psi_T$.

\begin{proposition}
Consider a Hilbert \cst-module  $E$ over a \cst-algebra $A$. Let $T$ be a normal 
element in $\cR(E)$. Let $f$ be a continuous function from within $\,\C$ into $\,\C$ defined on
$\si(T)$ and $g$ a continuous function from within $\,\C$ into $\,\C$ defined on $\si(f(T))$.
Then $g \circ f$ is defined on $\si(T)$ and $(g \circ f)(T) = g(f(T))$.
\end{proposition}

\bigskip

We end this section with a familiar alternative condition for normality.

\begin{proposition}
Consider a Hilbert \cst-module  $E$ over a \cst-algebra $A$. Let $T$ be an element in $\cR(E)$. 
Then $T$ is normal $\Leftrightarrow$ $T^* T = T T^*$.
\end{proposition}
\begin{demo}
\begin{itemize}
\item Suppose that $T$ is normal. Define the function $f$ from $\si(T)$ into $\,\C$ such that
$f(c) = c$ for every $c \in \si(T)$. So $\psi_T(f) = T$.

It is clear that $f^* f = f f^*$, so result \ref{prel1.res1} implies that
$$T^* T = \psi_T(f)^* \psi_T(f) = \psi_T(f^* f) = \psi_T(f f^*) = \psi_T(f) \psi_T(f^*) = T T^* \ . $$
\item Suppose that $T^* T = T T^*$. By  proposition \ref{prel1.prop1}, we have in this case  that $D(T^* T)$  
is a core for both $T^* T$ and $T T^*$.

Because $T^* T = T T^*$, we have moreover that $\lan T v , T v \ran = \lan T^* v , T^*v \ran$ for every $v \in D(T^* T)$.

Using the fact that $D(T^* T)$ is a core for both $T^* T$ and $T T^*$, we get now easily that $D(T) = D(T^*)$. 
It is now also easy to prove that $\lan T v , T v \ran = \lan T^* v , T^* v \ran$ for every $v \in D(T)$.
\end{itemize}
\end{demo}

\bigskip

\section{A further development of the functional calculus} \label{art1}

In this section we are going to prove some results which makes it possible to develop the functional calculus for 
normal elements a little bit further than in the previous section. 

Essentially, we are going to say something about calculation rules involving (unbounded) continuous functions 
of normal affiliated elements. 

We will work a little bit more general than the situation described above because this will be useful to us in later applications.

\bigskip

Except for the last result in this section, we will fix a Hilbert \cst-module $E$ over a \cst-algebra A, a locally compact 
space $X$ and a non-degenerate $^*$-homomorphism $\pi$ from $\text{C}_0(X)$ into $\cL(E)$.

\medskip

If $f$ belongs to $\text{C}(X)$, then $f$ is affiliated to $\text{C}_0(X)$ so we get a normal element $\pi(f)$ in $\cR(E)$.

\bigskip

We have the following easy lemma.

\begin{lemma}
Consider $f \in \text{C}(X)$.
\begin{enumerate}
\item We have for every $g \in \text{K}(X)$ that $g$ belongs to $D(f)$.
\item $\text{K}(X)$ is a core for $f$.
\end{enumerate}
\end{lemma}
\begin{demo}
By \cite{Wor6}, we know that $D(f) = \{ \, g \in \text{C}_0(X) \mid f g  \in
\text{C}_0(X) \, \}$ and $f(g) = f g$ for every $g \in D(f)$.

We have for every $g \in \text{K}(X)$ that $f g$ belongs to $\text{K}(X)$, so
$g$ belongs to $D(f)$ by the previous remark.

Take an approximate unit $(e_k)_{k \in K}$ for $\text{C}_0(X)$ in $\text{K}(X)$. 

Choose $h \in D(f)$. We have for every $k \in K$ that $h e_k$ belongs to
$\text{K}(X)$ and $f(h e_k) = f h e_k$. Because $h$ and $f h$ belong to
$\text{C}_0(X)$ we see that $(h e_k)_{k \in K}$ converges to $h$ and that $(f(h e_k))_{k
\in K}$ converges to $f h = f(h)$.
\end{demo}

\medskip

By theorem \ref{prel1.thm1} and the remark after it, this implies immediately the following lemma.

\begin{lemma}  \label{art1.lem3}
Consider $f \in \text{C}(X)$. Then the following properties hold.
\begin{itemize}
\item We have for every $g \in \text{K}(X)$ and  $v \in E$ that
$\pi(g) v$ belongs to $D(\pi(f))$ and $\pi(f)(\pi(g) v) = \pi(f g) \,v$.
\item The space $\pi(\text{K}(X)) E$ is a core for $\pi(f)$.
\end{itemize}
\end{lemma}

\medskip

\begin{lemma} \label{art1.lem1}
Consider $f \in \text{C}(X)$, $g \in \text{K}(X)$. Then
$\pi(g)$ is a left and right multiplier of $\pi(f)$ 
and $$ \pi(g)  \pi(f) \subseteq \pi(f)  \pi(g) = \pi(f g)$$
\end{lemma}
\begin{demo}
The previous lemma implies by result \ref{prel2.res1} that $\pi(g)$ is a right multiplier of $\pi(f)$ and that 
$\pi(f) \comp \pi(g) = \pi(f g)$.

The previous lemma implies of course also that $\pi(\overline{g})$ is a right multiplier of $\pi(\overline{f})$ 
and $\pi(\overline{f}) \comp \pi(\overline{g}) = \pi(\overline{f g})$.

Hence, corollary \ref{prel2.cor1} implies that $\pi(g)$ is a left multiplier of $\pi(f)$ and
$\pi(g) \comp \pi(f) = \pi(f g)$.
\end{demo}

\medskip

This lemma implies, as usual, the following result 

\begin{result} \label{art1.res1}
Consider a net $(e_k)_{k \in K}$ in $\text{K}(X)$ which is bounded and converges strictly to 1.
Let $f$ be an element in $\text{C}(X)$. Then we have the following properties.
\begin{enumerate}
\item We have for every $v \in D(\pi(f))$ that the net $(\pi(f e_k)\,v)_{k \in K}$ converges to $\pi(f) \,v$. 
\item Let $v \in E$. Then $v$ belongs to $D(\pi(f))$ $\Leftrightarrow$ the net $( \pi(f e_k)
\, v )_{k \in K}$ is convergent in $E$.
\end{enumerate}
\end{result}

Like in the Hilbert space case (see e.g. \cite{Stra}), we will use this truncating net $(\pi(e_k))_{k \in K}$ to 
get more detail into  the functional calculus of unbounded continuous functions. This result will be freely used in the sequel. 

\medskip

\begin{result} \label{art1.res2}
Consider $f_1,\ldots\!,f_n \in \text{C}(X)$. Then $\pi(f_1) \ldots \pi(f_n)$ is closable
and we have the equality $$\overline{\pi(f_1) \ldots \pi(f_n)} = \pi(f_1 \ldots f_n)$$
\end{result}
\begin{demo}
Take an approximate unit $(e_k)_{k \in K}$ for $\text{C}_0(X)$ in $\text{K}(X)$. 

Using lemma \ref{art1.lem1} and induction, it is easy to prove that
$$\pi(e_k) \pi(f_1) \ldots \pi(f_n) \subseteq \pi(f_1) \ldots \pi(f_n) \pi(e_n) = \pi(f_1 \ldots f_n \, e_k) 
\ \in \cL(E)$$ for $k \in K$.

\begin{itemize}
\item Choose $v \in D(\pi(f_1) \ldots \pi(f_n))$. By lemma \ref{art1.lem1} and the relation above,  we have for 
every $k \in K$ that $\pi(e_k) v \in D(\pi(f_1 \ldots f_n))$ and that
$$\pi(f_1 \ldots f_n) \, \pi(e_k)v = \pi(e_k) \,\, (\pi(f_1) \ldots \pi(f_n))\, v $$
This implies that the net $\bigl(\,\pi(f_1 \ldots f_n)\,\pi(e_k)v \bigl)_{k \in K}$ converges to $(\pi(f_1) \ldots 
\pi(f_n))\,v$. Combining this with the closedness of $\pi(f_1 \ldots f_n)$
and the fact that $(\pi(e_k)v)_{k \in K}$ converges to $v$, we see that
$v \in D(\pi(f_1 \ldots f_n))$ and that $\pi(f_1 \ldots f_n) \, v 
= (\pi(f_1) \ldots \pi(f_n))\,v$.

So we see that $\pi(f_1) \ldots \pi(f_n) \subseteq \pi(f_1 \ldots f_n)$ which implies that
$\pi(f_1) \ldots \pi(f_n)$ is closable and that
$\overline{\pi(f_1) \ldots \pi(f_n)} \subseteq \pi(f_1 \ldots f_n)$.
\item Choose $v \in  D(\pi(f_1 \ldots f_n))$. By the relation above, we have for every $k \in K$ that 
$\pi(e_k)\,v \in D(\pi(f_1) \ldots \pi(f_n))$. 

We know by result \ref{art1.res1} also that
the net $\bigl(\pi(f_1 \ldots f_n) \, \pi(e_k)v\bigr)_{k \in K}$ converges to $\pi(f_1 \ldots f_n)\,v$. 
We have of course also that the net $(\pi(e_k)\,v)_{k \in K}$ converges to $v$.

So we see that $D(\pi(f_1) \ldots \pi(f_n))$ is a core for $\pi(f_1 \ldots f_n)$.
\end{itemize}
\end{demo}

\medskip

We will need the following notation later on.

\begin{notation} \label{art1.not1}
Consider $f_1,\ldots\!,f_n \in \text{C}(X)$. Then we define the element $\pi(f_1) \comp \ldots \comp \pi(f_n)$ 
as the closure of $\pi(f_1) \ldots \pi(f_n)$. So $\pi(f_1) \comp \,\ldots\, \comp \pi(f_n) = \pi(f_1 \ldots f_n)$.
\end{notation}

\medskip

\begin{result} \label{art1.res3}
Consider $f,g \in \text{C}(X)$. Then $D(\pi(f) \pi(g)) = D(\pi(g)) \cap D(\pi(f g))$.
\end{result}
\begin{demo}
Take an approximate unit $(e_k)_{k \in K}$ for $\text{C}_0(X)$ in $\text{K}(X)$. 

Because $\pi(f) \pi(g) \subseteq \pi(f g)$, we have immediately that $D(\pi(f) \pi(g))$ is a subset of 
$D(\pi(f g)) \cap D(\pi(g))$.

Choose $v \in D(\pi(f g)) \cap D(\pi(g))$. Because $v \in D(\pi(f g))$, we have that the net
$(\pi(f g e_k) v)_{k \in K}$ converges to $\pi(f g) v$.

Lemma \ref{art1.lem1} implies that $\pi(f e_k)( \pi(g)  v) = \pi(f g e_k) v$ for every $k \in
K$. Consequently, the net \newline $\bigl(\pi(f e_k)( \pi(g) v )\bigr)_{k \in K}$ is convergent in $E$, 
implying that $\pi(g) v$ belongs to $D(\pi(f))$.

Hence,  $D(\pi(f) \pi(g)) = D(\pi(f g)) \cap D(\pi(g))$.
\end{demo}

\begin{corollary}
Consider $f \in \text{C}(X)$ and $g \in \text{C}_b(X)$. Then
$\pi(f) \pi(g) = \pi(f g)$.
\end{corollary}

\begin{corollary}
Consider $f \in \text{C}(X)$ and $g \in \text{C}_b(X)$ such that $f g \in \text{C}_b(X)$. Then $\pi(g)$ 
is a left and a right multiplier of $\pi(f)$ and
$$ \pi(g)  \pi(f) \subseteq \pi(f)  \pi(g) = \pi(f g)$$
\end{corollary}

\bigskip

The proofs of the following results are based on the same principle as the proof above. Therefore we leave them out.

\begin{result} \label{art1.res4}
Consider $f_1,\ldots\!,f_n \in \text{C}(X)$. Then $\pi(f_1) + \ldots + \pi(f_n)$ is closable
and we have the equality $$\overline{\pi(f_1) + \ldots + \pi(f_n)} = \pi(f_1 + \ldots + f_n)$$
\end{result}

\begin{corollary} \label{art1.cor1}
Consider $f,g \in \text{C}(X)$ such that $f$ or $g$ is bounded. Then $\pi(f)+\pi(g) =
\pi(f+g)$.
\end{corollary}

\medskip

\begin{result}
Consider $f \in \text{C}(X)$ and $c \in \C_0$. Then $\pi(c f) = c \pi(f)$.
\end{result}

\medskip

The next result is true because taking the adjoint commutes with every
non-degenerate $^*$-homomorphism.

\begin{result}
Consider $f \in \text{C}(X)$. Then $\pi(f)^* = \pi(\overline{f})$.
\end{result}

\medskip

\begin{result}
Let $f \in \text{C}(X)$. Then $\pi(f)^* \pi(f) = \pi(f) \pi(f)^* = \pi(|f|^2)$.
\end{result}
\begin{demo}
We know that $\pi(f)^* \pi(f)$ is  closed. Hence,
$$\pi(|f|^2) = \pi(\,\overline{f} f) = \overline{ \pi(\,\overline{f}\,) \pi(f) }
= \overline{\pi(f)^* \pi(f)} = \pi(f)^* \pi(f) \ . $$
Similarly, we have that  $\pi(|f|^2) = \pi(f) \pi(f)^*$.
\end{demo}

\bigskip

The next three results will be very useful to us in later sections.

\begin{lemma} \label{art1.lem2}
Consider $f,g \in \text{C}(X)$ such that there exists a positive number $r$ and a
compact subset $M$ of $X$ such that $|f(x)| \leq r \, |g(x)|$ for every $x \in X
\setminus M$. Then $D(\pi(g)) \subset D(\pi(f))$.
\end{lemma}
\begin{demo}
Take an approximate unit $(e_k)_{k \in K}$ for $\text{C}_0(X)$ in $\text{K}(X)$ such that $e_k = 1$ 
on $M$ for every $k \in K$ (This is possible because of the lemma of Urysohn).

\medskip

Choose $v \in D(\pi(g))$.

Take $k, l \in K$.  Because $e_k =1$ on $M$, $e_l =1$ on $M$ and because we
have that $|f(c)| \leq r \, |g(c)|$ for every $c \in X \setminus M$, we see that
$$|f e_k - f e_l| = |f| \, |e_k - e_l|
 \leq r \, |g| \, |e_k - e_l| =
|g e_k - g e_l| \ .$$

So we have that
\begin{eqnarray*}
\lan  \pi(f e_k)  v- \pi(f e_l)  v , \pi(f e_k) v - \pi(f e_l) v \ran 
& = & \lan (\pi(f e_k) - (f e_l))^* (\pi(f e_k) - \pi(f e_l)) v ,
v \ran \\
& = & \lan  \pi(|f e_k - f e_l|^2) v , v \ran  \\
& \leq & r^2 \, \lan \pi(|g e_k - g e_l|^2) v , v \ran \\
& = &  r^2 \, \lan (\pi(g e_k) - \pi(g e_l))^* (\pi(g
e_k) - \pi(g e_l)) v , v \ran \\
& = & r^2 \, \lan \pi(g e_k) v - \pi(g e_l) v , \pi(g e_k) v - \pi(g e_l) v \ran
\end{eqnarray*}
which implies that $\|\pi(f e_k) v - \pi(f e_l) v\| \leq r \,
\|\pi(g e_k) v - \pi(g e_l) v\|$.

\vspace{1mm}

Because $v$ belongs to $D(\pi(g))$, we know that $(\pi(g e_k)v)_{k \in K}$ is convergent in
$E$. Therefore, the previous estimation guarantees that $(\pi(f e_k) v)_{k \in K}$ is also
convergent in $E$, which in turn implies that $v$ belongs to $D(\pi(f))$.
\end{demo}

\begin{lemma} 
Consider $f,g \in \text{C}(X)$ such that there exists a positive number $r$ 
such that $|f(x)| \leq r \, |g(x)|$ for every $x \in X
$. Then $\lan \pi(f) v , \pi(f) v \ran \leq 
r^2 \, \lan \pi(g) v , \pi(g) v \ran$ for every $v \in D(\pi(g))$.
\end{lemma}
\begin{demo}
Remember from the previous lemma that $D(\pi(g)) \subseteq D(\pi(f))$.

Take an approximate unit $(e_k)_{k \in K}$ for $\text{C}_0(X)$ in $\text{K}(X)$.

Then $(\pi(f e_k) v)_{k \in K}$ converges to $\pi(f) v$ and
$(\pi(g e_k) v)_{k \in K}$ converges to $\pi(g) v$.

We have for every $k \in K$ that
$$\lan \pi(f e_k) v , \pi(f e_k) v \ran = \lan \pi(|f e_k|^2) v , v \ran
\leq r^2 \, \, \lan \pi(|g e_k|^2) v , v \ran = r^2 \,\, \lan \pi(g e_k) v , \pi(g e_k) v \ran$$
So the lemma follows now immediately.
\end{demo}

\begin{lemma} \label{art1.lem4}
Consider $f \in \text{C}(X)$, $M$ a compact subset of $X$ and $r$ a positive number such that
$|f(x)| > r$ for every $x \in M$.
Then there exist an element $g \in \text{K}(X)$ such that $g=1$ on $M$, $0 \leq g \leq 1$ and
$r^2 \, \lan \pi(g) v , \pi(g) v \ran \leq \lan \pi(f) v , \pi(f) v \ran$ for every $v \in  D(\pi(f))$.
\end{lemma}
\begin{demo} If $M = \emptyset$ or $r=0$, the lemma is trivially true. So suppose that $M \neq \emptyset$ and $r \neq 0$.

Define for every $x \in M$ the set $O_x = \{ \, y \in X \mid |f(y)| > r \, \}$, then $O_x$ is an open neigbourhood 
of $x$. By the compactness of $M$, there exist $x_1, \ldots\! , x_n \in M$ such that $M \subseteq O_{x_1} \cup \ldots \cup O_{x_n}$.

Put $U = O_{x_1} \cup \ldots \cup O_{x_n}$. Then $U$ is an open subset of $X$ such that
$M \subseteq U$ and $|f(y)| \geq  r$ for every $y \in U$.

By the lemma of Urysohn, there exist an element $g \in \text{K}(X)$ such that $g=1$ on $M$, $g = 0$ on $X \setminus U$ 
and $0 \leq g \leq 1$. Then it is clear that $|f| \geq r \, g$. The previous lemma now guarantees that $r^2 \, \lan \pi(g) v , \pi(g) v \ran
\leq \lan \pi(f) v , \pi(f) v \ran$ for every $v \in  D(\pi(f))$.
\end{demo}

\bigskip\medskip

We end this section with a result of which the proof is easy in the Hilbert space case because then $\overline{\bld T^*} = 
(\ker T)^\perp$. The bounded case has already been proven in proposition 3.7 of \cite{Lan}. We will use in fact this result in 
the proof of the next proposition.

\begin{proposition} \label{art1.prop1}
Consider Hilbert \cst-modules $E$,$F$ over a \cst-algebra $A$ and an element $T \in \cR(E,F)$.
Then we  have  that \ $\overline{\bld T^* T} = \overline{\bld T^*}$ and \ $\ker T^* T = \ker T$
\end{proposition}
\begin{demo} Take a bounded net $(e_i)_{i \in I}$ in $\text{K}(\C)$ such that 
$(e_i)_{i \in I}$ converges strictly to 1.

Take $j \in J$. Then $e_j(T^* T) E \subseteq D(T^* T)$ by lemma \ref{art1.lem1}, so we have certainly that  
$e_j(T^* T) E \subseteq D(T)$. Therefore is $e_j(T^* T)$ a right multiplier of $T$. So we have the element 
$T \, e_j(T^* T) \in \cL(E,F)$.

Hence, proposition 3.7 of \cite{Lan} implies that
$$\overline{\bld (T \, e_j(T^* T))^*} = \overline{\bld (T \, e_j(T^* T))^* (T \, e_j(T^* T))} $$

It is clear that  $e_j(T^* T) \, T^*  \subseteq (T \, e_j(T^* T))^*$.

\vspace{1mm}

This implies also that 
$$ e_j(T^* T) \, T^* T \, e_j(T^*T) \subseteq (T \, e_j(T^* T))^* (
T \, e_j(T^* T))$$

By result \ref{art1.res2}, we know that $T^* T \, e_j^2(T^* T)$ and $e_j(T^* T) \, T^* T \, e_j(T^*T)$ are equal 
on the intersection of their domains. But $e_j(T^* T)$ and $e_j^2(T^* T)$ are right multipliers of $T^* T$. So 
$$ T^* T \, e_j^2(T^* T) = e_j(T^* T) \, T^* T \, e_j(T^*T) \in \cL(E)$$

Hence, $$T^* T \, e_j^2(T^* T) =  (T \, e_j(T^* T))^* (
T \, e_j(T^* T))$$

Consequently, 
\begin{eqnarray*}
\bld(e_j(T^* T) \, T^*) & \subseteq & \overline{\bld (T \, e_j(T^* T))^*}
= \overline{\bld (T \, e_j(T^* T))^* (T \, e_j(T^* T))} \\
& =  & \overline{\bld (T^* T \, e_j^2(T^* T))}
\subseteq \overline{\bld T^* T} 
\end{eqnarray*}

\medskip

Choose $v \in D(T^*)$. It is clear that the net $( e_i(T^* T) \, T^*(v) )_{i \in I}$ converges to $T^*(v)$. 
By the equation above, we know that $e_i(T^* T) \, T^*(v) \in \overline{\bld T^* T}$ for every $i \in I$. 
Hence, $T^*(v)$ belongs to $\overline{\bld T^*T}$.

So we see that $\overline{\bld T^*} \subseteq \overline{\bld T^* T}$. The other inclusion is trivially true. 

\vspace{1mm}

The equality involving the kernel is straightforward to prove.
\end{demo}

\bigskip

\section{Invertible regular operators} \label{art2}

In the first part of this section, we introduce the proper notion of invertibility of regular operators and 
investigate its basic properties (see page 104 of \cite{Lan}). In the second part, we look at the proper 
notion of bounded invertibility. 

\begin{definition} \label{art2.def1}
Consider Hilbert \cst-modules $E$,$F$ over a \cst-algebra $A$ and let $T$ be an element in $\cR(E,F)$. 
Then we call $T$ invertible (in $\cR(E,F)$) $\Leftrightarrow$ $T$ is injective and $T^{-1}$ belongs to $\cR(F,E)$.
\end{definition}

If $T$ is invertible, then $T^{-1}$ is densely defined so $T$ has dense range. It is also immediately clear 
that $T^{-1}$ is invertible if $T$ is invertible.

\medskip

As usual, we have the following results.

\begin{lemma} \label{art2.lem1}
Consider Hilbert \cst-modules $E$,$F$ over a \cst-algebra $A$ and let $T$ be an element in $\cR(E,F)$. 
Then we have the following properties :
\begin{itemize}
\item If $T$ has dense range, then $T^*$ is injective.
\item If $T$ has dense range and is injective, then $(T^{-1})^* = (T^*)^{-1}$.
\end{itemize}
\end{lemma}
\begin{demo}
\begin{itemize}
\item Suppose that $T$ has dense range. Take $v \in D(T^*)$ such that $T^* v = 0$. Then we have for every 
$w \in D(T)$ that $\lan v , T w \ran = \lan T^* v , w \ran = 0$. Because $T$ has dense range, this implies that $v = 0$. 
\item Suppose that $T$ has dense range and that $T$ is injective. We know already that
$T^*$ is injective.
\begin{enumerate}
\item Choose $v \in D(T^*)$. Then we have for every $w \in D(T)$ that
$$\lan T^{-1} (T w) , T^* v \ran = \lan w , T^* v \ran = \lan T w , v \ran \ . $$
So $T^* v$ belongs to $D((T^{-1})^*)$ and $(T^{-1})^*(T^* v) = v$

This implies that $(T^*)^{-1} \subseteq (T^{-1})^*$.
\item Choose $v \in D((T^{-1})^*)$. Then we have for every $w \in D(T)$ that
$$\lan (T^{-1})^*(v) , T w \ran = \lan v , T^{-1}(T w) \ran = \lan v , w \ran \ .$$
So $(T^{-1})^*(v)$ belongs to $D(T^*)$ and
$T^*((T^{-1})^*(v))= v$.

This implies that $(T^{-1})^* \subseteq (T^*)^{-1}$.
\end{enumerate}
So we get that $(T^{-1})^* = (T^*)^{-1}$. 
\end{itemize}
\end{demo}

\medskip

\begin{result}
Consider Hilbert \cst-modules $E$,$F$ over a \cst-algebra $A$ and let $T$ be an element in $\cR(E,F)$. 
Then we have the following properties.
\begin{itemize}
\item $T$ is invertible $\Leftrightarrow$ $T^*$ is invertible
\item If $T$ is invertible, then $(T^*)^{-1} = (T^{-1})^*$.
\end{itemize}
\end{result}
\begin{demo} Suppose that $T$ is invertible.

Because $T$ has dense range and is injective, the previous lemma implies that $T^*$ is injective
and $(T^*)^{-1} = (T^{-1})^*$. Hence, the regularity of $T^{-1}$ implies that
$(T^*)^{-1}$ is regular.

The other implication in the first statement follows by symmetry.
\end{demo}

\medskip

Another useful characterization of invertibility is given in the following proposition.

\begin{proposition} \label{art2.prop1}
Consider Hilbert \cst-modules $E$,$F$ over a \cst-algebra $A$ and let $T$ be an 
element in $\cR(E,F)$. Then $T$ is invertible $\Leftrightarrow$ $T$ and $T^*$ have dense range.
\end{proposition}
\begin{demo}
\begin{trivlist}
\item[\,\,$\Rightarrow$] This follows immediately from the previous proposition and the remark after 
definition \ref{art2.def1}.
\item[\,\,$\Leftarrow$] Because $T^*$ has dense range, $T$ is injective by lemma \ref{art2.lem1}.

Because $T$ has dense range, we have that $T^{-1}$ is a densely defined closed A-linear operator from within $F$ into $E$.

By lemma \ref{art2.lem1}, we have moreover that $(T^{-1})^* = (T^*)^{-1}$. Because $T^*$ has dense range, 
this implies that $(T^{-1})^*$ is densely defined.

Denote the flip map from $E \oplus F$ to $F \oplus E$ by $U$. Then theorem 9.3 of \cite{Lan} implies that 
$$G(T^{-1}) + G(T^{-1})^\perp = U(G(T)) + U(G(T))^\perp = U(G(T)+G(T)^\perp) = U(E \oplus F) = F \oplus E \ . $$
Hence, proposition 9.5 of \cite{Lan} implies that $T^{-1}$ is regular.
\end{trivlist}
\end{demo}

\medskip

Because the image of $z_T$ and $T$ are the same, this proposition implies immediately the following result.

\begin{corollary}
Consider Hilbert \cst-modules $E$,$F$ over a \cst-algebra $A$ and let $T$ be an element in $\cR(E,F)$. 
Then $T$ is invertible $\Leftrightarrow$ $z_T$ is invertible (in $\cR(E,F)$!).
\end{corollary}

\medskip

Another immediate corollary is the following one.

\begin{corollary} \label{art2.cor1}
Consider a Hilbert \cst-module $E$ over a \cst-algebra $A$ and let $T$ be a normal element in $\cR(E)$. 
Then $T$ is invertible $\Leftrightarrow$ $T$ has dense range.
\end{corollary}
\begin{demo}
Suppose that $T$ has dense range. Then proposition \ref{art1.prop1} implies that $T T^*$ has dense range. 
So $T^* T$ has also dense range. Consequently has $T^*$ also dense range.
\end{demo}

\medskip

This implies of course also that strictly positive elements are just the invertible positive elements.

\bigskip

Using proposition \ref{art1.prop1}, we get immediately the next result.

\begin{proposition}
Consider Hilbert \cst-modules $E$,$F$ over a \cst-algebra $A$ and let $T$ be an elements in $\cR(E,F)$. 
Then $T$ is invertible $\Leftrightarrow$ $T^* T$ and $T T^*$ are invertible.
\end{proposition}

\medskip

Invertibility is transferred by non-degenerate $^*$-homomorphisms : 

\begin{proposition} \label{art2.prop2}
Consider a Hilbert \cst-module $E$ over a \cst-algebra $A$. Let $B$ be a \cst-algebra and $\pi$ a non-degenerate 
$^*$-homomorphism from $B$ into $\cL(E)$. Consider an invertible element $T$ affiliated with $B$. Then $\pi(T)$ 
is invertible and $\pi(T)^{-1} = \pi(T^{-1})$.
\end{proposition}
\begin{demo}
We have for every $b \in D(T)$ and $v \in E$ that $\pi(b) v$ belongs to $D(\pi(T))$ and
$\pi(T)(\pi(b) v) = \pi(T(b)) v$. Hence, the density of the range of $T$ and the non-degeneracy of $\pi$ imply 
that $\pi(T)$ has dense range. We get in a similar way that $\pi(T)^* = \pi(T^*)$ has dense range.

So we see that $\pi(T)$ is invertible.

\medskip

Choose $b \in D(T^{-1})$ and $v \in E$. Then $\pi(b) v$ belongs to $D(\pi(T^{-1}))$ and
$\pi(T^{-1})(\pi(b) v) = \pi(T^{-1}(b)) v$.

Because $T^{-1}(b)$ belongs to $D(T)$, we get that
$\pi(T^{-1}(b)) v$ belongs to $D(\pi(T))$ and that \newline 
$\pi(T)\bigl(\pi(T^{-1}(b))v\bigr) = \pi(T(T^{-1}(b))) v = \pi(b) v $.

This implies that $\pi(b) v$ belongs to $D(\pi(T)^{-1})$ and
$\pi(T)^{-1}(\pi(b) v) = \pi(T^{-1}(b)) v = \pi(T^{-1})(\pi(b) v)$.

Consequently, the closedness of $\pi(T)^{-1}$ and the fact that the set $\lan \, \pi(b) v \mid b \in D(T^{-1}), 
v \in E \, \ran$ is a core for $\pi(T^{-1})$ imply that
$\pi(T^{-1}) \subseteq \pi(T)^{-1}$.

\medskip

We have in a similar way that $\pi((T^*)^{-1}) \subseteq \pi(T^*)^{-1}$, which gives us that
$(\pi(T)^{-1})^* \subseteq (\pi(T)^{-1})^*$. Hence $\pi(T)^{-1} \subseteq \pi(T^{-1})$ by taking the adjoint of the last inclusion.

So we arrive at the conclusion that $\pi(T)^{-1} = \pi(T^{-1})$.
\end{demo}

Even if $\pi$ is injective, the converse of this implication is not true :
Define the non-degenerate $^*$-homomorphism $\pi$ from $\text{C}([0,1])$ into $\text{C}_0(]0,1[)$ such that 
$\pi(f) = f_{]0,1[}$ for every $f \in \text{C}([0,1])$. Then $\pi$ is injective but you can find easily a function 
$g \in \text{C}([0,1])$ such that $\pi(g)$ is invertible but $g$ isn't (cfr. the next result).

\bigskip

In the commutative case, this form of invertibility corresponds to the usual one.

\begin{proposition}
Consider a locally compact space $X$ and let $f$ be an element in $\text{C}(X)$. Then we have the following properties :
\begin{itemize}
\item We have that $f$ is invertible $\Leftrightarrow$ $f(x) \neq 0$ for every $x \in X$.
\item If $f$ is invertible, then $f^{-1}(x) = \frac{1}{f(x)}$ for every $x \in X$.
\end{itemize}
\end{proposition}
\begin{demo}
\begin{itemize}
\item Suppose there exists $x_0 \in X$ such that $f(x_0) = 0$. Then it is clear that
$g(x_0)=0$ for every $g \in \overline{f \, \text{C}_0(X)}$. This implies that $f \, \text{C}_0(X)$ is not dense in 
$\text{C}_0(X)$. So $f$ is not invertible.
\item Suppose that $f(x) \neq 0$ for every $x \in X$.

Take $h \in \text{K}(X)$ and define $g \in \text{K}(X)$ such that $g(x) = \frac{h(x)}{f(x)}$ for every $x \in X$. 
Then $g$ belongs to $D(f)$ and $f(g) = h$. 
This implies that $f$ has dense range. \newline 
We get in a similar way that $\overline{f}$ has dense range. So $f$ is invertible.

Define the function $k \in \text{C}(X)$ such that $k(x) = \frac{1}{f(x)}$ for every $x \in X$. Then it is easy to check 
that $f \, k = \io_{D(k)}$ and that $k \, f = \io_{D(f)}$. So $f^{-1} = k$.
\end{itemize}
\end{demo}

\medskip

The following corollary will be useful in the next section.

\begin{corollary} \label{art2.cor2}
Consider a locally compact space $X$ and let $f$ be an element in $\text{C}(X)$. Let $c$ be a complex number. 
Then $c \notin f(X)$ $\Leftrightarrow$ $f - c \, 1$ is invertible.
\end{corollary}

\bigskip\medskip

We will now consider a stronger form of invertibility which is natural to look at in connection with the 
spectrum of an element.

\begin{definition}
Consider two Hilbert \cst-modules $E$,$F$ over a \cst-algebra $A$ and let $T$ be 
an element in $\cR(E,F)$. Then we call $T$ adjointable invertible $\Leftrightarrow$ There exists an element 
$S \in \cL(F,E)$ such that $S T \subseteq T S = 1$ 
\end{definition}

\begin{corollary}
Consider two Hilbert \cst-modules $E$ and $F$ over a \cst-algebra $A$ and let $T$ 
be an element in $\cR(E,F)$. Then $T$ is adjointable invertible $\Leftrightarrow$ $T$ is invertible and 
$T^{-1}$ belongs to $\cL(F,E)$.
\end{corollary}

\begin{remark} \rm
Let $E$ be a Hilbert \cst-module over a \cst-algebra $A$ and let $T$ an element in $\cL(E)$. Then it is clear 
that $T$ is adjointable invertible if and only if $T$ is invertible in the \cst-algebra $\cL(E)$.
\end{remark}

\medskip

A useful characterization for adjointable invertibility is the following generalization of a well known Hilbert result.

\begin{result} \label{art2.res1}
Consider two Hilbert \cst-modules $E$ and $F$ over a \cst-algebra $A$ and let $T$ 
be an element in $\cR(E,F)$.  Then $T$ is adjointable invertible $\Leftrightarrow$ $\text{Ran}\,T = F$ and 
$\text{Ran}\,T^*$ is dense in $E$ $\Leftrightarrow$
$\text{Ran}\,T = F$ and $\text{Ran}\,T^* = E$
\end{result}
\begin{demo}
\begin{itemize}
\item Suppose that $T$ is adjointable invertible. It is clear that $\bld T =F$.

Now  $T$ is invertible and $T^{-1}$ belongs to $\cL(F,E)$. We know that $T^*$ is invertible and 
$(T^*)^{-1} = (T^{-1})^* \in \cL(E,F)$, so $\text{Ran}\,T^* = F$.
\item Suppose that $\text{Ran}\,T = F$ and $\text{Ran}\,T^*$ is dense in $E$.
By \ref{art2.prop1}, we know that $T$ is invertible. Because $\text{Ran}\,T = F$, we get that $D(T^{-1}) = F$.
Combining this with the closedness of $T^{-1}$, the closed graph theorem implies that
$T^{-1}$ is bounded.

We know also that $T^*$ is invertible and that $(T^{-1})^* = (T^*)^{-1}$. Because $\text{Ran}\,T^*$ is dense in $E$, 
this implies that $(T^{-1})^*$ is densely defined. 

Because  $T^{-1}$ is bounded, we have also that $(T^{-1})^*$ is bounded.
So combining this with the closedness of $(T^{-1})^*$, we must have that $D((T^{-1})^*)$ closed.
Hence, $D((T^{-1})^*) = E$. 

From this all, we conclude that $T^{-1}$ belongs to $\cL(F,E)$.
\end{itemize}
\end{demo}

\medskip

This implies immediately the following result.

\begin{corollary}
Consider two Hilbert \cst-modules $E$ and $F$ over a \cst-algebra $A$ and let $T$ 
be an element in $\cR(E,F)$. Then $T$ is adjointable invertible $\Leftrightarrow$ $T^*$ is adjointable invertible.
\end{corollary}

\medskip

Combining result \ref{art2.res1} with corollary \ref{art2.cor1} , we get also the following one.

\begin{corollary}
Consider a Hilbert \cst-module $E$ over a \cst-algebra $A$ and let $T$ 
be a normal element in $\cR(E)$. Then $T$ is adjointable invertible $\Leftrightarrow$ 
$\text{Ran}\,T = E$.
\end{corollary}

\bigskip

We can again connect the adjoint invertibility of $T$ to the invertibility of its $z$-transform.
The proof of this result follows immediately from the fact that $\bld T = \bld z_T$ and
$\bld T^* = \bld (z_T)^*$ and result \ref{art2.res1}.

\begin{result} \label{art2.res2}
Consider two Hilbert \cst-modules $E$ and $F$ over a \cst-algebra $A$ and let $T$ 
be an element in $\cR(E,F)$. Then $T$ is adjointable invertible $\Leftrightarrow$ 
$z_T$ is adjointable invertible.
\end{result}

\medskip

This characterization of adjoint invertibility makes the proof of the following result very easy. However, 
a little bit of care has to be taken for the second equivalence. You need the fact that invertibility for 
an element in $M(\pi(B))$ is equivalent to the invertibility of this element in $\cL(E)$.

\begin{result} \label{art2.res3}
Consider a Hilbert \cst-module $E$ over a \cst-algebra $A$. Let $B$ be a \cst-algebra and $\pi$ a 
non-degenerate $^*$-homomorphism from $B$ into $\cL(E)$. Let $T$ be a regular operator in $E$. Then we have the following properties.
\begin{itemize}
\item If $T$ is adjointable invertible, then $\pi(T)$ is adjointable invertible.
\item Suppose moreover that $\pi$ is injective. 

Then $T$ is adjointable invertible $\Leftrightarrow$ $\pi(T)$ is adjointable invertible.
\end{itemize}
\end{result}

\bigskip\medskip

Another application of result \ref{art2.res2} is the proof of the fact that the definition of the spectrum of 
a normal regular operator given by Woronowicz coincides with the usual definition in the Hilbert space case. First we need an easy lemma.

\begin{lemma}
Consider a Hilbert \cst-module $E$ over a \cst-algebra $A$ and let $T$ be a normal element in $\cR(E)$. 
Then we have for every $c \in \C$ that $\si(T - c \, 1) = \si(T) - c$.
\end{lemma}
\begin{demo}
Define $f \in \text{C}(\C)$ such that $f(\lambda) = \lambda - c$ for $\lambda \in \C$. Then 
corollary \ref{art1.cor1} implies that $f(T) = T - c \, 1$. Hence, $$\si(T - c \, 1) = \overline{f(\si(T))} = 
\overline{\si(T) - c} = \si(T) - c \ .$$
\end{demo}

\medskip

\begin{proposition} \label{art2.prop3}
Consider a Hilbert \cst-module $E$ over a \cst-algebra $A$ and let $T$ be a normal element in $\cR(E)$.
Let $c$ be a complex number. Then $c \notin \si(T)$ $\Leftrightarrow$ $T - c \, 1$ is adjointable invertible
\end{proposition}
\begin{demo} Using the previous lemma, result \ref{prel3.res1} and result \ref{art2.res2}, we get that
\begin{eqnarray*}
c \notin \si(T) & \Leftrightarrow &
0 \notin \si(T-c\,1) \Leftrightarrow 0 \notin \si(z_{T-c\,1})
\Leftrightarrow z_{T-c\,1} \text{ is invertible in } \cL(E) \\
& \Leftrightarrow & T - c\, 1 \text{ is adjointable invertible }
\end{eqnarray*}
\end{demo}

\bigskip

This suggests the following natural and consistent definition of the spectrum of a general regular operator.

\begin{definition}
Consider a Hilbert \cst-module $E$ over a \cst-algebra $A$ and let $T$ be a normal element in $\cR(E)$. 
We define the resolvent of $T$ as the set $\rho(T) = \{ \, c \in \C \mid T - c \, 1 \text{ is adjointable invertible } \,\}$
and we define the spectrum of $T$ as the set $\si(T) = \C \setminus \rho(T)$.
\end{definition}

\medskip

The usual Hilbert space techniques allow us to prove the next two results.

\begin{result}
Consider a Hilbert \cst-module $E$ over a \cst-algebra $A$ and let $T$ be a normal element in $\cR(E)$. 
Then $\si(T)$ is closed and the mapping
$\rho(T) \rightarrow \cL(E) : c \mapsto (T-c)^{-1}$ is continuous.
\end{result}
\begin{demo}
Choose $d \in \rho(T)$. Define $M = \|(T - d \, 1)^{-1}\| \in \R^+$.

\vspace{1mm}

We have for every $c \in \C$ that $(T - d \, 1)^{-1} \, E = D(T - d \, 1) = D(T - c \, 1)$ so result \ref{prel2.res1} 
implies that $(T - c \, 1) \, (T - d \, 1)^{-1}$ belongs to $\cL(E)$.
Define the function $f$ from $\C$ into $\cL(E)$ such that $f(c) =$ $(T - c\,1) \, (T - d \, 1)^{-1}$ for every $c \in \C$. Then $f(d) = 1$.

We have for every $c_1, c_2 \in \C$ that
\begin{eqnarray*}
\|f(c_1) - f(c_2) \| & = & \| (T - c_1 \, 1) \, (T - d \, 1)^{-1} - (T - c_2 \, 1) \, (T - d \, 1)^{-1} \| \\
& = & \| (c_2 - c_1) \, (T - d \, 1)^{-1}\| = |c_1 - c_2| \, M
\end{eqnarray*}
so the function  $f$ is continuous.

Choose $c \in \C$ such that $|c - d| \leq \frac{1}{M+1}$. Then the previous inequality implies that $\|f(c) - 1\| < 1$ which 
in turn implies that $f(c)$ is invertible in $\cL(E)$.

It is clear that $f(c) (T - d \, 1) = T - c \, 1$. Hence we get that $T - c \, 1$ is injective
and that
$$(T - c \, 1)^{-1} = (T - d \, 1)^{-1} \, f(c)^{-1}$$
So $(T - c \, 1)^{-1}$ belongs to $\cL(E)$. By definition, 
we have that $c$ belongs to $\rho(T)$.

\medskip

The above equation and the continuity of the invertibility operation in $\cL(E)$ imply also immediately that the function 
$\rho(T) \rightarrow \cL(E) : c \mapsto (T-c)^{-1}$ is continuous.
\end{demo}

\medskip

\begin{result}
Consider a Hilbert \cst-module $E$ over a \cst-algebra $A$ and let $T$ be a normal element in $\cR(E)$. Then the mapping 
$\rho(T) \rightarrow \cL(E) : c \mapsto (T-c)^{-1}$ is analytic.
\end{result}

As usual, this follows from the fact that
$(T - c \, 1)^{-1} - (T - d \, 1)^{-1} = (c - d) \, (T - c \, 1)^{-1}\, (T - d \, 1)^{-1}$
for $c,d \in \rho(T)$.

\bigskip

Result \ref{art2.res3} implies also the following result.

\begin{result}
Consider a Hilbert \cst-module $E$ over a \cst-algebra $A$. Let $B$ be a \cst-algebra and $\pi$ a non-degenerate 
$^*$-homomorphism from $B$ into $\cL(E)$. Let $T$ be a regular operator in $E$. Then we have the following properties.
\begin{itemize}
\item $\si(\pi(T)) \subseteq \si(T)$
\item If $\pi$ is injective, then $\si(\pi(T)) = \si(T)$.
\end{itemize}
\end{result}

\bigskip

\section{The functional calculus revisited}  \label{art3}

Consider a normal regular operator $T$ in a Hilbert \cst-module. Until now, we had a functional calculus for $T$ on 
the spectrum $\si(T)$. But there are cases where we want a functional calculus for $T$ on $\si(T) \setminus K$ where 
$K$ is a certain well-behaved finite subset of $\,\C$.

\medskip

A typical situation occurs when $T$ is strictly positive in which case we want to define any complex power of $T$. 
So we need in this case a functional calculus on $\si(T) \setminus \{0\}$.

\bigskip

We will always work with functional calculi on a special kind of subsets of $\,\C$. Therefore we will use the following terminology.

\begin{terminology}
Let $G$ be a subset of $\,\C$. Then we call $G$ almost closed if there exists  a finite subset $K$ of $\,\C$ such 
that $G \cup K$ is closed. 
\end{terminology}

It is clear that any almost closed subset is locally compact. By adding or substracting a finite number of points of an 
almost closed set, we keep an almost closed set.

\medskip

\begin{lemma}  \label{art3.lem1}
Let $G$ be an almost closed subset of $\,\C$ and $f \in \text{C}_0(G)$. Then there exists an element 
$g \in \text{C}_0(\C)$ such that $f \subseteq g$.
\end{lemma}
\begin{demo}
There exists a finite subset $K$ of $\,\C \setminus G$ such that $G \cup K$ is closed.
Define $C = \{\,h_{G \cup K} \mid h \in \text{C}_0(\C)\,\}$. Because $G \cup K$ is closed, $C$ will be a subset of 
$\text{C}_0(G \cup K)$. So $C$ is a sub \cst-algebra $\text{C}_0(G \cup K)$ (as the image under an obvious $^*$-homomorphism).
We have also clearly that
\begin{itemize}
\item There exists for every $c \in G \cup K$ an element $k \in C$ such that $k(c) \neq 0$.
\item There exists for every $c,d \in G \cup K$  with $c \neq d$, an element $k \in C$ such that $k(c) \neq k(d)$.
\end{itemize}
Hence, the Stone Weierstrass theorem implies that $C = \text{C}_0(G \cup K)$.

Now define $\tilde{f} \in \text{C}_0(G \cup K)$ such that $\tilde{f} = 0$ on $F$ and
$\tilde{f} = f$ on $G$ \ (the closedness of $F$ is here used to get the continuity of $\tilde{f}$ in points of $G$). 
Then there exists $g \in \text{C}_0(\C)$ such that
$g_{G \cup K} = \tilde{f}$. Hence, $f \subseteq g$.
\end{demo}

\medskip

\begin{result} \label{art3.res1}
Consider a Hilbert \cst-module over a \cst-algebra $A$ and an almost closed subset $G$ of $\,\C$ and $f \in \text{C}_0(G)$. 
Let $\pi$, $\th$ be non-degenerate $^*$-homomorphisms from $\text{C}_0(G)$ into $\cL(E)$. 

If $\pi(\io_G) = \th(\io_G)$, then $\pi = \th$.
\end{result}
\begin{demo}
By proposition \ref{prel3.prop1}, we have for every $g \in \text{C}_0(\C)$ that
$\pi(g_G) = \pi(g \!\circ\! \io_G) = \pi(g(\io_G)) = g(\pi(\io_G))$ and similarly,
$\th(g_G) = g(\pi(\io_G))$. So we see that $\pi(g_G) = \th(g_G)$.
Now the previous lemma implies that $\pi = \th$.
\end{demo}

\bigskip

We will soon prove the  key result to introduce these new functional calculi but first we need the following lemma.

\begin{lemma}
Consider a Hilbert \cst-module $E$ over a \cst-algebra $A$ and let $T$ be a normal element in $\cR(E)$. Consider 
$\lambda_1,\ldots\!,\lambda_n \in \C$ such that $T - \lambda_1\,1, \ldots\! , T-\lambda_n\,1$ have dense range.

Define the function $g$ from $\text{C}(\si(T))$ such that $g(c) = (c-\lambda_1) \ldots (c-\lambda_n)$ for every 
$c \in \si(T)$. Then $g(T)$ has dense range.
\end{lemma}
\begin{demo} Take an approximate unit $(e_k)_{k \in K}$ for $\text{C}_0(\si(T))$ in $\text{K}(\si(T))$.

Define for every $r \in \{1,\dots\!,n\}$ the function $g_r \in \text{C}(\si(T))$ such that
$g_r(c) = (c-\lambda_1) \ldots (c-\lambda_r)$ for every $c \in \si(T)$.
\begin{enumerate}
\item Because $g_1(T) = T - \lambda_1 \, 1$, we have by assumption that $g_1(T)$ has dense range.
\item Choose $r \in \{1,\ldots\!,n-1\}$ and suppose that $g_r(T)$ has dense range.

Define the function $h \in \text{C}(\si(T))$ such that $h(c) = c - \lambda_{r+1}$ for every $c \in \si(T)$. \newline 
Then $h(T) = T-\lambda_{r+1} \, 1$, so $h(T)$ has dense range.

\medskip

Take $k \in K$. By lemma \ref{art1.lem1}, we know that $e_k(T) E \subseteq D(g_{r+1}(T))$.
We have moreover that
\begin{eqnarray*}
\overline{g_{r+1}(T)(e_k(T) E)} & = & 
\overline{(g_{r+1} \, e_k)(T) E} = \overline{(g_{r+1} \, e_k)(T) D(h(T))} \\
& = & \overline{(g_r \, h \, e_k)(T) D(h(T))}
= \overline{(g_r \, e_k)(T)\bigl(h(T) D(h(T))\bigr)} \ . 
\end{eqnarray*}
So, using the fact that $h(T)$ has dense range, we see that
\begin{eqnarray*}
\overline{g_{r+1}(T)(e_k(T) E)} & = & \overline{(g_r \,e_k)(T) E}
\supseteq (g_r e_k)(T) E  \\
& \supseteq & \ (g_r \, e_k)(T) D(g_r(T)) = e_k(T)\bigl(g_r(T) D(g_r(T))\bigr) \ . 
\end{eqnarray*}
Hence, the assumed densitity of the range of $g_r(T)$ implies that
$$e_k(T) E \subseteq \overline{g_{r+1}(T)(e_k(T) E)} \subseteq \overline{\bld g_{r+1}(T)} \ . $$

\medskip

Because $(e_k(T))_{k \in K}$ converges strongly to 1, this implies that 
$\overline{\bld g_{r+1}(T)} = E$. 
\end{enumerate}

By induction, we can now conclude that $g(T)=g_n(T)$ has dense range. 
\end{demo}

\medskip

Now we are ready to prove the main result of this section.

\begin{proposition} \label{art3.prop4}
Consider a Hilbert \cst-module $E$ over a \cst-algebra $A$ and let $T$ be a normal element in $\cR(E)$. Consider a 
finite subset $K$ of $\,\C$ such that $T - \lambda \, 1$ is invertible for every $\lambda \in K$. 

Then there exists a unique non-degenerate $^*$-homomorphism $\pi$ from $\text{C}_0(\si(T) \setminus K)$ into $\cL(E)$ 
such that $\pi(\io_{\si(T) \setminus K} ) = T$.
We have moreover for every $f \in \text{C}(\si(T))$ that
$f(T) = \pi(f_{\si(T) \setminus K})$.
\end{proposition}
\begin{demo} The unicity is guaranteed by result \ref{art3.res1} so let us turn to the existence.

There is nothing to prove if $G \cap K = \emptyset$ so suppose that
$G \cap K \neq \emptyset$. Then there exist different complex numbers $\lambda_1,\ldots\!,\lambda_n$ such that 
$G \cap K = \{\lambda_1,\ldots\!,\lambda_n\}$.
Then $G \setminus K = G \setminus \{\lambda_1,\ldots\!,\lambda_n\}$ and we have for every 
$i \in \{1,\ldots,n\}$ that $T - \lambda_i \, 1$ has dense range.

\medskip

For every $f \in \text{C}_0(\si(T) \setminus\{\lambda_1,\ldots\!,\lambda_n\})$, we define the function 
$\tilde{f} \in \text{C}_0(\si(T))$ such that 
$$\tilde{f}=f \text{ on } \si(T) \setminus\{\lambda_1,\ldots\!,\lambda_n\} \hspace{1.5cm} \text{ and } 
\hspace{1.5cm} \tilde{f}(\lambda_i) = 0 \text{ for every } i \in \{1,\ldots\!,n\} \ . $$

Then the mapping $\text{C}_0(\si(T) \setminus\{\lambda_1,\ldots\!,\lambda_n\}) \rightarrow
\text{C}_0(\si(T)) : f \mapsto \tilde{f}$ is clearly a $^*$-homomorphism.

Now define the mapping  $\pi$ from $\text{C}_0(\si(T) \setminus\{\lambda_1,\ldots\!,\lambda_n\})$ into 
$\cL(E)$ such that $\pi(f) = \tilde{f}(T)$ for every $f \in \text{C}_0(\si(T) \setminus \{\lambda_1,\ldots\!,\lambda_n\})$. 
Then $\pi$ is a $^*$-homomorphism.

\begin{itemize}
\item First we prove that $\pi$ is non-degenerate.

Define the function $g \in \text{C}(\si(T))$ such that $g(c) = (c-\lambda_1) \ldots (c-\lambda_n)$ for every $c \in \si(T)$. 
By the previous lemma, we know that $g(T)$ has dense range.

Take an approximate unit $(e_k)_{k \in K}$ for $\text{C}_0(\si(T))$ in $\text{K}(\si(T))$.

\medskip

Take $k \in K$ and define $h_k$ as the restriction of $g \, e_k$ to $\si(T) \setminus \{\lambda_1,\ldots\!,\lambda_n\}$. 

Because $g(\lambda_1) = \ldots = g(\lambda_n) = 0$, we have 
that $h_k$ belongs to $\text{C}_0(\si(T) \setminus \{\lambda_1,\ldots\!,\lambda_n\})$
and that $\widetilde{h_k} = g \, e_k$, so $\pi(h_k) = (g \, e_k)(T)$.

This implies that 
$$e_k(T)\bigl(g(T) D(g(T))\bigr) = (g \, e_k)(T) D(g(T)) = \pi(h_k) D(g(T)) \subseteq \pi(\text{C}_0(\si(T) \setminus 
\{\lambda_1,\ldots\!,\lambda_n\}))\,E $$
So the density of the range of $g(T)$ implies that
$e_k(T) E \subseteq \overline{\pi(\text{C}_0(\si(T) \setminus \{\lambda_1,\ldots\!,\lambda_n\}))\,E}$.

Because $(e_k(T))_{k \in K}$ converges strongly to 1, this implies that
$\pi(\text{C}_0(\si(T) \setminus \{\lambda_1,\ldots\!,\lambda_n\})) \, E$ is dense in $E$.

\item Define the mapping $\th$ from $\text{C}_0(\si(T))$ into $\text{C}_b(\si(T) \setminus \{\lambda_1,\ldots\!,\lambda_n\})$ 
such that we have for every $h \in \text{C}_0(\si(T))$ that
$\th(h)$ is the restriction of $h$ to $\si(T) \setminus \{\lambda_1,\ldots\!,\lambda_n\}$.

Then $\th$ is non-degenerate and we have for every $h \in \text{C}(\si(T))$ that
$\th(h)$ is the restriction of $h$ to $\si(T) \setminus \{\lambda_1,\ldots\!,\lambda_n\}$.

\medskip

Choose $g \in \text{C}_0(\si(T))$. 

Take an approximate unit $(e_k)_{k \in K}$ for $\text{C}_0(\si(T) \setminus \{\lambda_1,\ldots\!,\lambda_n\})$ in 
$\text{K}(\si(T) \setminus \{\lambda_1,\ldots\!,\lambda_n\})$.

Choose $k \in K$. Then  $\th(g) \, e_k$ belongs to $\text{K}(\si(T) \setminus \{\lambda_1,\ldots\!,\lambda_n\})$ 
and $(\th(g)\,e_k)\golf = g \, \widetilde{e_k}$.

So we get that
$$\pi(\th(g)) \, \pi(e_k) = \pi(\th(g) \, e_k) 
 = (g \, \widetilde{e_k})(T) = g(T) \, \widetilde{e_k}(T) = g(T) \pi(e_k) \ .$$
Because $(\pi(e_k))_{k \in K}$ converges strongly to 1,
we get that $\pi(\th(g)) = g(T) = \psi_T(g)$.

Hence, proposition \ref{prel1.prop2}  implies for every $f \in \text{C}(\si(T))$ that
$\pi(\th(f)) = (\pi \th)(f) = \psi_T(f) = f(T)$.

We find in particular that $\pi(\io_{\si(T) \setminus \{\lambda_1,\ldots\!,\lambda_n\}})
= \pi(\th(\io_{\si(T)})) = \io_{\si(T)}(T) = T$.
\end{itemize}
\end{demo}

\medskip 

The rest of this section is devoted to introduce the proper notations and prove quickly the usual results for 
functional calculi. The proofs of these results do not differ very much of the proofs of their corresponding results 
for the ordinary functional calculus.

\bigskip

\begin{terminology}
Consider a Hilbert \cst-module $E$ over a \cst-algebra $A$ and let $T$ be a normal element in $\cR(E)$. Let $f$ be a 
function from within $\,\C$ into $\,\C$. Then we call $f$ compatible with $T$
$\Leftrightarrow$ $f$ is continuous and there exists a finite subset $K$ of $\,\C$ such that $T - \lambda\,1$ is 
invertible for every $\lambda \in K$ and such that $f$ is defined on $\si(T) \setminus K$.
\end{terminology}

\medskip

We want to define $f(T)$ for such compatible functions $f$. First we have to resolve a small consistency problem. 
This will be done in the following lemma.

\begin{lemma} \label{art3.lem2}
Consider a Hilbert \cst-module $E$ over a \cst-algebra $A$ and let $T$ be a normal element in $\cR(E)$. 
Consider finite subsets $K$,$L$ of $\,\C$ such that $T - \lambda\,1$ is invertible for every $\lambda \in K \cup L$. 
\begin{itemize}
\item Let $\pi$ denote the unique non-degenerate $^*$-homomorphism $\pi$ from $\text{C}_0(\si(T) \setminus K)$ 
into $\cL(E)$ such that  $\pi(\io_{\si(T) \setminus K}) = T$.
\item  Let $\th$ denote the unique non-degenerate $^*$-homomorphism $\pi$ from $\text{C}_0(\si(T) \setminus L)$ into $\cL(E)$ such that
$\th(\io_{\si(T) \setminus L}) = T$.
\end{itemize}
Then we have for every continuous function $f$ from within $\,\C$ into $\,\C$ which is defined on both $\si(T) 
\setminus K$ and $\si(T) \setminus L$ that $\pi(f_{\si(T) \setminus K}) = \th(f_{\si(T) \setminus L})$.
\end{lemma}
\begin{demo}
Let $\eta$ denote the unique non-degenerate $^*$-homomorphism $\pi$ from $\text{C}_0(\si(T) \setminus (K \cup L))$ into 
$\cL(E)$ such that $\eta(\io_{\si(T) \setminus (K \cup L)}) = T$.

Define $\rho$ to be the $^*$-homomorphism from $\text{C}_0(\si(T) \setminus K)$ into 
$\text{C}_b(\si(T) \setminus (K \cup L))$ such that we have for every $g \in \text{C}_0(\si(T) \setminus K)$ that $\rho(g)$ is 
the restriction of $g$ to $\si(T) \setminus (K \cup L)$.
Then $\rho$ is non-degenerate and we have for every $g \in \text{C}(\si(T) \setminus K)$ that
$\rho(g)$ is the restriction of $g$ to $\si(T) \setminus (K \cup L)$.

So we have that $(\eta \rho)(\io_{\si(T) \setminus K}) = \eta(\rho(\io_{\si(T) \setminus K}))
= \eta(\io_{\si(T) \setminus (K \cup L)}) = T$, which implies that $\pi = \eta \rho$

This implies that $\pi(f_{\si(T) \setminus K})
= \eta(\rho(f_{\si(T) \setminus K})) = \eta(f_{\si(T) \setminus (K \cup L)})$.

\medskip

We get in a similar way that $\th(f_{\si(T) \setminus L}) = \eta(f_{\si(T) \setminus (K \cup L)})$. So the lemma follows.
\end{demo}

\medskip

Proposition \ref{art3.prop4} (or the previous lemma with $L = \emptyset$) implies moreover the following lemma.

\begin{lemma}
Consider a Hilbert \cst-module $E$ over a \cst-algebra $A$ and let $T$ be a normal element in $\cR(E)$. Consider a finite 
subset $K$ of  $\,\C$ such that $T - \lambda\,1$ is invertible for every $\lambda \in K$.

Let $\pi$ denote the unique non-degenerate $^*$-homomorphism $\pi$ from $\text{C}_0(\si(T) \setminus K)$ into $\cL(E)$ such 
that  $\pi(\io_{\si(T) \setminus K}) = T$. 
Then we have for every continuous function $f$ from within $\,\C$ into $\,\C$ which is defined on $\si(T)$ that 
$\pi(f_{\si(T) \setminus K}) = f(T)$.
\end{lemma}

\medskip

Now we are in a position to give the following consistent definition.

\begin{definition} \label{art3.def1}
Consider a Hilbert \cst-module $E$ over a \cst-algebra $A$ and let $T$ be a normal element in $\cR(E)$. Let $f$ be a function 
from within $\,\C$ into $\,\C$ which is compatible with $T$.
Then we define the normal element $f(T)$ in $\cR(E)$ as follows.

Take a finite subset $K$ of $\,\C$ such that $T - \lambda\,1$ is invertible for every $\lambda \in K$ and such that $f$ is 
defined on $\si(T) \setminus K$. 
Let $\pi$ denote the unique non-degenerate $^*$-homomorphism $\pi$ from $\text{C}_0(\si(T) \setminus K)$ into $\cL(E)$ such 
that $\pi(\io_{\si(T) \setminus K}) = T$.

Then we define $f(T) = \pi(f_{\si(T) \setminus K})$.
\end{definition}

If $f$ is a bounded function, it is clear that $f(T)$ belongs to $\cL(E)$.

\bigskip

\begin{terminology} \label{art3.def2}
Consider a Hilbert \cst-module $E$ over a \cst-algebra $A$ and let $T$ be a normal element in $\cR(E)$. Let $G$ be a subset 
of $\,\C$. Then we say that $G$ is compatible with $T$ if and only if  $G$ is almost closed and there exist a finite subset 
$K$ of $\,\C$ such that
\begin{enumerate}
\item The element $T - \lambda\,1$ is invertible for every $\lambda \in K$.
\item $\si(T) \setminus K \subseteq G$
\end{enumerate}
\end{terminology}

\medskip

\begin{remark} \rm
The following properties are immediately clear.
\begin{itemize}
\item The sets $\si(T)$ and $\,\C$ are always compatible with $T$.
\item An almost closed subset of $\,\C$ which contains a subset of $\,\C$ which is compatible with 
$T$ is itself compatible with $T$.
\item If $T$ is invertible, then $\,\C \setminus \{0\}$ and $\si(T) \setminus \{0\}$ are compatible with $T$.
\item IF $T$ is strictly positive, then $\R^+_0$ is compatible with $T$.
\end{itemize}
\end{remark}

\medskip

If $G$ is compatible with $T$, it is also clear that any element $f \in \text{C}(G)$ is compatible with $T$ and so we can look 
at the normal element $f(T)$ in $\cR(E)$.

\begin{proposition}
Consider a Hilbert \cst-module $E$ over a \cst-algebra $A$ and let $T$ be a normal element in $\cR(E)$. Let $G$ be a subset of 
$\,\C$ which is compatible with $T$. 

Define the mapping $\pi$ from $\text{C}_0(G)$ into $\cL(E)$ such that
$\pi(f) = f(T)$ for every $f \in \text{C}_0(G)$. 

Then $\pi$ is the unique non-degenerate $^*$-homomorphism from $\text{C}_0(G)$
into $\cL(E)$ such that $\pi(\io_G) = T$. We call $\pi$ the functional calculus of $T$ on $G$. 
We have moreover that $\pi(f) = f(T)$ for every $f \in \text{C}(G)$.
\end{proposition}
\begin{demo}
Take a finite subset $K$ of $\,\C$ such that
\begin{enumerate}
\item The element $T - \lambda\,1$ is invertible for every $\lambda \in K$.
\item $\si(T) \setminus K \subseteq G$
\end{enumerate}
Let $\th$ denote the unique non-degenerate $^*$-homomorphism from $\text{C}_0(\si(T) \setminus K)$ into $\cL(E)$ such that 
$\th(\io_{\si(T) \setminus K}) = T$.

Define the non-degenerate $^*$-homomomorphism $\eta$ from $\text{C}_0(G)$ into 
$\text{C}_b(\si(T) \setminus K)$ such that we have for every $f \in \text{C}_0(G)$ that $\eta(f)$ is equal to the restriction 
of $f$ to $\si(T) \setminus K$. Then we have for every $f \in \text{C}(G)$ that $\eta(f)$ is equal to the restriction of $f$ to $\si(T) \setminus K$.

Then we have by definition for every $f \in \text{C}_0(G)$ that $(\th \eta)(f) = \th(\eta(f))
= f(T) = \pi(f)$ so $\th \eta = \pi$.

This implies for every $f \in \text{C}(G)$ that $\pi(f) = \th(\eta(f)) = f(T)$. We have in particular that $\pi(\io_G) = 
\th(\eta(\io_G)) = \th(\io_{\si(T) \setminus K}) = T$.
\end{demo}

\medskip

Looking at definition \ref{art3.def1} and terminology \ref{art3.def2}, we get immediately the following result.

\begin{result} \label{art3.res3}
Consider a Hilbert \cst-module $E$ over a \cst-algebra $A$ and let $T$ be a normal element in $\cR(E)$. Let $G$,$H$ subsets 
of $\,\C$ which are compatible with $T$ and such that $G \subseteq H$. Then we have for every $f \in \text{C}(H)$ that $f(T) = (f_G)(T)$.
\end{result}

\medskip

We have of course the following injectivity property.

\begin{result} \label{art3.res2}
Consider a Hilbert \cst-module $E$ over a \cst-algebra $A$ and let $T$ be a normal element in $\cR(E)$. Let $G$ be a subset 
of $\,\C$ which is compatible with $T$ and such that $G \subseteq
\si(T)$. Then the functional calculus of $T$ on $G$ is injective.
\end{result}
\begin{demo} 
Call $\pi$ the functional calculus of $T$ on $G$. Take $f \in \text{C}_0(G)$ such that
$\pi(f) = 0$. Define $g \in \text{C}_0(\si(T))$ such that $g = f$ on $G$ and $g = 0$ on $\si(T) \setminus G$. Then 
$\psi_T(g) = \pi(f) =0$, so the injectivity of $\psi_T$ implies that
$g = 0$, so $f = 0$.
\end{demo}

\medskip

Every compatible subset of $\,\C$ gives rise to a compatible subset of $\,\C$ of the kind above.

\begin{result}
Consider a Hilbert \cst-module $E$ over a \cst-algebra $A$ and let $T$ be a normal element in $\cR(E)$. Let $G$ be a subset 
of $\,\C$ which is compatible with $T$. Then $G \cap \si(T)$ is compatible with $T$.
\end{result}
\begin{demo}
Because $G$ is compatible with $T$, there exists a finite subset $K$ of $\,\C$ such that
\begin{enumerate}
\item The element $T - \lambda\,1$ is invertible for every $\lambda \in K$.
\item $\si(T) \setminus K \subseteq G$
\end{enumerate}
Then it is clear that $G \cap \si(T) = \si(T) \setminus K$, so $G \cap \si(T)$ is compatible with $T$.
\end{demo}

\medskip

\begin{result}
Consider a Hilbert \cst-module $E$ over a \cst-algebra $A$ and let $T$ be a normal element  in $\cR(E)$. Let $G$ be a subset of 
$\,\C$ which is compatible with $T$ and $f$ an element in 
$\text{C}(G)$. 

Then  $\si(f(T)) = \overline{f(G \cap \si(T))} \subseteq \overline{f(G)}$.
\end{result}
\begin{demo}
Call $\pi$ the functional calculus of $T$ on $G \cap \si(T)$. Then result \ref{art3.res2} implies that $\pi$ is faithful.  
Put $g = f_{G \cap \si(T)} \in \C(G \cap \si(T))$. Then $f(T) = g(T) = \pi(g)$.

By the remark after proposition \ref{prel3.prop1}, this implies that $\si(f(T)) = \si(g(T)) = \si(\pi(g)) = 
\overline{g(G \cap \si(T))} = \overline{f(G \cap \si(T))}$.
\end{demo}

\bigskip

\begin{proposition} \label{art3.prop1}
Consider a Hilbert \cst-module over a \cst-algebra $A$. Let $B$ be a \cst-algebra and $\pi$  a non-degenerate 
$^*$-homomorphism from $B$ into $\cL(E)$. Consider a normal element $T$ affiliated with $B$ and $G$ a subset of 
$\,\C$ which is compatible with $T$. Then we have the following properties.
\begin{itemize}
\item $G$ is compatible with $\pi(T)$
\item We have for every $f \in C(G)$ that $f(\pi(T)) = \pi(f(T))$
\end{itemize}
\end{proposition}
\begin{demo}
Take a finite subset $K$ of $\,\C$ such that
\begin{enumerate}
\item The element $T - \lambda\,1$ is invertible for every $\lambda \in K$.
\item $\si(T) \setminus K \subseteq G$
\end{enumerate}
By proposition \ref{art2.prop2}, we have for every $\lambda \in K$ that $\pi(T) - \lambda \, 1$ is invertible. Furthermore, 
proposition \ref{prel3.prop1} implies that $\si(\pi(T)) \subseteq \si(T)$ so $\si(\pi(T)) \setminus K \subseteq \si(T) 
\setminus K \subseteq G$. So $G$ is compatible with $\pi(T)$.

Let $\th$ denote the functional calculus of $T$ on $G$. 
Then we have that $(\pi \th)(\io_G) = \pi(\th(\io_{G})) = \pi(T)$. 
So $\pi \th$ is the functional calculus of $\pi(T)$ on $G$.

This implies that 
$$f(\pi(T)) = (\pi \th)(f) = \pi(\th(f)) = \pi(f(T)) $$
for every $f \in C(G)$.
\end{demo}

\medskip

\begin{result} \label{art3.res4}
Consider a locally compact space $X$ and $f$ an element in $\text{C}(X)$. Let $G$ be an almost closed  subset of $\,\C$. Then the following holds.
\begin{itemize}
\item $G$ is compatible wih $f$ $\Leftrightarrow$ $f(X) \subseteq G$
\item Suppose that $G$ is compatible with $f$. Then we have that $g(f) = g \!\circ\! f$ for every $g \in \text{C}(G)$.
\end{itemize}
\end{result}
\begin{demo}
Suppose that $G$ is compatible with $f$. Then there exists a finite $K$ of $\,\C$ such that
\begin{itemize}
\item We have for every $\lambda \in \C$ that $f - \lambda \, 1$ is invertible.
\item $\si(f) \setminus K \subseteq G$.
\end{itemize}
By corollary \ref{art2.cor2}, we know for every $\lambda \in K$ that $\lambda \notin f(X)$. Because moreover $\si(f) = \overline{f(X)}$, 
we get that $f(X) \subseteq \si(f) \setminus K$. So $f(X) \subseteq G$.

Define the $^*$-homomorphism $\pi$ from $\text{C}_0(G)$ into $\text{C}_b(X)$ such that $\pi(g) = g \!\circ\! f$ for $g \in \text{C}_0(G)$. 
Then $\pi$ is non-degenerate and $\pi(g) = g\!\circ\! f$ for $g \in \text{C}(G)$.

We have in particular that $\pi(\io_G) = f$, so $\pi$ is the functional calculus of $f$ on $G$.

\medskip

Suppose on the other hand that $f(X) \subseteq G$. 
Because $G$ is almost closed, there exists a finite subset $K$ of $\,\C \setminus G$ such that $G \cup K$ is closed. 
Hence, $\si(f) = \overline{f(X)} \subseteq G \cup K$.
Therefore, $\si(f) \setminus K \subseteq G$.

We have for every $\lambda \in K$ that $\lambda \notin G$, so $\lambda \notin f(X)$ which by corollary \ref{art2.cor2} 
implies that $f - \lambda \, 1$ is invertible.

So we see that $G$ is compatible with $f$.
\end{demo}

\medskip

If we combine the previous result with proposition \ref{art3.prop1}, we get the following familiar one.

\begin{proposition} \label{art3.prop2}
Consider a Hilbert \cst-module $E$ over a \cst-algebra $A$ and let $T$ be a normal element  in $\cR(E)$. Let $F$ be a subset 
of $\,\C$ which is compatible with $T$ and $f$ an element
in $\text{C}(F)$. Consider an almost closed subset $G$ of $\,\C$  such that $f(F) \subseteq G$.
Then
\begin{itemize}
\item  The set $G$ is compatible with $f(T)$.
\item We have that $g(f(T)) = (g\!\circ\!f)(T)$ for every $g \in \text{C}(G)$.
\end{itemize}
\end{proposition}

\bigskip

\begin{remark} \rm
Consider a Hilbert space $H$ and a normal element $T$ in $\cR(H)$. Let $G$ be a subset of $\,\C$ which is compatibel with $H$. 
Then $G$ is certainly Borel measurable. Let $\cB(G)$ denote the collection of Borel subsets of $G$.

Denote the spectral measure of $T$ on $\,\C$ by $F$ and call $E$ the restriction of $F$ to  $\cB(G)$. 

By assumption, there exists a finite subset $K$ of $\,\C$ such that
\begin{itemize}
\item We have for every $\lambda \in K$ that $T - \lambda \, 1$ is invertible.
\item $\si(T) \setminus K \subseteq G$
\end{itemize}
We have for every $\lambda \in K$ that $\lambda$ is not an eigenvalue of $T$ which implies that
$F(\{\lambda\}) = 0$. Because $K$ is finite, this implies that $F(K) = 0$. Hence we get that $F((G \setminus \si(T)) \cup 
(\si(T) \setminus G)) = 0$. 

So $E$ is a spectral measure on $G$ such that $\int \io_G \, dE = T$. We call $E$ the spectral measure of $T$ on $G$.
\end{remark}

\bigskip

\begin{proposition} \label{art3.prop3}
Consider a Hilbert space $H$ and $T$ a normal  operator in $H$. Let $G$ be a subset of $\,\C$ which is compatible with $T$ 
and call $E$ the spectral measure of $T$ on $G$. Then $f(T) = \int f \,dE$ for every $f \in \text{C}(G)$.
\end{proposition}
\begin{demo}
Define the $^*$-homomorphism $\pi$ from $\text{C}_0(G)$ into $\cB(H)$ such that
$\pi(f) = \int f \, dE$ for $f \in \text{C}_0(G)$.

Take an approximate unit $(e_i)_{i \in I}$ for $\text{C}_0(G)$ in $\text{K}(G)$.
\begin{trivlist}
\item[\ \,\,$\bullet$] First we prove that $\pi$ is non degenerate.

\vspace{1mm}

Choose $v \in H$. Then we have a regular Borel measure $E_v$ on $G$.

The set $\{ \, e_i \mid i \in I \, \}$ is an upwardly directed set of continuous positive functions such that the equality 
$\sup\,\{ \, e_i(c) \mid i \in I \,\} = 1$ holds for every $c \in G$.

Hence, the regularity of $E_v$ implies that
\ $\sup\,\{\,\int e_i \, d E_v \mid i \in I \, \} = \int 1 \, dE_v = \|v\|^2$.

\vspace{1mm}

So we see that $\bigl(\, \int e_i \, d E_v \, \bigr)_{i \in I}$ converges to 1.

We have moreover for every $i \in I$ that
\begin{eqnarray*}
& & \| v - \pi(e_i)v \|^2 =  \| v - \bigl(\,\int e_i \, dE\,\bigr)\, v \|^2
= \| \bigl(\,\int \,1 - e_i \,\, dE\,\bigr) \, v \|^2 \\
& & \spat = \int |1-e_i|^2 \, d E_v \leq \int |1-e_i| \, d E_v 
= \int \,1 - e_i \,\, d E_v = 1 - \int e_i \, d E_v
\end{eqnarray*}

This implies that $(\pi(e_i)v)_{i \in I}$ converges to $v$.
\item[\ \,\, $\bullet$] Take $g \in \text{C}(G)$. We will prove now that $\pi(g) = \int g \, dE$.

\vspace{1mm}

Fix $j \in I$. Then the spectral functional calculus rules imply that
$\bigl(\,\int e_j \, dE\,\bigr) \, H \subseteq D\bigl(\,\int g \, dE\,\bigr)$ and that
$$\bigl(\,\int e_j \, dE\,\bigr) \, \bigl(\,\int g \, dE\,\bigr)
\subseteq \bigl(\,\int g \, dE\,\bigr)\, \bigl(\,\int e_j \, dE\,\bigr)
=\int g e_j \, dE$$
This impies immediately that $\pi(e_j) \, H \subseteq D\bigl(\,\int g \, dE\,\bigr)$
and that
$$\pi(e_j) \, \bigl(\,\int g \, dE\,\bigr)
\subseteq \bigl(\,\int g \, dE\,\bigr)\, \pi(e_j)
=\pi(g e_j)$$
By the functional calculus rules of section \ref{art1}, we have moreover that $\pi(e_j) \, H$ subseteq $\pi(g)$ and that
$$\pi(e_j)\,\pi(g) \subseteq \pi(g) \, \pi(e_j) = \pi(g e_j)$$

\medskip

Choose $v \in D\bigl(\,\int g \, dE\,\bigr)$. By the remarks above, we know for every $i \in I$ that $\pi(e_i) v \in D(\pi(g))$ and that
$$\pi(g)(\pi(e_i) v) = \pi(g e_i) v = \pi(e_i) \, \bigl(\,\int g \, dE\,\bigr)\,v$$
This implies that $(\pi(e_i)v)_{i \in I}$ converges to $v$ and that
$\bigl(\pi(g)(\pi(e_i)v)\bigl)_{i \in I}$ converges to $\bigl(\,\int g \, dE\,\bigr)\,v$.

So the closedness of $\pi(g)$ implies that $v \in D(\pi(g))$ and
$\pi(g) \, v = \bigl(\,\int g \, dE\,\bigr)\,v$.

\vspace{1mm}

Consequently, $\pi(g) \subseteq \int g \, dE$. We get in a similar way that
$\int g \, dE \subseteq \pi(g)$. Hence, $\int g \, dE = \pi(g)$.
\end{trivlist}
\medskip
So we have in particular that $\pi(\io_G) = \int \io_G \, dE = T$ which implies that
$\pi$ is the functional calculus of $T$ on $G$.
\end{demo}

\bigskip

\section{Powers of positive regular operators} \label{art4}

In this section, we show that the power theory of positive closed operators has an obvious generalization to a power theory of positive 
regular operators. By the previous sections, the proofs of the results in this section can be copied from the corresponding proofs in 
the Hilbert space case and are always simple consequences of calculation rules from sections \ref{art1} and \ref{art3}. As a consequence, 
this section contains not many proofs

For the sake of completeness, we will give a list of the basic calculation rules.

\bigskip

For the most part of this section, we will fix a Hilbert \cst-module $E$ over a \cst-algebra $A$.

\medskip

Let us start first with a basic result.

\begin{result} \label{art4.res1}
Consider a normal regular operator $T$ in $E$. Let $n \in \N \cup \{0\}$ and define the continuous function $f$ from $\,\C$ into $\,\C$ 
such that $f(c) = c^n$ for every $c \in \C$. Then $T^n = f(T)$.
\end{result}
\begin{demo} We proceed by induction.

Define for every $m \in \N \cup \{0\}$ the function $f_m$
from $\,\C$ into $\,\C$ such that $f_m(c) = c^m$ for every $c \in \C$.
\begin{itemize}
\item We have that $f_0(T) = 1 = T^0$.
\item Let $m \in \N \cup \{0\}$ and suppose that $f_m(T) = T^m$.

By definition of the functional calculus, the equality $f_1(T) = T$ holds.

We have that $f_{m+1}(T) = (f_{m} f_{1})(T) \supseteq f_{m}(T) f_{1}(T)
= T^m T = T^{m+1}$.
Using lemma \ref{art1.lem2}, it is not so difficult to see that $D(f_{m+1}(T)) \subseteq
D(f_1(T))$. So we get that
\begin{eqnarray*}
D(T^{m+1}) & = & D(f_m(T) f_1(T)) = D( (f_m f_1)(T) ) \cap D( f_1(T)) \\
&  = & D(f_{m+1}(T)) \cap D(f_1(T)) = D(f_{m+1}(T)) \ .
\end{eqnarray*}
Consequently, $f_{m+1}(T) = T^{m+1}$.
\end{itemize}
\end{demo}

\begin{corollary}
Consider a normal regular operator $T$ in $E$ and  $n \in \N \cup \{0\}$. Then $T^n$ is a normal regular operator in $E$.
\end{corollary}

\medskip

\begin{result}
Consider an invertible normal regular operator $T$ in $E$. Let $n \in \Z$ and define the continuous function $f$ from $\,\C_0$ into 
$\,\C_0$ such that $f(c) = c^n$ for every $c \in \C_0$. Then $T^n = f(T)$.
\end{result}
\begin{demo}
If $n \geq 0$, then the result follows from the previous result. Therefore suppose that $n < 0$. Define the continuous function 
$g,h$ from $\,\C_0$ into $\,\C_0$ such that $g(c) = \frac{1}{c}$ and $h(c) = c^{-n}$ for every $c \in \C_0$. Then $h \!\circ\! g 
= f$, so $f(T) = h(g(T))$ by proposition \ref{art3.prop2}.

Hence, the previous result implies that $f(T) = g(T)^{-n}$. \hspace{0.5cm} (*)

Define the continuous function $k$ on $\,\C_0$ such that $k(c) = c$ for every $c \in \C_0$ , so $k(T) = T$. We have moreover that 
$g k = 1$. By results \ref{art1.res2} and  \ref{art1.res3}, this implies that 
$$g(T) \, T = g(T) \, k(T) \subseteq (g k)(T) = 1$$
and 
$$D(g(T)\,T) = D( g(T) \, k(T) ) = D( (g k)(T) )
\cap D(k(T)) =  E \cap D(T) = D(T) $$
We get in a similar way that $T \,g(T) \subseteq 1$ and $D(T \,g(T)) = D(g(T))$.

So we see that $g(T) = T^{-1}$. By (*), this gives us that  $f(T) = (T^{-1})^{-n} = T^n$.
\end{demo}

\medskip

\begin{corollary}
Consider an invertible normal regular operator $T$ in $E$ and  $n \in \Z$. Then $T^n$ is an invertible normal regular operator in $E$.
\end{corollary}

\bigskip

Using these results and the results of the previous sections, we are now able to define powers of positive affiliated elements and 
prove the most important results.

\begin{definition}
Consider a positive  regular operator $T$ in $E$. Let $s$ be a positive number and define the continuous function $f$ from $\R^+$ 
into $\R^+$ such that $f(t) = t^s$ for every $t \in \R^+$. We define $T^s = f(T)$, so $T^s$ is a positive element affiliated with $A$.
\end{definition}

Due to result \ref{art4.res1}, this definition is consistent with the usual notion of $T^s$
when $s$ belongs to $\N \cup \{0\}$. Due to proposition \ref{art3.prop3}, this definition is also consistent with the usual 
definition of $T^s$ if $E$ is a Hilbert space.

\bigskip

The following calculation rules hold:

\begin{proposition}
Consider a  positive  regular operator $T$ in $E$.  Then the following elementary properties hold:
\begin{enumerate}
\item We have for every $s,t \in \R^+$ that $(T^s)^t = T^{s t}$.
\item We have for every  $s,t \in \R^+$ with $s \leq t$ that $D(T^t) \subseteq D(T^s)$.
\item We have for every  $s,t \in \R^+$ that $T^{s+t} = T^s \, T^t$.
\end{enumerate}
\end{proposition}

The first result follows from proposition \ref{art3.prop2} and the second from lemma
\ref{art1.lem2}. The third equality follows from results \ref{art1.res2} and \ref{art1.res3}, using the second
result.

\medskip

\begin{remark} \rm
Let $S$,$T$ be  positive regular operators in $E$ and $r$ a strictly positive number. The first result of the previous proposition implies 
that $S=T$ $\Leftrightarrow$ $S^r = T^r$.

Therefore, we have for every regular positive  operator $T$ in $E$ that $T^\frac{1}{2}$ is
the unique positive element in $\cR(E)$ such that $(T^\frac{1}{2})^2 = T$.
\end{remark}

\bigskip

\begin{result}
Consider a Hilbert \cst-modules $E$ over a \cst-algebra $A$ and a positive element $T \in \cR(E)$. Then we  have for every $\al \in \R^+_0$ 
that $\overline{\bld T}= \overline{\bld T^\al}$.
\end{result}
\begin{demo}By proposition \ref{art1.prop1}, we have that $\overline{\bld T} = \overline{\bld T^{2^n}}$ for every $n \in \Z$.

We have for $\be, \gamma \in \R^+_0$ with $\be \leq \gamma$ that
$T^\gamma = T^{\be + (\gamma - \be)} = T^\be \, T^{\gamma - \be}$ which implies that
$\overline{\bld T^\gamma} \subseteq \overline{\bld T^\be}$.

Suppose that $\al \geq 1$. Then there exists $m \in \N$ such that $\al \leq 2^m$.
So we get that 
$$\overline{\bld T} = \overline{\bld T^{2^m}} \subseteq \overline{\bld T^\al}
\subseteq \overline{\bld T}$$
so $\overline{\bld T^\al} = \overline{\bld T}$.

The case $\al \leq 1$ is treated in a similar way.
\end{demo}

\medskip

\begin{corollary}
Consider a Hilbert \cst-modules $E$ over a \cst-algebra $A$ and a positive element $T \in \cR(E)$. Then we  have for every $\al \in \R^+_0$ 
that $\ker T = \ker T^\al$.
\end{corollary}

This follows because $\ker T = (\bld T)^\perp$ and $\ker T^\al = (\bld T^\al)^\perp$.

\bigskip\medskip

An important application of the functional calculus of section \ref{art3}, can be found in the following definition.

\begin{definition}
Consider  a strictly positive regular operator $T$ in $E$.  Let $z$ be a complex number. Define the function $f$ from $\R^+_0$ into $\,\C_0$ 
such that $f(t) = t^z$ for every $t \in \R^+_0$. Then we define $T^z = f(T)$, so $T^z$ is an invertible normal regular operator in $E$.
\end{definition}

\begin{proposition}
Consider a strictly positive regular operator $T$ in $E$. Then we have the following calculation rules.
\begin{enumerate}
\item We have for every $s \in \R$ that $T^s$ is a strictly positive regular operator in $E$.
\item We have for every $s \in \R$ that $T^{is}$ is a unitary element of $\cL(E)$.
\item Let $s$ be a real number and $z$ a complex number. Then $(T^s)^z = T^{s z}$.
\item Let $n$ be an integer and $z$ a complex number. Then $(T^z)^n = T^{n z}$.
\item We have for every $z \in \C$ that $(T^z)^* = T^{\overline{z}}$.
\item Consider complex numbers $y,z$ such that $\text{Im\,}y$ lies between 0 and $\text{Im\,}z$.
Then $D(T^z) \subseteq D(T^y)$.
\item Consider $y,z \in \C$. Then $T^y\,T^z$ is closable and the closure is equal to $T^{y+z}$.
We have moreover that $D(T^y \, T^z) = D(T^z) \cap D(T^{y+z})$  and
$\bld(T^y \, T^z) = \bld(T^y) \cap \bld(T^{y+z})$.
\item Consider $y,z  \in \C$ such that $y$ and $z$ lie at the same side of the real line.
Then $T^y \, T^z = T^{y+z}$.
\item We have for every $z \in \C$ and $s \in \R$ that $T^{z+is} = T^z T^{is} = T^{is} T^z$.
\end{enumerate}
\end{proposition}

\medskip

\begin{corollary}
Consider a positive regular operator $T$ in $E$ which is adjointable invertible. Then
\begin{enumerate}
\item We have for every $z \in \C$ with $\text{Re\,}z \leq 0$ that $T^z$ belongs to $\cL(E)$.
\item We have for every $z \in \C$ with $\text{Re\,}z \geq 0$ that $T^z$ is adjointable invertible.
\end{enumerate}
\end{corollary}

\medskip

\begin{proposition}
Consider a strictly positive regular operator $T$ in $E$.  Then we have that the function \newline $\R \rightarrow \cL(E) 
: s \mapsto T^{is}$ is a strongly continuous unitary group homomorphism.
\end{proposition}

\bigskip\medskip

As usual, we can switch between strictly positive elements and selfadjoint elements using the exponential and the logarithm.

\begin{definition}
Consider a normal regular operator $T$ in $E$. Define the function $f$ from $\,\C$ into $\,\C$ such that $f(c) = \text{e}^{\,c}$ 
for every $c \in \C$. Then we define $\text{e}^{\,T} = f(T)$, so $\text{e}^{\,T}$ is a normal operator in $E$.
\end{definition}

\vspace{1mm}

Proposition \ref{prel1.prop3} implies immediately the following result.

\begin{corollary}
Consider a selfadjoint regular operator $T$ in $E$. Then $\text{e}^{\,T}$ is a strictly positive regular operator in $E$.
\end{corollary}

\bigskip

\begin{definition}
Consider a strictly positive regular operator in $E$. Define the function $f$ from $\R^+_0$ into $\R$ such that $f(t) = \ln t$ 
for every $t \in \R^+_0$. Then we define $\ln T = f(T)$, so $\ln T$ is a selfadjoint regular operator in $E$.
\end{definition}

By using cuttings of the plane, we can of course define the logarithm also for a wider class of normal regular operators.

\medskip

As usual, proposition \ref{art3.prop2}  implies the following result.

\begin{result}
We have the following properties :
\begin{enumerate}
\item We have for every selfadjoint element $T \in \cR(E)$ that $\ln(\text{e}^{\,T}) = T$.
\item We have for every strictly positive element $T \in \cR(E)$ that $\text{e}^{\,\ln T} = T$.
\end{enumerate}
\end{result}

\bigskip\medskip

So we see in fact that the functional calculus theory on Hilbert spaces can be rather succesfully generalized to a functional 
calculus theory on Hilbert \cst-modules.

\bigskip\bigskip

\begin{proposition} \label{art4.prop1}
Consider Hilbert \cst-modules $E$,$F$ over a \cst-algebra $A$ and let  $T$ be an element  
in $\cR(E,F)$. Then we define $|T| = (T^* T)^\frac{1}{2}$, so $|T|$ is a positive 
element in $\cR(E)$. We have moreover that $D(T) = D(|T|)$ and that $\lan T(v) , T(w) \ran = \lan |T|(v) , |T|(w) \ran$ for 
every $v,w \in D(T)$.
\end{proposition}
\begin{demo} Because $T^* T = |T|^2$, the set $D(T^* T)$ is a core for $T$ and $|T|$.
It is straightforward to check that $\lan T(v) , T(w)
\ran = \lan |T|(v) , |T|(w) \ran$ for every $v,w \in D(T^* T)$.

This implies also that $\|T(v)\| = \| \, |T|(v)\|$ for every $v \in D(T^* T)$. Using this
equality and the fact that $D(T^* T)$ is a core for both $T$ and $|T|$, it is easy to check
that $D(T) = D(|T|)$.

Also the last equality can now be easily deduced.
\end{demo}

\medskip

\begin{result}
Consider Hilbert \cst-modules $E$,$F$ over a \cst-algebra $A$ and an element $T \in \cR(E,F)$.
Then  $\overline{\bld |T|} = \overline{\bld T^*}$ and $\ker |T| = \ker T$ 
\end{result}

The result concerning the image follows by applying proposition \ref{art1.prop1} to $|T|$ and $T^*$. The result concerning the 
kernels follows from proposition \ref{art4.prop1}.

\bigskip

By result \ref{prel1.res1} and the functional calculus rules, we have immediately the following result.

\begin{result}
Consider a Hilbert \cst-module $E$ over a \cst-algebra $A$. Let $B$ be 
a \cst-algebra and $\pi$ a non-degenerate $^*$-homomorphism from $B$ into $\cL(E)$. Then we have for every element $T$ affiliated 
with $B$ that $\pi(|T|) = |\pi(T)|$.
\end{result}

\bigskip\bigskip

In chapter 9 of \cite{Lan}, Lance proved the next result (Remember that 
the adjoint invertibility of $1+T^* T$ is equivalent to the fact that   
$\bld (1+T^* T) = E$).

\begin{proposition}
Consider Hilbert \cst-modules $E$,$F$ over a \cst-algebra $A$ and $T \in \cR(E,F)$.
Then $1 + T^* T$ is adjointable invertible, $1-z_T^* z_T$ is invertible and
$(1 + T^* T)^{-1} = 1 -z_T^* z_T$. 
\end{proposition}

\medskip

He proved also the next result (this has in fact everything to do with the equivalence of the two approaches to regular operators).

\begin{proposition} \label{art4.prop2}
Consider Hilbert \cst-modules $E$,$F$ over a \cst-algebra $A$ and $T \in \cR(E,F)$.
Then 
$$z_T = T \, (1 + T^* T)^{-\frac{1}{2}} \hspace{1.5cm} \text{ and } \hspace{1.5cm}
T = z_T \, (1 -z_T^* z_T)^{-\frac{1}{2}} $$
\end{proposition}

\bigskip

To be more precise, it is proven in \cite{Lan} that
$$z_T = T \, \bigl((1 + T^* T)^{-1}\bigr)^{\frac{1}{2}} \hspace{1.5cm} \text{ and } \hspace{1.5cm} T = z_T \, 
\bigl((1 -z_T^* z_T)^{\frac{1}{2}}\bigr)^{-1} $$
but these two results are the same by the results of this section.

\bigskip

If $T$ is a normal, the results above and the functional calculus rules imply that
$z_T$ is a function of $T$.

\begin{corollary}
Consider a Hilbert \cst-module $E$ over a \cst-algebra $A$ and  a normal element $T$ in $\cR(E)$. Define the continuous function 
$f$ from $\,\C$ into $\,\C$ such that $f(c) = \frac{c}{1 + |c|^2}$ for every $c \in \C$. Then we have that $z_T = f(T)$.
\end{corollary}

\bigskip

\begin{remark} \label{art4.rem1} \rm
Define $D = \{ \, c \mid c \in \C \text{ such that } |c| < 1 \, \}$. Then we can define the function $g$ from $D$ into $\,\C$ such 
that $g(c) = \frac{c}{1 - |c|^2}$ for $c \in D$. So $g$ is nothing but the inverse of $f$. Then we should also have that $T = g(z_T)$ 

\vspace{1mm}

But this needs a more general kind of functional calculus than introduced in section \ref{art3}. The finite set $K$ appearing in 
this section should be replaced by a closed set (which in this case would be the unit circle). 

One should find a general (and easy to check) condition  on elements of $K$ (or probably, on the set $K$) such that proposition 
\ref{art3.prop4} remains true. 

\vspace{1mm}

The condition stated in proposition \ref{art3.prop4} is not strong enough. Just look at the case of a normal operator in a Hilbert 
space which has no eigenvalues in which case $K$ can be any closed subset of the spectrum.

In the classical commutative case however, the condition in proposition \ref{art3.prop4} is strong enough because it is equivalent 
to saying that the range of the function $T$ does not meet $K$.
\end{remark}

\bigskip

\section{The Fuglede-Putnam theorem for Hilbert \cst-modules}

A well known result in Hilbert space theory is the Fuglede-Putnam theorem. We will prove in this section the Fuglede-Putnam theorem 
for regular normal operators between Hilbert \cst-modules. The proof is modelled on the proof for the Hilbert space case 
(see exercices 5 and 8 on page 334 of \cite{Conw}) but a little bit of care has to be taken.

\bigskip

The next result guarantees that the procedure of cutting off unbounded operators remains a nice (and useful) operation in the Hilbert 
\cst-module framework.

\begin{lemma} \label{art5.lem1}
Consider a Hilbert \cst-module $E$ over a \cst-algebra $A$ and a normal element $T \in \cR(E)$.
Let $G$ be a subset of $\,\C$ such that $G$ is compatible with $T$. Take $f \in \text{K}(G)$
and call $K$ the support of $f$.
Define $F$ as the closure of $f(T) E$ in $E$, so $F$ is a sub Hilbert \cst-module of $E$.
Then we have the following properties.
\begin{itemize}
\item We have that $F \subseteq D(T), D(T^*)$, \ $T_F$ is a normal element in $\cL(F)$ and $(T_F)^* = (T^*)_F$.
\item The set $G$ is compatible with $T_F$.
\item Consider $g \in \text{C}(G)$. Then $F \subseteq D(g(T))$,\ $g(T)_F$ belongs to $\cL(F)$, $\|g(T)_F\| \leq \|g_K\|$ and $g(T)_F = g(T_F)$. 
\end{itemize}
\end{lemma}
\begin{demo} 
\begin{trivlist}
\item[\ \,$\bullet$] Choose $g \in \text{C}(G)$. We have that $|g f| \leq \|g_K\| \,\, |f|$.

Take $v \in E$. Then we have that $f(T) v$ belongs to $D(g(T))$ and $g(T)(f(T) v) = (g f)(T) v$.
So we see that
\begin{eqnarray*}
\lan g(T)(f(T) v) , g(T)(f(T) v) \ran & = & \lan (g f)(T)\,v , (g f)(T)\,v \ran
= \lan |g f|^2(T)\, v , v \ran \\
& \leq & \|g_K\|^2 \, \lan |f|^2(T) \, v , v \ran 
=  \|g_K\|^2 \, \lan f(T) v , f(T) v \ran  
\end{eqnarray*}
which implies that $\|g(T) (f(T) v)\| \leq  \|g_K\| \, \|f(T) v\|$. 

Hence, using the closedness of $g(T)$, this allows us to conclude  that $F \subseteq D(g(T))$ and that $\|g(T) w\| \leq \|g_K\| \, \|w\|$ for 
$w \in F$. So $g(T)_F$ is a bounded operator from $F$ into $E$ such that $\|g(T)_F\| \leq \|g_K\|$.

We have for every $v \in D(g(T))$ that $g(T)_F(f(T) v) = g(T)(f(T) v)
= f(T)(g(T) v)$ which implies that $g(T)_F(F) \subseteq F$. So $g(T)_F$ is a bounded operator on $F$.

\item[\ \,$\bullet$] Choose $g \in \text{C}(G)$. Then we have bounded operators $g(T)_F$,
$\overline{g}(T)_F$ on $F$. Because $g(T)^* = \overline{g}(T)$, it follows easily that $g(T)_F$ belongs to $\cL(F)$ and that $(g(T)_F)^* = 
\overline{g}(T)_F$

\item[\ \,$\bullet$] We get in particular that $F \subseteq D(T),D(T^*)$, that $T_F$ belongs to $\cL(F)$ and that $(T_F)^* = (T^*)_F$.
Because $T^* T = T T^*$, this implies also immediately that $T_F$ is normal.

\item[\ \,$\bullet$] Define the mapping $\pi$ from $\text{C}_0(G)$ into $\cL(F)$ such that
$\pi(g) = g(T)_F$ for every $g \in \text{C}_0(G)$. Then we get immediately that
$\pi$ is a $^*$-homomorphism.

Furthermore,
\begin{eqnarray*}
& & (\pi(\text{C}_0(G)) \, f(F) E\,)\streep  =  
[\,g(T) f(T) v \mid g \in \text{C}_0(G), v \in E \,] \\
& & \spat  =  [\,f(T) g(T) v \mid g \in \text{C}_0(G), v \in E \,] 
=  [\,f(T) w \mid w \in E\,] = F 
\end{eqnarray*}
which implies that $\pi$ is non-degenerate.

Take $g \in \text{C}(G)$. By lemma \ref{art1.lem3}, we know that $\pi(\text{K}(G)) F$ is a core for $\pi(g)$. We have moreover for every 
$v \in F$ and $h \in \text{K}(G)$ that
$$\pi(g)(\pi(h)v) = \pi(g h) v = (g h)(T) v = g(T)(h(T) v) = g(T) (\pi(h) v) = g(T)_F (\pi(h) v)$$

From this, it follows that $\pi(g) = g(T)_F$.

\medskip

So we have in particular that $\pi(\io_G) = T_F$. By proposition \ref{art3.prop1}, this implies that $T_F$ is compatible with $G$. It implies 
also that $\pi$ is the functional calculus of $T_F$. Hence, $g(T_F) = \pi(g) = g(T)_F$ for $g \in \text{C}(G)$.
\end{trivlist}
\end{demo}

\medskip

It is possible (and not too hard) to prove a similar result if $f$ belongs to $\text{C}(G)$, but $g(T)_F$ will then belong to $\cR(F)$ 
(and $D(g(T)_F) = D(g(T)) \cap F$).

\bigskip

We will need the previous lemma in its full generality later on. Let us start of with the following easy implication. 

\begin{lemma} \label{art5.lem2}
Consider a Hilbert \cst-module $E$ over a \cst-algebra $A$ and let $T$ be a normal element in $\cR(E)$. Let $r$ be a  strictly positive 
number, $f$ an element in $\text{K}(\C)$ such that
$f(c)=0$ for every $c \in \C$ with $|c| \geq r$. Define $F$ as the closure of $f(T) E$ in $E$.
Then we have for every $n \in \N$ and $w \in F$  that 
$w$ belongs to $D(T^n)$ and $\|T^n (w) \| \leq r^n \, \|w\|$.
\end{lemma}

This lemma follows from the previous lemma applied to the function $g \in \text{C}(\C)$ defined by $g(c) = c^n$ for $c \in \C$. We know 
for this function that $g(T) = T^n$ by result \ref{art4.res1}.

\medskip

Now we will prove a last lemma which is some kind of converse of the previous one.

\begin{lemma} \label{art5.lem3}
Consider a Hilbert \cst-module $E$ over a \cst-algebra $A$ and let $T$ be a normal element in $\cR(E)$. Let $r$ be a strictly positive 
number and $f$ an element in $\text{C}_b(\C)$ such that $f(c) = 1$ for every $c \in \C$ with $|c|\leq r$. Consider an element $w$ in $E$ 
such that we have for every $n \in \N$ that $w$ belongs to $D(T^n)$ and $\|T^n(w)\| \leq r^n \, \|w\|$. Then $w = f(T) w$.
\end{lemma}
\begin{demo} Define $D = \{ \, c \in \C \mid |c| \leq r \, \}$.

Take an element $h \in \text{K}(\C)$ such that $h$ has its support in 
$\,\C \setminus D$ and such that $0 \leq h \leq 1$.

Denote the support of $h$ by $M$, so $M$ is a compact subset of $\,\C \setminus D$. Then there exists a strictly positive number $s$ such 
that $s > r$ and such that $|c| > s$ for every $c \in M$. 

\begin{list}{}{\setlength{\leftmargin}{.4 cm}}

\item Take $n \in \N$. 

Define the element $g_n \in \text{C}(\C)$ such that $g_n(c) = c^n$ for every $c \in \C$. Then result \ref{art4.res1} implies that 
$T^n = g_n(T)$. So we have by assumption that $w$ belongs to $D(g_n(T))$ and that we have the inequality $\|g_n(T) w\| \leq r^n 
\, \|w\|$. \hspace{0.5cm} (*)

It is clear that $|g_n(c)| > s^n$ for every $c \in M$. By lemma \ref{art1.lem4}, this implies the existence of  an element 
$h_n \in \text{K}(\C)$ such that $0 \leq h_n \leq 1$,
$h_n = 1$ on $M$ and $\lan g_n(T) u , \lan g_n(T) u \ran \geq s^n \, \lan h_n(T) u , \lan h_n(T) u \ran$ for every $u \in D(g_n(T))$.

So, using (*), we get that
$$r^n \, \|w\| \geq \|g_n(T) w\| \geq  s^n \, \|h_n(T) w\| \geq s^n \, \|h(T) w \|$$
where we used the fact that $h \leq h_n$ in the last inequality.

Hence, we see that $\|h(T) w\| \leq (\frac{r}{s})^n \, \|w\|$.

\end{list}

Because $\frac{r}{s} < 1$, this implies that $\|h(T) w\| = 0$.
So we get that $h(T) w = 0$.

\medskip

Now take an approximate unit $(e_k)_{k \in K}$ for $\text{C}_0(\C \setminus D)$ in $\text{K}(\C \setminus D)$. We define for every 
$k \in K$ the element $d_k \in \text{K}(\C)$ such that $d_k = 0$ on $D$ and $d_k = e_k$ on $\,\C \setminus D$. By the first part of 
the proof, we know that $d_k(T) w = 0$ for every $k \in K$.

\medskip

Choose $p \in \C_0(\C)$. Denote the restriction of $f$ and $p$ to $\,\C \setminus D$ by 
$\tilde{f}$ and $\tilde{p}$ respectively. 

Then $(1 - \tilde{f}) \tilde{p}$ belongs to $\text{C}_0(\C \setminus D)$ (the first factor takes care of the behaviour in the 
neighbourhood of the boundary of $D$, the second factor
takes care of the behaviour for large complex numbers).

Because $f = 1$ on $D$, we have moreover for every $k \in K$ that
$$\| (1-f)d_k\,p - (1-f) p \| = \|(1-\tilde{f}) e_k \, \tilde{p} - (1 - \tilde{f}) \tilde{p}\| \ .$$
This implies that $\bigl((1- f)d_k\, p \bigr)_{k \in K}$ converges to $(1-f) p$.

So we see that $\bigl(\,(1- f)d_k\,\bigr)_{k \in K}$ is a bounded net which converges strictly  to $1-f$. Consequently, the net 
$\bigl(\,[(1-f)d_k](T)\,\bigr)_{k \in K}$ converges strongly to $(1-f)(T)$.

\medskip

We have for every $k \in K$ that $ [(1-f)d_k](T) w = (1-f)(T) d_k(T) w = 0$. This implies that
$(1-f)(T) w =0$, so $f(T) w = w$.
\end{demo}

\medskip

First we need the bounded form of the Fuglede-Putnam theorem. The proof of this one can be directly copied from the proof of the 
Hilbert space form (see e.g. the proof of theorem 6.7 of \cite{Conw}).

\begin{proposition} \label{art5.prop1}
Consider  Hilbert \cst-modules $E$,$F$ over a \cst-algebra $A$ and $U$ a bounded linear mapping from $E$ into $F$. Let $S$ be  a 
normal element in $\cL(E)$ and $T$ a normal element in $\cL(F)$ such that $U \, S = T \, U  $. Then $U \, S^* = T^* \, U$.
\end{proposition}

\medskip

\begin{corollary} \label{art5.cor1}
Consider  Hilbert \cst-modules $E$,$F$ over a \cst-algebra $A$ and $U$ a bounded linear mapping from $E$ into $F$. Let $S$ be  
a normal element in $\cL(E)$ and $T$ a normal element in $\cL(F)$ such that $U \, S = T \, U  $. Then we have for every 
$f \in \text{C}(\si(S) \cup \si(T))$ that $U \, f(S) = f(T) \, U$.
\end{corollary}

\bigskip

Now we will prove the unbounded forms of these results.

\begin{proposition} \label{art5.prop2}
Consider  Hilbert \cst-modules $E$,$F$ over a \cst-algebra $A$ and $U$ a bounded linear mapping from $E$ into $F$. Let $S$ 
be  a normal element in $\cR(E)$ and $T$ a normal element in $\cR(F)$ such that $U \, S \subseteq T \, U  $. Then $U \, S^* \subseteq T^* \, U$.
\end{proposition}
\begin{demo} We may suppose that $\|U\| \leq 1$. Define for every $m \in N$ the set $D_m = \{ \, c \in \C \mid |c| \leq m \, \}$.

\medskip

Choose $m \in \N$ and take an element $e_m \in K(\C)$ such that $0 \leq e_m \leq 1$, such that
$e_m = 1$ on $D_m$ and  \newline $e_m = 0$ on $\,\C \setminus D_{m+1}$ (such an element clearly exists).
\begin{itemize}
\item  We define $E_m$ as the closure of $e_m(S) E$ in $E$. So $E_m$ is a sub Hilbert \cst-module of  $E$. By lemma \ref{art5.lem1}, 
we know that $E_m \subseteq D(S), D(S^*)$. Denote the restrictions of $S$ and $S^*$ to $E_m$ by $S_m$ and $(S^*)_m$ respectively.

This same lemma implies then moreover that $S_m$ belongs to $\cL(E_m)$ and that
$(S_m)^* = (S^*)_m$.

\item Analogously, we define $F_m$ as the closure of $e_m(T) F$ in $F$.
Then we have that $F_m \subseteq D(T), D(T^*)$ and we denote the restrictions of
$T$ and $T^*$ to $F_m$ by $T_m$ and $(T^*)_m$ respectively.

Then we have again that $T_m$ belongs to $\cL(F_m)$ and that $(T_m)^* = (T^*)_m$.
\end{itemize}

Choose $v \in D(S^*)$.

\medskip

We will further reduce the problem to a bounded one. Therefore take $m \in \N$. 

\begin{list}{}{\setlength{\leftmargin}{.4 cm}}

\item Take $w \in E_m$. Choose $n \in \N$. 

By lemma \ref{art5.lem2}, we know that $w$ belongs to $D(S^n)$ and that $\|S^n(w)\| \leq (m+1)^n \, \|w\|$.

Because $U \, S \subseteq T \, U$, we have also that $U \, S^n \subseteq T^n \, U$. This implies that $U(w)$ belongs to $D(T^n)$ and that
$$\|T^n(U(w))\| = \|U(S^n(w))\| \leq \|S^n(w)\| \leq (m+1)^n \, \|w\| \ .$$

By lemma \ref{art5.lem3}, this implies that $U(w) = e_{m+1}(T)\,U(w)$, so $U(w)$ belongs to $F_{m+1}$.

\end{list}

Define $U_m$ as the restriction of $U$ to $E_m$. Then the above argument shows that
$U_m$ is a bounded linear operator from $E_m$ into $F_{m+1}$.

The fact that $U \, S \subseteq T \, U$ implies then immediately that
$U_m \, S_m = T_{m+1} \, U_m$. So proposition \ref{art5.prop1} implies that $U_m \, (S_m)^* = (T_{m+1})^* \, U_m$.
Hence we get that $U_m \, (S^*)_m = (T^*)_{m+1} \, U_m$.

\medskip

We have that $U(e_m(T) v)$ belongs to $F_{m+1}$, so $U(e_m(T) v)$ belongs to $D(T^*)$. Furthermore,
\begin{eqnarray*}
T^*(U (e_m(S) v)) & = & (T^*)_{m+1}(U_m (e_m(S) v)) = U_m((S^*)_m (e_m(S) v)) \\
 & = & U(S^*(e_m(S) v)) = U \bigl(e_m(S)(S^*(v))\bigr) 
\end{eqnarray*}

So we see that $\bigl(U(e_m(S) v)\bigr)_{m=1}^\infty$ converges to $U(v)$ and that
$\bigl(T^* (U (e_m(S)v))\bigr)_{m=1}^\infty$ converges to $U(S^*(v))$.

Therefore, the closedness of $T^*$ implies that $U(v)$ belongs to $D(T^*)$ and that
$T^*(U(v)) = U(S^*(v))$.
\end{demo}

\medskip

\begin{lemma} \label{art5.lem4}
Consider  Hilbert \cst-modules $E$,$F$ over a \cst-algebra $A$ and a bounded linear mapping $U$ from $E$ into $F$. Let $B$ be a 
\cst-algebra, $\pi$ a non-degenerate $^*$-homomorphism from $B$ into $\cL(E)$ and $\th$ a non-degenerate $^*$-homomorphism from 
$B$ into $\cL(F)$ such that $U \, \pi(b) = \th(b) \, U$ for every $b \in B$
Then $U \, \pi(T) \subseteq \th(T) \, U$ for every $T \in \cR(B)$
\end{lemma}
\begin{demo}
Choose $b \in D(T)$ and $v \in E$. 

Now  $U (\pi(b) v) = \th(b) (U v)$, so $U (\pi(b) v)$ belongs to $D(\th(T))$ and
$$\th(T)(U (\pi(b) v)) = \th(T)(\th(b) (U v)) = \th(T(b)) (U v) = U ( \pi(T(b)) v) \ .$$

We also have that $\pi(b) v$ belongs to $D(\pi(T))$ and $\pi(T) ( \pi(b)  v) = \pi(T(b)) v$.
So we see that $\th(T)(U (\pi(b) v))$ $= U( \pi(T) (\pi(b) v))$.

\medskip

Using the closedness of $\th(T)$ and the fact that
$\lan \, \pi(b) v \mid b \in D(T), v \in E \, \ran$ is a core for $\pi(T)$, we get now easily that $U \, \pi(T) \subseteq \th(T) \, U$.
\end{demo}

\medskip

Now we can prove the final form of the Fuglede-Putnam theorem.

\begin{theorem} \label{art5.thm1}
Consider  Hilbert \cst-modules $E$,$F$ over a \cst-algebra $A$ and a bounded linear mapping $U$ from $E$ into $F$. Let $S$ be  
a normal element in $\cR(E)$ and $T$ a normal element in $\cR(F)$ such that $U \, S \subseteq T \, U  $.
Consider a subset $G$ of $\,\C$ which is compatible with both $S$ and $T$. Then we have for every $f \in \text{C}(G)$ that 
$U \, f(S) \subseteq f(T) \, U$.
\end{theorem}
\begin{demo}
By proposition \ref{art5.prop2}, we know that also  $U \, S^* \subseteq T^* \, U$. This will imply that $U \,(S^* S) \subseteq (T^* T) \, U$.
In turn, this will imply that $U \, (1+S^* S) \subseteq (1+T^* T) \, U$.
From this, we get that $U \, (1+S^* S)^{-1} \subseteq (1+T^* T)^{-1} \, U$.

But $(1+S^* S)^{-1}$ and $(1+T^*T)^{-1}$ are positive elements in $\cL(E)$ and $\cL(F)$ respectively. So $U \, (1+S^* S)^{-1} = (1+T^* T)^{-1} \, U$. 
Hence, corollary \ref{art5.cor1} implies that $U \, \bigl((1+S^* S)^{-1}\bigr)^\frac{1}{2} = \bigl((1+T^* T)^{-1}\bigr)^\frac{1}{2} \, U$.

This will in turn imply that $U \, S \, \bigl((1+S^* S)^{-1}\bigr)^\frac{1}{2} \subseteq T \, \bigl((1+T^* T)^{-1}\bigr)^\frac{1}{2} \, U$.

So, using proposition \ref{art4.prop2}, we get that $U z_S \subseteq z_T \, U$ which implies that  $U  z_S = z_T \, U$.

\medskip

Let us now take a function $g \in \text{C}_0(G)$. Then there exist $h \in \text{C}_0(\C)$ such that $h_G = g$.

Define $D = \{\,z \in \C \mid |z| \leq 1 \, \}$ and define the the homeomorphism $J$ from 
$\,\C$ to $D^0$ such that $J(c) = \frac{c}{(1+|c|^2)^\frac{1}{2}}$ for every $c \in \C$.
Then there exists a unique element $k \in \text{C}(D)$ such that $k = 0$ on $\partial D$ and
such that $k\!\circ\!J = h$. By equality 1.19 of \cite{Wor6}, we have that $k(z_S) = h(S)$ and $k(z_T) = h(T)$.

Because $U z_S = z_T \, U$, we know by corollary \ref{art5.cor1} that $U \, k(z_S) = k(z_T) \, U$. So we get that $U \, h(S) = h(T) \, U$ 
which implies that $U \, g(S) = g(T) \, U$.

\medskip

The proposition follows by the previous lemma.
\end{demo}

\medskip\medskip

\begin{corollary}
Consider  Hilbert \cst-modules $E$,$F$ over a \cst-algebra $A$ and a unitary operator $U$ from $E$ into $F$. Let $S$ be  a normal element 
in $\cR(E)$ and $T$ a normal element in $\cR(F)$ such that $U \, S \subseteq T \, U  $.
Consider a subset $G$ of $\,\C$ which is compatible with both $S$ and $T$. Then we have for every $f \in \text{C}(G)$ that $U \, f(S) = f(T) \, U$.
\end{corollary}
\begin{demo}
By the previous proposition, we know that $U \, f(S) \subseteq f(T) \, U$.

This same proposition implies also that $U \,\overline{f}(S) \subseteq \overline{f}(T) \, U$. So, taking the adjoint of this last inclusion, we get that
$U^* \, f(T) \subseteq f(S) \, U^*$. In turn, this implies that
$f(T) \, U \subseteq U \, f(S)$. 

The corollary follows.
\end{demo}

\begin{corollary}
Consider a Hilbert \cst-module  $E$ over a \cst-algebra $A$ and let $S$ and $T$ be  normal elements in $\cR(E)$. If $S \subseteq T$, then $S = T$.
\end{corollary}

\bigskip

\section{Natural left and right multipliers of a regular operator}

In this section, we will show that any regular operator has enough well-behaved left and right multipliers. 

\medskip

The functional calculus for normal operators gurantees immediately the existence of well behaved left and right multipliers. Recall the 
following result from section \ref{art1}.

\begin{result} \label{art6.res1}
Consider a Hilbert \cst-module $E$ over a \cst-algebra $A$ and a normal element $T \in \cR(E)$.
Let $G$ be a subset of $\,\C$ which is compatible with $T$ and let $f$ be an element in $\text{C}(G)$. Consider moreover an element 
$g \in \text{K}(G)$. Then $g(T)$ is a left and right multiplier of $f(T)$ and 
$$g(T) \, f(T) \subseteq f(T) \, g(T)  = (f g)(T)$$
\end{result}

\medskip

\begin{result} \label{art6.res2}
Consider a Hilbert \cst-module $E$ over a \cst-algebra $A$ and a normal element $T \in \cR(E)$.
Let $G$ be a subset of $\,\C$ which is compatible with $T$ and let $f$ be an element in $\text{C}(G)$. Consider moreover a bounded 
net $(e_i)_{i \in I}$ in $\text{K}(G)$ which converges strictly to 1. Then the following properties hold.
\begin{itemize}
\item We have for every $v \in D(f(T))$ that $\bigl( \,f(T) \, e_i(T)\,v \,)_{i \in I}$ converges to $f(T) v$.
\item Consider $v \in E$. Then $v$ belongs to $D(f(T))$ $\Leftrightarrow$
The net $\bigl(\,f(T) \, e_i(T) \,v \,\bigr)_{i \in I}$ is convergent.
\end{itemize}
\end{result}

\medskip

In this section, we will prove an abnormal version of the previous result.
We will at the same time be a little bit more general.

\bigskip

\begin{lemma} \label{art6.lem1}
Consider Hilbert \cst-modules $E$,$F$,$G$ over a \cst-algebra $A$. Then we have the following properties :
\begin{itemize}
\item Consider $T \in \cR(F,G)$ and $x \in \cL(E,F)$. \newline Then $x$ is a right multiplier of $T$ $\Leftrightarrow$ $x$ 
is a right multiplier of $|T|$.
\item Consider $T \in \cR(E,F)$ and $x \in \cL(F,G)$. \newline Then $x$ is a left multiplier of $T$ $\Leftrightarrow$ $x$ 
is a left multiplier of $|T^*|$.
\end{itemize}
\end{lemma}

The first statement follows immediately from the fact that $D(T) = D(|T|)$ and result \ref{prel2.res1}. The second follows from 
the first by using corollary \ref{prel2.cor1}.

\medskip

\begin{lemma} \label{art6.lem2}
Consider Hilbert \cst-modules $E$,$F$ over a \cst-algebra $A$ and an element  $T \in \cR(E,F)$.
Let $(x_i)_{i \in I}$ be a net in $\cL(E)$ which converges strongly to 1 and such that
we have for every $i \in I$ that $x_i$ is a left and right multiplier of $|T|$ and
that $x_i \, |T| \subseteq |T| \, x_i$. Then we have the following properties.
\begin{enumerate}
\item We have for every $v \in D(T)$ that $(T  x_i v)_{i \in I}$ converges to $T(v)$.
\item Consider $v \in E$. Then $v$ belongs to $D(T)$ $\Leftrightarrow$ The net 
$(T  x_i v)_{i \in I}$ is convergent.
\end{enumerate}
\end{lemma}
\begin{demo}
Choose $v \in D(T)$. Then $v$ belongs also to $D(|T|)$.

We have moreover for every $i \in I$ that 
$$\|T(x_i\,v) - T(v) \| = \|T(x_i \, v - v)\| = \|\,|T|(x_i v - v)\|  = \|x_i\, |T|(v) - |T|(v) \|$$
which guarantees that the net $(T(x_i\, v))_{i \in I}$ converges to $T(v)$.

\medskip

The other implication of the second statement follows immediately from the closedness of $T$.
\end{demo}

\medskip

A similar results holds for right multipliers.

\begin{lemma}
Consider Hilbert \cst-modules $E$,$F$ over a \cst-algebra $A$ and an element  $T \in \cR(E,F)$.
Let $(x_i)_{i \in I}$ be a net in $\cL(E)$ which converges strongly$\,^*$ to 1 and such that
we have for every $i \in I$ that $x_i$ is a left and right multiplier of $|T^*|$ and that
$x_i \, |T^*| \subseteq |T^*| \, x_i$. Then we have the following properties.
\begin{enumerate}
\item We have for every $v \in D(T)$ that $( (x_i \comp T)\, v )_{i \in I}$ converges to $T(v)$.
\item Consider $v \in E$. Then $v$ belongs to $D(T)$ $\Leftrightarrow$ The net 
$((x_i \comp T) \, v)_{i \in I}$ is convergent.
\end{enumerate}
\end{lemma}
\begin{demo}
The first statement follows immediately because $(x_i \comp T) \, v = x_i \, T(v)$ for 
$i \in I$ and $v \in D(T)$.

\medskip

Take $v \in E$ and $w \in F$ such that $((x_i \comp T) \, v)_{i \in I}$ converges to $w$.

Choose $u \in D(T^*)$. We have for every $i \in I$ that 
that $x_i^*$ is a left and right multiplier of $|T^*|$ and that $x_i^* \, |T^*| \subseteq |T^*| \, x_i^*$. We have moreover 
for every $i \in I$ that
$$ \lan (x_i \comp T) \, v , u \ran = \lan v , (x_i \comp T)^* \,  u \ran
= \lan v , (T^*  x_i^*)\, u \ran $$

If $i \rightarrow \infty$, the left hand side of this equation converges to $\lan w , u \ran$, whilst the right hand side 
converges to $\lan v , T^* u \ran$ by the previous lemma.
So we get that $\lan w , u \ran = \lan v , T^* u \ran$.

Because $T^{**} = T$, this implies that $v \in D(T)$ and $T(v) = w$.
\end{demo}

\bigskip\medskip

The next lemma is a key result in this section.

\begin{lemma}
Consider Hilbert \cst-modules $E$,$F$ over a \cst-algebra $A$ and an element $T \in \cR(E,F)$.
Let $G$ be a subset of $\R^+$ which is compatible with both $T^* T$ and $T T^*$ and let $f$ be an element in $\text{C}(G)$. 

Consider $g \in \text{K}(G)$ and define $C$ as the closure of $g(T^* T) E$ in $E$, so $C$ is a sub Hilbert \cst-module of $E$.
Then we have that $C \subseteq D(T \, f(T^* T))  , D(f(T T^*) \, T)$ and that
$(T \, f(T^* T)) \, v = (f(T T^*) \, T) \, v$ for $v \in C$.
\end{lemma}
\begin{demo} 
By lemma \ref{art5.lem1}, we know that $C \subseteq D(|T|)$ and that $|T|_C$ belongs to $\cL(C)$. Because $D(T) = D(|T|)$, we 
have that $C \subseteq D(T)$.  We have moreover that 
$$\|T_C(v)\| = \|T(v)\| = \|\,|T|(v)\| \leq \|\,|T|_C\| \, \|v\|$$ for $v \in C$. So we see that $T_C$ is a bounded $A$-linear 
mapping from $C$ into $F$.

\vspace{1mm}

By lemma \ref{art5.lem1}, we know also that $C \subseteq D(T^*T)$. The same lemma implies also that $(T^* T)_C$ is a positive 
element in $\cL(C)$.

Because $T(T^* T) = (T T^*) T$, we have immediately that  $T_C \, (T^* T)_C = (T T^*) T_C$.

We have also that $G$ is compatible with $(T^* T)_C$ so the 
Fuglede-Putnam theorem (theorem \ref{art5.thm1}) guarantees that $T_C \, f((T^* T)_C) \subseteq f(T T^*) \, T_C$, so $T_C \, 
f((T^* T)_C) = f(T T^*) \, T_C$.

\vspace{1mm}

Again, lemma \ref{art5.lem1} gives us that $C \subseteq D(f(T^* T))$, that $f(T^* T)_C$ belongs to $\cL(C)$ and that 
$f(T^* T)_C = f((T^* T)_C)$. Consequently, $T_C \, f(T^* T)_C = f(T T^*) \, T_C$. 
The lemma follows.
\end{demo}

\bigskip

So we have in particular the following result.

\begin{result} \label{art6.res3}
Consider Hilbert \cst-modules $E$,$F$ over a \cst-algebra $A$ and an element $T \in \cR(E,F)$.
Let $G$ be a subset of $\R^+$ which is compatible with both $T^* T$ and $T T^*$ and let $f$ be an element in $\text{C}(G)$. 
Consider moreover $g \in \text{K}(G)$ and $v \in E$. Then we have that $g(T^* T) v$ belongs to $D(T \, f(T^* T)) \cap 
D(f(T T^*) \, T)$, that $(f g)(T^* T)$ belongs to $D(T)$ and that
$$(f(T T^*) \, T) \,(g(T^* T) v) = (T \, f(T^* T)) \,(g(T^* T) v) = T(\, (f g)(T^* T) v\,)$$
\end{result}

\medskip

A first application of this lemma concerns the abnormal version of result \ref{art6.res1}.

\begin{result} \label{art6.res4}
Consider Hilbert \cst-modules $E$,$F$ over a \cst-algebra $A$ and an element $T \in \cR(E,F)$.
Let $G$ be a subset of $\R^+$ which is compatible with both $T^* T$ and $T T^*$ and let $f$ be an element in $\text{K}(G)$. 
Then the element $f(T^* T)$ is a right multiplier of $T$, the element $f(T T^*)$ is a left multiplier of $T$ and
$$f(T T^*) \, T \subseteq T \, f(T^* T) $$
\end{result}
\begin{demo}
By result \ref{art6.res1}, we know that $f(T^* T)$ is a right multiplier of $|T|$, so lemma \ref{art6.lem1} implies that 
$f(T^* T)$ is a right multiplier of $T$. We get in a similar way that $f(T T^*)$ is a left multiplier of $T$.

Moreover, result \ref{art6.res3} implies for every $g \in \text{K}(G)$ and $v \in E$ that $g(T^* T) v$ belongs to 
$D(T \, f(T^* T)) \cap D(f(T T^*) \, T)$ and that
$$(T \, f(T^* T)) \, g(T^* T) v =  (f(T T^*) \, T) \, g(T^* T) v  $$
Now the boundedness of $f(T T^*) \, T$ and $T \, f(T^* T)$ implies easily the above inclusion.
\end{demo}

\medskip

Refering to lemma \ref{art6.lem2}, we get also immediately the following result.

\begin{result}
Consider Hilbert \cst-modules $E$,$F$ over a \cst-algebra $A$ and an element $T \in \cR(E,F)$.
Let $G$ be a subset of $\R^+$ which is compatible with both $T^* T$ and $T T^*$.  Consider moreover a bounded net 
$(e_i)_{i \in I}$ in $\text{K}(G)$ which converges strictly to 1. Then the following properties hold.
\begin{itemize}
\item We have for every $v \in D(T)$ that $\bigl(\,T \, e_i(T^* T) \,v\, \bigr)_{i \in I}$ converges to $T v$.
\item Consider $v \in E$. Then $v$ belongs to $D(T)$ $\Leftrightarrow$
The net $\bigl(\,T \, e_i(T^* T) \,v\, \bigr)_{i \in I}$ is convergent.
\end{itemize}
\end{result}

\bigskip

In the next part of this section, we want to prove a version of result \ref{art6.res4} where $f \in \text{C}(G)$. 
We will of course get unbounded operators so a little bit of care has to be taken.

\begin{proposition}
Consider Hilbert \cst-modules $E$,$F$ over a \cst-algebra $A$ and an element $T \in \cR(E,F)$.
Let $G$ be a subset of $\R^+$ which is compatible with both $T^* T$ and $T T^*$ and let $f$ be an element in $\text{C}(G)$. 
Then $T \, f(T^* T)$ and $f(T T^*) \, T$ are closable and we define
$T \comp f(T^* T)$ as the closure of $T \, f(T^* T)$ and we define $f(T T^*) \comp T$ as the closure of $f(T T^*) \, T$. 
Then we have that $T \comp f(T^* T) = f(T T^*) \comp T$.
\end{proposition}
\begin{demo}
By proposition \ref{art4.prop1}, it is clear that $D(T \, f(T^* T)) = D(|T| \, f(T^* T))$ and that
$$\|(T \, f(T^* T))(v)\| =\|(|T| \, f(T^* T))(v)\|$$ for $v \in 
D(T \, f(T^* T))$. Result \ref{art1.res2} implies moreover that $|T| \, f(T^* T)$ is closable, so we get also   
that $T \, f(T^* T)$ is closable. We denote the closure of 
$|T| \, f(T^* T)$ by $R$ and the closure of $T \, f(T^* T)$ by $S$.

Then the previous remarks imply that $D(R) = D(S)$ and that $\|R(v)\| = \|S(v)\|$ for $v \in D(S)$. \hspace{0.5cm} (*)

\medskip

Choose $v \in D(f(T T^*) \, T)$. Then $v \in D(T)$ and $T(v) \in D(f(T T^*))$.

Take a bounded net $(e_i)_{i \in I}$ in $\text{K}(G)$ such that $(e_i)_{i \in I}$ converges strictly to $1$. 

\vspace{1mm}

Fix $j \in I$. Using result \ref{art6.res3}, we get that $e_j(T^* T) v$ belongs to $D(S) \cap D(f(T^* T) \, T)$ and 
that \newline $S(e_j(T^* T)v) = (f(T T^*) \, T)(e_j(T^* T) v)$.

Result \ref{art6.res4}  gives us that $e_j(T T^*) \, T \subseteq T \, e_j(T^* T)$, so we get that
$T(e_j(T^* T) v) = e_j(T T^*) \, T(v)$. Because $T(v)$ belongs to $f(T T^*)$, this in turn implies that 
\begin{eqnarray*}
S(e_j(T^* T) v) & = & (f(T T^*) \, T)(e_j(T^* T) v) = f(T T^*)(e_j(T T^*)\, T(v))\\
&  = & e_j(T T^*)\bigl(f(T T^*)(T(v))\bigr) = e_j(T T^*)\,(f(T T^*) \, T)(v) \ .
\end{eqnarray*}
Hence, we see that $(e_i(T^* T) v)_{i \in I}$ converges to $v$ and that
$\bigl(S(e_i(T^* T) v)\bigr)_{i \in I}$ converges to $(f(T T^*) \, T)(v)$.

Therefore, the closedness of $S$ implies that $v \in D(S)$ and
$S(v) = (f(T T^*) \, T)(v)$.

\medskip

So we have proven that $f(T T^*) \, T \subseteq S$. This implies immediately that $f(T T^*) \, T$ is closable. 

Define the set $C = \lan \, g(T^* T) v \mid g \in \text{K}(G), v \in E \, \ran$.
Then result \ref{art1.res2} and lemma \ref{art1.lem3} imply that $C$ is a core for $R$. So (*) implies that $C$ is 
also a core for $S$. 
We know moreover that $C \subseteq D(f(T T^*) \, T)$ by result \ref{art6.res3}.
From this all, we conclude that the closure of $f(T T^*) \, T$ is equal to $S$.
\end{demo}

\medskip

Notice that we have also proven the following result :

\begin{result}
Consider Hilbert \cst-modules $E$,$F$ over a \cst-algebra $A$ and an element $T \in \cR(E,F)$.
Let $G$ be a subset of $\R^+$ which is compatible with both $T^* T$ and $T T^*$ and let $f$ be an element in $\text{C}(G)$. 
Define the set $C = \lan \, g(T^* T) v \mid g \in \text{K}(G), v \in E \, \ran$. Then $C \subseteq D(T \, f(T^* T)) \cap 
D(f(T T^*) \, T)$ and $C$ is a core for $T \comp f(T^* T)$.
\end{result}

So we have also that $D(T \, f(T^* T)) \cap D(f(T T^*) \, T)$ is a core for $T \comp f(T^* T)$.

\bigskip

\begin{result}
Consider Hilbert \cst-modules $E$,$F$ over a \cst-algebra $A$ and an element $T \in \cR(E,F)$.
Let $G$ be a subset of $\R^+$ which is compatible with both $T^* T$ and $T T^*$ and let $f$ be an element in $\text{C}(G)$. 
Then $T \comp f(T^* T) = f(T T^*) \comp T$ belongs to $\cR(E,F)$. 
\end{result}
\begin{demo}
By result \ref{art6.res3}, we have that $T \comp f(T^* T)$ is a densely defined closed $A$-linear operator from within $E$ 
into $F$. We have also immediately that 
$$\overline{f}(T^* T) \, T^* \subseteq (T \, f(T^* T))^* = (T \comp f(T^* T))^* $$
This implies by result \ref{art6.res3} that $(T \comp f(T^* T))^*$ is also densely defined.
We get moreover that
$$1 + \overline{f}(T^*T) \, T^* T \, f(T^* T)  \subseteq 1 + (T \comp f(T^* T))^* (T \comp f(T^* T))$$ 
By result \ref{art1.res2}, we know that  $\overline{f}(T^*T)\,T^* T\,f(T^* T)$ is closable and that its closure $P$ is equal 
to $(\overline{f} \, \io_G f)(T^* T)$.
So $P = (|f|^2 \, \io_G)(T^* T)$ which implies that $P$ is positive. Using  proposition \ref{art2.prop3}, this implies that \newline
$\bld (1 + P) = E$.

So we get that $1 + \overline{f}(T^*T)\,T^* T\,f(T^* T)$ has dense range, which in turn implies  that the element \newline 
$1 + (T \comp f(T^* T))^* (T \comp f(T^* T))$ has dense range.
Hence, we get by definition that $T \comp f(T^* T)$ belongs to $\cR(E,F)$.
\end{demo}

\bigskip

\begin{result}
Consider Hilbert \cst-modules $E$,$F$ over a \cst-algebra $A$ and an element $T \in \cR(E,F)$.
Let $G$ be a subset of $\R^+$ which is compatible with both $T^* T$ and $T T^*$ and let $f$ be an element in $\text{C}(G)$. 
Consider  $g \in \text{K}(G)$. Then we have the following properties.
\begin{itemize}
\item  The element $g(T^* T)$ is a right multiplier of $T \comp f(T^* T)$, $(f g)(T^* T)$ is a right multiplier of $T$. 
\item The element $g(T T^*)$ is a left multiplier of $T \comp f(T^* T)$, $(g f)(T T^*)$ is a right multiplier of $T$. 
\item We have the following chain of equalities
$$(T \comp f(T^* T)) \comp g(T^* T) = T \comp (f g)(T^* T) = 
(g f)(T T^*) \comp T = g(T T^*) \comp (T \comp f(T^* T))   $$
\end{itemize}
\end{result}
\begin{demo}
\begin{itemize}
\item The first statement and the first equality in the chain above follow easily from results \ref{art6.res1} and \ref{art6.res4}.
\item By result \ref{art6.res4}, we know that $(f g)(T T^*)$ is a right multiplier of $T$. We have moreover for every $v \in D(f(T T^*) \, T)$ that
$$g(T T^*)\, (f(T T^*) \, T)(v) = (g f)(T T^*) \, T(v) = ((g f)(T T^*) \comp T)(v) $$
This implies for every $v \in D(f(T T^*) \comp T)$ that
$g(T^* T)\, (f(T T^*) \comp T)(v) = ((g f)(T T^*) \comp T)(v)$.

By definition, this gives us that $g(T T^*)$ is a left multiplier of $f(T T^*) \comp T$
and that $g(T T^*) \comp (f(T T^*) \comp T)$ $= (g f)(T T^*) \comp T$.
\end{itemize}
Now result \ref{art6.res4} joins the two equalities.
\end{demo}

\medskip

This allows now to prove easily the following result.

\begin{proposition}
Consider Hilbert \cst-modules $E$,$F$ over a \cst-algebra $A$ and an element $T \in \cR(E,F)$.
Let $G$ be a subset of $\R^+$ which is compatible with both $T^* T$ and $T T^*$ and let $f$ be an element in $\text{C}(G)$.  
Then $$(T \comp f(T^* T))^* = (f(T T^*) \comp T)^* = \overline{f}(T^* T) \comp T^* = T^* \comp \overline{f}(T^* T) $$
\end{proposition}
\begin{demo} We have that $$\overline{f}(T^* T) \, T^* \subseteq (T \, f(T^* T))^* = (T \comp f(T^* T))^* $$ which implies that 
$\overline{f}(T^* T) \comp T^* \subseteq (T \comp f(T^* T))^*$.

\medskip

Choose $v \in D\bigl((T \comp f(T^* T))^*\bigr)$.

Take a bounded net $(e_i)_{i \in I}$ in $\text{K}(G)$ which converges strictly to 1.

Take $j \in I$. By the previous lemma, we know that $e_j(T T^*) v$ belongs to $D(T^* \comp \overline{f}(T T^*))$, 
that $(\overline{f} e_j)(T T^*)$ is a right multiplier of $T^*$ and that
$$(T^* \comp \overline{f}(T T^*))(e_j(T T^*) v) = (T^* \comp (\overline{f} e_j)(T T^*))(v) \hspace{2.8cm} (a)$$
Because $(\overline{f} e_j)(T T^*)$ is a right multiplier of $T^*$, we know that
$(f \overline{e_j})(T T^*)$ is a left multiplier of $T$ and that
$$T^* \comp (\overline{f} e_j)(T T^*) = ((f \overline{e_j})(T T^*) \comp T)^*
= (T \comp (f \overline{e_j})(T^* T))^* \hspace{1.5cm} \text{(b)}$$
Again, we know by the previous result that $\overline{e_j}(T^* T)$ is a right multiplier of $T \comp f(T^* T)$ and that
$$T \comp (f \overline{e_j})(T^* T) = (T \comp f(T^* T)) \comp \overline{e_j}(T^* T) $$
This implies that $e_j(T^* T)$ is a left multiplier of $(T \comp f(T^* T))^*$ and that 
$$e_j(T^* T) \comp (T \comp f(T^* T))^* = \bigl((T \comp f(T^* T)) \comp \overline{e_j}(T^* T)\bigr)^* = 
(T \comp (f \overline{e_j})(T^* T))^* =  T^* \comp (\overline{f} e_j)(T T^*) $$
where we used (b) in the last equality. So, using (a), we get that
$$(T^* \comp \overline{f}(T T^*))(e_j(T T^*) v) = (e_j(T^* T) \comp (T \comp f(T^* T))^*)(v) 
= e_j(T^* T) \ (T \comp f(T^* T))^*(v) $$

Hence, we see that $(e_i(T T^*) v)_{i \in I}$ converges to $v$ and that
$\bigl((T^* \comp \overline{f}(T T^*))(e_i(T T^*) v)\bigr)_{i \in I}$ converges to  \newline 
$(T \comp f(T^* T))^*(v)$.
Consequently, the closedness of $T^* \comp \overline{f}(T T^*)$ implies that
$v$ belongs to $D(T^* \comp \overline{f}(T T^*))$.
\end{demo}

\bigskip

We will use notation  \ref{art1.not1} in the next propositions.

\begin{proposition}
Consider Hilbert \cst-modules $E$,$F$ over a \cst-algebra $A$ and an element $T \in \cR(E,F)$.
Let $G$ be a subset of $\R^+$ which is compatible with both $T^* T$ and $T T^*$ and let $f$ be an element in $\text{C}(G)$.  
Then $$(T \comp f(T^* T))^* (T \comp f(T^* T))
= \bigl(T^* T\bigr) \comp \bigl(|f|^2(T^* T)\bigr)  \text{\ \  and \ \ } 
(T \comp f(T^* T))(T \comp f(T^* T))^* = \bigl(T T^*\bigr) \comp \bigl(|f|^2(T T^*)\bigr)$$
\end{proposition}
\begin{demo}
We have that
$$(T \comp f(T^* T))^* (T \comp f(T^* T)) = (f(T^* T) \comp T^*)(T \comp f(T T^*)) \supseteq
f(T^* T) T^* T f(T^* T)$$
By result \ref{art1.res2}, we know that $f(T^* T) T^* T f(T^* T)$ is closable and that the closure is equal to
$(\io_G \, |f|^2)(T^* T)$. Hence, the previous inclusion implies that
$$(T \comp f(T^* T))^* (T \comp f(T^* T)) \supseteq (\io_G \, |f|^2)(T^* T)$$
But both operators are selfadjoint, so the inclusion must be an equality.
Hence, result \ref{art1.res2} implies that
$$(T \comp f(T^* T))^* (T \comp f(T^* T)) = (\io_G \, |f|^2)(T^* T) = T^* T \comp |f|^2(T^* T)$$

Using the previous proposition, the second equality can be brought in the same form as the first one 
(with the role of $T$ and $T^*$ interchanged).
\end{demo}

\medskip

Using proposition \ref{art3.prop2} and result \ref{art1.res2}, this implies the following result.

\begin{corollary}
Consider Hilbert \cst-modules $E$,$F$ over a \cst-algebra $A$ and an element $T \in \cR(E,F)$.
Let $G$ be a subset of $\R^+$ which is compatible with both $T^* T$ and $T T^*$ and let $f$ be an element in $\text{C}(G)$.  
Then $$|T \comp f(T^* T)| = |T| \comp |f|(T^* T) \hspace{2cm} \text{ and }  \hspace{2cm}
|(T \comp f(T^* T))^*| = |T^*| \comp |f|(T T^*)$$
\end{corollary}

\bigskip\bigskip

Combining this result with lemma \ref{art6.lem1}, we get the following one.

\begin{result}
Consider Hilbert \cst-modules $E$,$F$ over a \cst-algebra $A$ and an element $T \in \cR(E,F)$.
Let $G$ be a subset of $\R^+$ which is compatible with both $T^* T$ and $T T^*$ and let $f$ be an element in $\text{C}(G)$. 
Consider a bounded net $(e_i)_{i \in I}$ in $\text{K}(G)$ which converges strictly to 1. Then the following holds.
\begin{enumerate}
\item We have for every $v \in D(T \comp f(T^* T))$ that $\bigl(\,(T \comp f(T^* T)) \, e_i(T^* T)\,v \,\bigr)_{i \in I}$ 
converges to $(T \comp f(T^* T))\,v$
\item Consider $v \in E$. Then $v$ belongs to $D(T \comp f(T^* T))$ $\Leftrightarrow$
The net $\bigl(\,(T \comp f(T^* T)) \, e_i(T^* T)\,v \,\bigr)_{i \in I}$ is convergent.
\end{enumerate}
\end{result}

\bigskip\bigskip

It is possible to get some further information about the domains by the following result.

\begin{result}
Consider Hilbert \cst-modules $E$,$F$ over a \cst-algebra $A$ and an element $T \in \cR(E,F)$.
Let $G$ be a subset of $\R^+$ which is compatible with both $T^* T$ and $T T^*$.
Then the following holds.
\begin{itemize}
\item Consider $f \in \text{C}(G)$. Then we have the following properties.
\begin{enumerate} 
\item $D(T \, f(T^* T)) = D(|T| \, f(T^* T)) = D(T \comp f(T^* T)) \cap D(f(T^* T))
= D(|T| \comp f(T^* T)) \cap D(f(T^* T))$
\item $D(T \comp f(T^* T)) = D(|T| \comp f(T^* T))$
\end{enumerate}
\item Consider $f,g \in \text{C}(G)$. Then we have for every $v \in D(T \comp f(T^* T))$ and
$w \in D(T \comp g(T^* T))$ that
$$\lan (T \comp f(T^* T))(v) , (T \comp g(T^* T))(w) \ran = \lan (|T| \comp f(T^* T))(v) ,
(|T| \comp g(T^* T))(w) \ran $$
\end{itemize}
\end{result}
\begin{demo}
\begin{itemize}
\item Choose $f \in \text{C}(G)$. Then proposition \ref{art4.prop1} implies immediately that
$D(T \, f(T^* T)) = D(|T| \, f(T^* T))$ and that $\|(T \, f(T^* T))(v)\| = \|(|T| \, f(T^* T))(v)\|$ for $v \in D(T \, f(T^* T))$. 
This implies immediately that
$D(T \comp f(T^* T)) = D(|T| \comp f(T^* T))$ and that
$\|(T \comp f(T^* T))(v)\| = \|(|T| \comp f(T^* T))(v)\|$ for $v \in D(T \comp f(T^* T))$.

Result \ref{art1.res3} gives us that $D(|T| \, f(T^* T)) = D(|T| \comp f(T^* T)) \cap D(f(T^* T))$.
\item This follows now easily from the results in the first part of the proof and proposition \ref{art4.prop1}.
\end{itemize}
\end{demo}

\medskip

\begin{remark} \rm
We give quickly an example how this can be useful.
Therefore suppose that $f$ is an element in $\text{C}(G)$ such that there exists a positive number $M$ satisfying $|f(c)| \leq M$ and 
$\sqrt{c} \, |f(c)| \leq M$ for $c \in \si(T^* T) \cap G$. Then $D(|T| \comp f(T^* T)) \cap D(f(T^* T)) = E$. 

So in this case $T \comp f(T^* T) = T \, f(T^* T) \in \cL(E,F)$.
\end{remark}

\bigskip \bigskip

A situation where the above terminology applies is the connection between $T$ and its    $z$-transform $z_T$. Consider  
Hilbert \cst-module $E$,$F$  over a \cst-algebra $A$ and $T \in \cR(E,F)$. Proposition \ref{art4.prop2} implies immediately the following results.
\begin{itemize}
\item $z_T = T \comp (1 + T^* T)^{-\frac{1}{2}} = (1 + T T^*)^{-\frac{1}{2}} \comp T $
\item $T = z_T \comp (1 - z_T^* z_T)^{-\frac{1}{2}} = (1 - z_T z_T^*)^{-\frac{1}{2}} \comp z_T$
\end{itemize}

\bigskip

\section{\cst-algebras of adjointable operators on Hilbert \cst-modules} \label{art7}

In this section, we will fix a Hilbert \cst-module $E$ over a \cst-algebra. At the same time, we will consider a non-degenerate 
sub-\cst-algebra $B$ of $\cL(E)$. We will look at an embedding of $\cR(B)$ into $\cR(E)$.

\medskip

Concerning the multiplier algebra, we have the following well known result:
$$M(B) = \{\, x \in \cL(E) \mid \text{ We have for every } b \in B \text{ that }
x b \text{ and } b x \text{ belong to } B \, \}$$

\medskip

As pointed out in \cite{Wor6} for Hilbert spaces, we can also embed $\cR(B)$ in $\cR(E)$.

\begin{definition}
Call $\pi$ the inclusion of $B$ into $\cL(E)$, then $\pi$ is a non-degenerate
$^*$-homomorphism from $B$ into $\cL(E)$. Let $T$ be an element affiliated to $B$. Then we
define $\tilde{T} = \pi(T)$, so $\tilde{T}$ is a regular operator in $E$.
\end{definition}

\medskip

Because $\pi$ is injective, we know immediately that the mapping $\cR(B) \rightarrow
\cR(E) : T \mapsto \tilde{T}$
 is injective.

\medskip

We have also immediately that $\tilde{x} = x$ for every $x \in M(B)$.

\bigskip

So we get the following determining property of the embedding :

\begin{result}
We  have  that $z_T = z_{\tilde{T}}$ for every element $T$ affiliated with $B$.
\end{result}

\medskip

Because this embedding arises from an injective non-degenerate $^*$-homomorphism, most of
the properties of this embedding follow from the theory of injective non-degenerate $^*$-homomorphisms.  We will now prove some extra results.

\bigskip

Looking at theorem \ref{prel1.thm1}, we have immediately the following result :

\begin{result} \label{art7.res1}
Consider an element $T$ affiliated to $B$. Then
\begin{itemize}
\item We have for every $x \in D(T)$ and $v \in E$ that $x v \in D(\tilde{T})$ and $\tilde{T}(x v) = T(x) v$.
\item $D(T) E$ is a core for $\tilde{T}$
\end{itemize}
\end{result}

Let $D$ be a core for $T$ and $K$ a dense subspace of $E$. Then $D K$ is a core for
$\tilde{T}$.

\bigskip

We have also some kind of converse of this result.

\begin{proposition} \label{art7.prop1}
Consider an element $T$ affiliated with $B$ and an element $x$ in $B$. Then
\begin{itemize}
\item If $x$ belongs to $D(T)$, then $x$ is a right multiplier of $\tilde{T}$ and
$\tilde{T}\comp x = T(x)$.
\item We have that $x$ belongs to $D(T)$ $\Leftrightarrow$ $x$ is a right multiplier of
$\tilde{T}$ and $\tilde{T} \comp x$ belongs to $B$.
\end{itemize}
\end{proposition}
\begin{demo}
The first statement follows immediately from the previous result.

Choose $x \in B$ such that $x$ is a right multiplier of $\tilde{T}$ and 
$\tilde{T} \comp x$ belongs to $B$.

Take $y \in D(T^*)$. Choose $v,w \in E$. 

Because $x$ is a right multiplier of $\tilde{T}$, we get that $x w$ belongs to $D(\tilde{T})$ and $\tilde{T}(x w) = 
(\tilde{T} \comp x) \, w$. 
By results \ref{prel1.res1} and \ref{art7.res1}, we have also that $y v$ belongs to $D(\tilde{T}^*)$ and that
$\tilde{T}^*(y v) = T^*(y) \, v$. Hence,
$$\lan x^* \, T^*(y) \, v , w \ran  =  \lan T^*(y) \, v , x w \ran
= \lan \tilde{T}^* (y v) , x w \ran  
= \lan y v , \tilde{T} (x w)  \ran 
= \lan y v , (\tilde{T} \comp x) \, w \ran
= \lan (\tilde{T} \comp x)^* y v , w \ran
$$
Consequently, $x^* \, T^*(y) = (\tilde{T} \comp x)^* y$.

Because $T = T^{**}$, this implies that $x \in D(T)$ and $T(x) = \tilde{T} \comp x$.
\end{demo}

\medskip

\begin{corollary}
Consider an element $T$ affiliated with $B$ and $x$ an element in $M(B)$. Then
\begin{itemize}
\item If $x$ is a right multiplier of $T$, then $x$ is a right multiplier of $\tilde{T}$ and
$\tilde{T} \comp x = T \comp x$.
\item We have that $x$ is a right multiplier of $T$ $\Leftrightarrow$ $x$ is a right
multiplier of $\tilde{T}$ and $\tilde{T} \comp x$ belongs to $M(B)$.
\end{itemize}
\end{corollary}
\begin{demo}
\begin{itemize}
\item Suppose that $x$ is a right multiplier of $T$.

Choose $b \in B$ and $v \in E$. Then $x b$ belongs to $D(T)$, so the previous proposition
implies that $x b$ is right multiplier of $\tilde{T}$ and $\tilde{T} \comp (x b) = T(x b) = 
(T \comp x) \, b$. So we get that $x b v$ belongs to $D(\tilde{T})$ and $\tilde{T}(x b v) = (T \comp x) \, b v$.

Because $B E$ is dense in $E$ and $\tilde{T}$ is closed, this implies easily for every $w
\in E$ that $x w$ belongs to $\tilde{T}$ and that $\tilde{T}(x w) = (T \comp x) \, w$.
This implies that $x$ is a right multiplier of $\tilde{T}$ and $\tilde{T} \comp x = T \comp x$.

We have in particular that $\tilde{T} \comp x$ belongs to $M(B)$.

\item Suppose that $x$ is right multiplier of $\tilde{T}$ and that $\tilde{T} \comp x$ belongs to $M(B)$.

Choose $b \in B$. Then $x b$ is a right multiplier of $\tilde{T}$ and $\tilde{T} \comp (x b) =
(\tilde{T} \comp x) \, b$ which clearly belongs to $B$. Hence, the previous proposition implies that $x b \in D(T)$ and 
$T(x b) = \tilde{T} \comp (x b) = (\tilde{T}\comp x)\, b$.

So we get that $x$ is a right multiplier of $T$ and $T\comp x = \tilde{T}\comp x$.
\end{itemize}
\end{demo}

\medskip

By using this result and corollary \ref{prel2.cor1} , we have also the following result.

\begin{corollary}
Consider an element $T$ affiliated with $B$ and an element $x$ in $M(B)$. Then
\begin{itemize}
\item If $x$ is a left multiplier of $T$, then $x$ is a left multiplier of $\tilde{T}$ and
$x \comp \tilde{T} = x \comp T$.
\item We have that $x$ is a left multiplier of $T$ $\Leftrightarrow$ $x$ is a left
multiplier of $\tilde{T}$ and $x \comp \tilde{T}$ belongs to $M(B)$.
\end{itemize}
\end{corollary}

\bigskip

Looking at example 4 of \cite{Wor6}, we have also the following result.

\begin{result}
Consider a regular operator $T$ in $E$. Then there exists an element $S$ affiliated with
$B$ such 

\vspace{1mm}
that $\tilde{S} = T$ \ \ \ $\Leftrightarrow$
\begin{minipage}[t]{10cm}
\begin{enumerate}
\item $z_T$ belongs to $M(B)$
\item $(1 - z_T^* z_T)^\frac{1}{2} B$ is dense in $B$
\end{enumerate}
\end{minipage}
\end{result}

\medskip

If there exists such an $S$, we have immediately that $z_T = z_S$, so $z_T$ will certainly
satisfy the two mentioned conditions.

If $z_T$ satisfies these two conditions, there exists an element $S$ affiliated with $B$
such that $z_S =z_T$. So we have that $z_{\tilde{S}} = z_S = z_T$ which implies that
$\tilde{S} = T$.

\bigskip

This implies immediately the following result.

\begin{proposition}
Consider a Hilbert \cst-module over $E$ over a \cst-algebra $A$. Then the mapping
$\cR(\cK(E)) \rightarrow \cR(E) : T  \mapsto \tilde{T}$ is a bijection.
\end{proposition}

\medskip

Remember that even for injective non-degenerate $^*$-homomorphisms, the invertibility of an element and its image where not equivalent. 
It is however true for the bijection between $\cR(\cK(E))$ and $\cR(E)$.

\vspace{2mm}

\begin{lemma}
Consider a Hilbert \cst-module over $E$ over a \cst-algebra $A$. Let $T$ be an element
affiliated with $\cK(E)$ and $v \in D(\tilde{T})$. Then we have for every $w \in E$ that
$\th_{v,w}$ belongs to $D(T)$ and $T(\th_{v,w}) = \th_{\tilde{T}(v),w}$.
\end{lemma}
\begin{demo}
Take $w \in E$. Then we have for every $u \in E$ that $\th_{v,w}(u) = v \, \lan u , w \ran$, which implies that $\th_{v,w}(u)$ belongs 
to $D(\tilde{T})$ and $\tilde{T}(\th_{v,w}(u)) = \tilde{T}(v) \, \lan u , w \ran = \th_{\tilde{T}(v),w}(u)$.

So we see that $\th_{v,w}$ is a right multiplier of $\tilde{T}$ and $\tilde{T} \comp \th_{v,w} =
\th_{\tilde{T}(v),w}$. By proposition \ref{art7.prop1}, this implies that $\th_{v,w}$ belongs to $D(T)$ and
$T(\th_{v,w}) = \th_{\tilde{T}(v),w}$.
\end{demo}

\medskip

So we get easily the following result.

\begin{proposition}
Consider a Hilbert \cst-module over $E$ over a \cst-algebra $A$ and $T$ an element affiliated with $\cK(E)$. Then $T$ is invertible 
$\Leftrightarrow$ $\tilde{T}$ is invertible.
\end{proposition}

\medskip

\begin{corollary} \label{art7.cor1}
Consider a Hilbert \cst-module over $E$ over a \cst-algebra $A$ and $T$ a normal element affiliated with $\cK(E)$. 
Let $G$ be an almost closed subset of $\,\C$. Then $G$ is compatible with $T$ $\Leftrightarrow$ $G$ is compatible with $\tilde{T}$.
\end{corollary}

\bigskip

\section{Representing Hilbert \cst-modules on Hilbert spaces} \label{art8}

Sometimes it is useful to represent Hilbert \cst-modules on Hilbert spaces. In this section, we describe a natural way to do so 
and show that this procedure is well behaved.

\medskip

\begin{notation}
Consider a Hilbert \cst-module $E$ over a \cst-algebra $A$ and $\om \in A^*_+$. 
As in the case of positive functionals on a \cst-algebra, there is a natural way to construct a Hilbert space $E_\om$ together with 
a bounded linear map $E \rightarrow E_\om : v \mapsto \overline{v}$ such that
\begin{enumerate}
\item $\overline{E}$ is dense in $E_\om$.
\item We have for $v,w \in E$ that $\lan \overline{v} , \overline{w} \ran = \om(\lan v, w \ran)$.
\end{enumerate}
It is clear that $\|\overline{v}\| \leq \|\om\| \, \|v\|$ for every $v \in E$.
\end{notation}

\medskip

In a next step, we want to represent operators on $E$ by operators on $E_\om$.

\begin{notation}
Consider Hilbert \cst-modules $E$,$F$ over a \cst-algebra $A$ and $\om \in A^*_+$.

Let $T$ be a linear map  from $E$ into $F$ such that there exists a positive number $M$ such that $\lan T v , T v \ran \leq 
M \, \lan v , v \ran$ for every $v \in H$.
Then there exists a unique bounded linear map $T_\om$ from $E_\om$ into $F_\om$ such that
$T_\om \, \overline{v} = \overline{T v}$ for $v \in E$. It is clear that $\|T_\om\| \leq M^\frac{1}{2}$
\end{notation}

This definition is in particular applicable if $T$ is an bounded $A$-linear map $T$ (e.g. if 
$T \in \cL(E,F)$) in which case $\lan T v , T v \ran \leq \|T\|^2 \, \lan v, v \ran$ for every $v \in E$. We have in this case 
that $\|T_\om\| \leq \|T\|$.

For the rest of this section, we will only be interested in the case that $T$ belongs to $\cL(E,F)$.

\bigskip 

The following properties are straightforward to prove.

\begin{result}
Consider Hilbert \cst-modules $E$,$F$ over a \cst-algebra $A$ and $\om \in A^*_+$.
Then 
\begin{enumerate}
\item The mapping $\cL(E,F) \rightarrow \cL(E_\om,F_\om)$ is linear 
\item We have for every $T \in \cL(E,F)$ that $(T_\om)^* = (T^*)_\om$
\item The mapping $\cL(E,F) \rightarrow \cL(E_\om,F_\om)$ is bounded and has norm $\leq 1$.
\item The mapping $\cL(E,F) \rightarrow \cL(E_\om,F_\om)$ is strongly continuous on bounded sets.
\item The mapping $\cL(E,F) \rightarrow \cL(E_\om,F_\om)$ is strongly$\,^*$ continuous on bounded sets.
\end{enumerate}
\end{result}

\begin{result}
Consider Hilbert \cst-modules $E$,$F$,$G$ over a \cst-algebra $A$ and $\om \in A^*_+$.
Then we have for every $S \in \cL(E,F)$, $T \in \cL(F,G)$ that $(T S)_\om = T_\om \, S_\om$.
\end{result}

\medskip

Let us now consider the case where $T$ is regular operator.

\begin{definition}
Consider Hilbert \cst-modules $E$,$F$ over a \cst-algebra $A$ and $\om \in A^*_+$.
Let $T$ be an element in $\cR(E,F)$. Then there exists a unique 
densely defined closed linear map $T_\om$ from within $E_\om$ into $F_\om$ such that
$\overline{D(T)}$ is a core for $T_\om$ and $T_\om  \, \overline{v} = \overline{T v}$ for $v \in E$.
\end{definition}

\medskip

It is easy to see that $\lan \, \overline{T v} , \overline{w} \, \ran = \lan \, \overline{v} , \overline{T^* w} \, \ran$ 
for $v \in D(T)$, $w \in D(T^*)$. Because $\overline{D(T^*)}$ is dense in $F_\om$, this equality implies easily the existence 
of the mapping $T_\om$ above.

\medskip

This definition implies also easily that $\overline{D}$ is a core for $T_\om$ if $D$ is a core for $D(T)$.

\bigskip

In a following proposition, we prove another expected characterization of $T_\om$.

\begin{proposition} \label{art8.prop1}
Consider Hilbert \cst-modules $E$,$F$ over a \cst-algebra $A$ and $\om \in A^*_+$.
Let $T$ be an element in $\cR(E,F)$. Then $(z_T)_\om = z_{T_\om}$.
\end{proposition}
\begin{demo}
Because $\|z_T\| \leq 1$, we know that $\|(z_T)_\om\| \leq 1$.
We have moreover that $1 - (z_T)_\om^* (z_T)_\om = (1 - z_T^* z_T)_\om$.
Because the mapping $\cL(E) \rightarrow \cB(E_\om) : S \rightarrow S_\om$ is a $^*$-homomorphism, this implies that 
$((1-z_T^* z_T)_\om)^\frac{1}{2}
= \bigl( (1-z_T^* z_T)^\frac{1}{2} \bigr)_\om$, 
so we see that $(1 - (z_T)_\om^* (z_T)_\om)^\frac{1}{2} =  \bigl( (1-z_T^* z_T)^\frac{1}{2} 
\bigr)_\om$. \hspace{0.5cm} (*)

\vspace{1mm}

Hence, we get that $(1 - (z_T)_\om^* (z_T)_\om)^\frac{1}{2} \, \overline{E}
= \bigl( (1 - z_T^* z_T)^\frac{1}{2} E \bigr)\streep$.
Because $(1-z_T^* z_T)^\frac{1}{2} E$ is dense in $E$, this implies that
$(1 - (z_T)_\om^* (z_T)_\om)^\frac{1}{2} E_\om$ is dense in $E_\om$.

So there exists a densely defined closed linear operator $S$ from within $E_\om$ into $F_\om$ such that $z_S = (z_T)_\om$. 

Because $\overline{E}$ is dense in $E_\om$, we have that
$(1 - z_S^* z_S)^\frac{1}{2} \, \overline{E}$ is a core for $S$. Now (*)  implies that
$$(1 - z_S^* z_S)^\frac{1}{2} \, \overline{E} =  \bigl( (1-z_T^* z_T)^\frac{1}{2} \bigr)_\om \, \overline{E} = 
\bigl( (1 - z_T^* z_T)^\frac{1}{2} E \bigr)\streep = \overline{D(T)} $$
which implies that $\overline{D(T)}$ is a core for $S$.

\vspace{1mm}

Choose $v \in D(T)$. Then there exist $w \in E$ such that
$v = (1-z_T^* z_T)^\frac{1}{2} w$. So $T v = z_T w$.

By (*), we have then also that $\overline{v} = \bigl((1-z_T^* z_T)^\frac{1}{2}\bigr)_\om \, \overline{w} = 
(1 - (z_T)_\om^* (z_T)_\om)^\frac{1}{2} \, \overline{w} =(1 - z_S^* z_S)^\frac{1}{2} \, \overline{w}$.

This implies that $S \, \overline{v} = z_S \, \overline{w} = (z_T)_\om \, \overline{w} = 
\overline{z_T w}= \overline{T v} = T_\om \, \overline{v}$.

Therefore we get that $S = T_\om$. So $(z_T)_\om = z_S = z_{T_\om}$.
\end{demo}

\bigskip

Because of this property, the following two results are very easy to prove.

\begin{result} \label{art8.res2}
Consider Hilbert \cst-modules $E$,$F$ over a \cst-algebra $A$ and $\om \in A^*_+$.
Let $T$ be an element in $\cR(E,F)$. Then we have that $(T_\om)^* = (T^*)_\om$.
\end{result}

\medskip

\begin{result}
Consider a Hilbert \cst-module $E$ over a \cst-algebra $A$ and $\om \in A^*_+$.
Let $T$ be an element in $\cR(E)$.
\begin{enumerate}
\item If $T$ is normal, then $T_\om$ is normal.
\item If $T$ is selfadjoint, then $T_\om$ is selfadjoint.
\item If $T$ is positive, then $T_\om$ is positive.
\item If $T$ is strictly positive, then $T_\om$ is strictly positive.
\end{enumerate}
\end{result}

\bigskip

There is another useful way to get to $T_\om$ for $T \in \cR(E)$, using the bijection 
$\cR(\cK(E)) \rightarrow \cR(E) : S \mapsto \tilde{S}$.

\begin{result} \label{art8.res1}
Consider a Hilbert \cst-module $E$ over a \cst-algebra $A$ and $\om \in A^*_+$.
Define the function $\pi$ from $\cK(E)$ into $\cK(E_\om)$ such that
$\pi(x) = x_\om$ for every $x \in \cK(E)$. Then $\pi$ is a non-degenerate $^*$-homomorphism of
$\cK(E)$ into $\cK(E_\om)$ such that $\pi(x) = x_\om$ for $x \in \cL(E)$.

Let $T$ be an element in $\cR(E)$. So there exists a unique element $S \in \cR(\cK(E))$ such that $\tilde{S} = T$. 
Then we have that $\pi(S) = T_\om$.
\end{result}

This follows immediately from the fact that $z_{\pi(S)} = \pi(z_S) = (z_S)_\om = (z_T)_\om
= z_{T_\om}$.

\bigskip

We can use this result to prove easily the following one.

\begin{proposition} \label{art8.prop2}
Consider a Hilbert \cst-module $E$ over a \cst-algebra $A$ and $\om \in A^*_+$.
Let $T$ be a normal element in $\cR(E)$. Let $G$ be a subset of $\,\C$ which is compatible with $T$. Then we have the 
following properties.
\begin{itemize}
\item The set $G$ is compatible with $T_\om$.
\item We have for every $f \in \text{C}(G)$ that $f(T_\om) = f(T)_\om$.
\end{itemize}
\end{proposition}
\begin{demo}
Define the non-degenerate $^*$-homomorphism $\pi$ from $\cK(E)$ into $\cK(E_\om)$ such that
$\pi(x) = x_\om$ for every $x \in \cK(E)$.

Take $S \in \cR(\cK(E))$ such that $\tilde{S} = T$. Then the previous result implies that $\pi(S) = T_\om$.

By the remark after proposition \ref{prel1.prop3} and corollary \ref{art7.cor1}, we know that $S$ is normal and that 
$G$ is compatible with $S$. This implies that $G$ is compatible with $\pi(S)$. So $G$ is compatible with $T_\om$.

Take $f \in \text{C}(G)$. By proposition \ref{art3.prop1}, we have that $\widetilde{f(S)} = f(\tilde{S}) = f(T)$. 
So the previous result implies that $f(T)_\om = \pi(f(S))$.
Hence,  $f(T)_\om = \pi(f(S)) = f(\pi(S)) = f(T_\om)$.
\end{demo}

\bigskip

Using result \ref{art8.res1}, also the next results are easy consequences of already known results.

\begin{result} \label{art8.res3}
Consider Hilbert \cst-modules $E$,$F$ over a \cst-algebra $A$ and $\om \in A^*_+$.
Let $T$ be an element in $\cR(E,F)$. Then we have that
$(T_\om)^* (T_\om) = (T^* T)_\om$ and $(T_\om ) (T_\om)^* = (T T^*)_\om$
\end{result}

\vspace{2mm}

\begin{result}
Consider Hilbert \cst-modules $E$,$F$ over a \cst-algebra $A$ and $\om \in A^*_+$.  
Let $T$ be an element in $\cR(E,F)$. Then
\begin{itemize}
\item If $T$ is invertible, then $T_\om$ is invertible and $(T_\om)^{-1} = (T^{-1})_\om$.
\item If $T$ is adjointable invertible, then $T_\om$ is bounded invertible.
\end{itemize}
\end{result}

\vspace{2mm}

\begin{corollary}
Consider a Hilbert \cst-module $E$ over a \cst-algebra $A$ and $\om \in A^*_+$.
Let $T$ be an element in $\cR(E)$. Then $\si(T_\om) \subseteq \si(T)$.
\end{corollary}

\bigskip

Proposition \ref{art8.prop2} implies immediately the following result.

\begin{result}
Consider a Hilbert \cst-module $E$ over a \cst-algebra $A$ and $\om \in A^*_+$.
Let $T$ be a normal element $\cR(E)$. Then we have the following properties.
\begin{itemize}
\item We have for every $n \in \N \cup \{0\}$ that $(T_\om)^n = (T^n)_\om$.
\item If $T$ is invertible, then we have for every $n \in \Z$ that $(T_\om)^n = (T^n)_\om$.
\item If $T$ is positive, then we have for every $s \in \R^+$ that $(T_\om)^s = (T^s)_\om$.
\item If $T$ is strictly positive, then we have for every $z \in \C$ that
$(T_\om)^z = (T^z)_\om$.
\end{itemize}
\end{result}

\bigskip

The proof of the following lemma is straightforward.

\begin{lemma} \label{art8.lem1}
Consider two Hilbert \cst-modules $E$,$F$ over a \cst-algebra $A$ and let $K$ be a subset of $A^*_+$ which is separating for $A$.
Consider elements $S,T$ in $\cL(E,F)$. Then $S = T$ $\Leftrightarrow$ We have for every $\om \in K$ that $S_\om = T_\om$.
\end{lemma}

\medskip

Combining this lemma with proposition \ref{art8.prop1}, we get immediately the following result.

\begin{proposition}
Consider two Hilbert \cst-modules $E$,$F$ over a \cst-algebra $A$ and let $K$ be a subset of $A^*_+$ which is separating for $A$. 
Consider elements $S,T$ in $\cR(E,F)$. Then $S = T$ $\Leftrightarrow$ We have for every $\om \in K$ that $S_\om = T_\om$.
\end{proposition}

\bigskip

Combining this proposition with results \ref{art8.res2} and \ref{art8.res3}, we get easily the following result.

\begin{corollary}
Consider a Hilbert \cst-module $E$ over a \cst-algebra $A$ and let $K$ be a subset of $A^*_+$ whih is separating for $A$. Consider 
an element $T$ in $\cR(E)$. Then we have the following properties.
\begin{enumerate}
\item $T$ is normal $\Leftrightarrow$ We have for every $\om \in K$ that $T_\om$ is normal.
\item $T$ is selfadjoint $\Leftrightarrow$ We have for every $\om \in K$ that $T_\om$ is selfadjoint.
\item $T$ is positive $\Leftrightarrow$ We have for every $\om \in K$ that $T_\om$ is positive.
\end{enumerate}
\end{corollary}

\medskip

This proposition can then also be used to prove easily the following result.

\begin{proposition}
Consider a Hilbert \cst-module $E$ over a \cst-algebra $A$ and let $S$ and $T$ be strictly positive elements in $\cR(E)$. 
Then $S = T$ $\Leftrightarrow$ We have for every $t \in \R$ that
$S^{it} = T^{it}$.
\end{proposition}

\medskip

In \cite{Kus}, we will use this Hilbert space procedure in the proof of proposition \ref{art9.prop1}. This was in fact the main 
reason to include this section.

\bigskip

\section{Commuting normal operators}

In this section, we look at an obvious notion of commuting normal regular operators and introduce the basic properties.

\medskip

Let us start with a familiar definition of strongly commuting operators.

\begin{definition}
Consider a Hilbert \cst-module $E$ over a \cst-algebra $A$ and let $S$,$T$ be normal elements in $\cR(E)$. Then we say that $S$ 
(strongly) commutes with $T$ $\Leftrightarrow$ We have for every $f,g \in \text{C}_0(\C)$ that $f(S) g(T) = g(T) f(S)$.
\end{definition}

It is clear that we get a symmetric definition in this way.

\medskip

\begin{result} \label{art9.res1}
Consider a Hilbert \cst-module $E$ over a \cst-algebra $A$ and let $S$,$T$ be normal elements in $\cR(E)$. Let $F$,$G$ be subsets 
of  $\,\C$ such that $F$ is compatible with $S$ and $G$ is compatible with $T$. Then we have that $S$ commutes with $T$ $\Leftrightarrow$
We have for every $f \in \text{C}_0(F)$ and $g \in \text{C}_0(G)$ that $f(S) g(T) = g(T) f(S)$.
\end{result}
\begin{demo}
\begin{itemize}
\item Suppose that $S$ commutes with $T$. Take $f \in \text{C}_0(F)$ and $g \in \text{C}_0(G)$. By lemma \ref{art3.lem1}, there exist 
$k, h \in \text{C}_0(\C)$ such that $f \subseteq h$ and
$g \subseteq k$. This implies that $f(S) g(T) = k(S) h(T) = h(T) k(S) = g(T) f(S)$.
\item Suppose that $f(S) g(T) = g(T) f(S)$ for every $f \in \text{C}_0(F)$ and $g \in \text{C}_0(G)$. Applying lemma \ref{art5.lem4} 
twice, this implies that $f(S) g(T) = g(T) f(S)$ for every $f \in \text{C}_b(F)$ and $g \in \text{C}_b(G)$.

\vspace{1mm}

So we get for every $h,k \in \text{C}_0(\C)$ that $h(S) k(T) = (h_F)(S) (k_G)(T) = (k_G)(T) (h_F)(S) = k(T) h(T)$.
\end{itemize}
\end{demo}

\medskip

Using lemma \ref{art5.lem4}, we get immediately the following results.

\begin{result} \label{art9.res2}
Consider a Hilbert \cst-module $E$ over a \cst-algebra $A$ and let $S$,$T$ be commuting  normal elements in $\cR(E)$. Let $F$,$G$ be subsets 
of  $\,\C$ such that $F$ is compatible with $S$ and $G$ is compatible with $T$. Then the following properties hold :
\begin{enumerate}
\item We have for $f \in \text{C}_b(F)$ and $g \in \text{C}_b(G)$ that $f(S) g(T) = g(T) f(S)$.
\item We have for $f \in \text{C}_b(F)$ and $g \in \text{C}(G)$ that $f(S) g(T) \subseteq g(T) f(S)$.
\item We have for $f \in \text{C}(F)$ and $g \in \text{C}_b(G)$ that $g(T) f(S) \subseteq f(S) g(T)$.
\end{enumerate}
\end{result}

\medskip

Combining result \ref{art9.res1} with the Fuglede-Putnam theorem (theorem \ref{art5.thm1}), we get immediately the following equivalent 
characterization of commuting normal operators.

\begin{proposition}  \label{art9.prop3}
Consider a Hilbert \cst-module $E$ over a \cst-algebra $A$ and let $S$,$T$ be normal elements in $\cR(E)$ and let $F$ be a subset of 
$\,\C$ which is compatible with $S$. Then $S$ commutes with $T$ $\Leftrightarrow$ We have for every $f \in \text{C}_0(F)$ that
$f(S) \, T  \subseteq T \, f(S)$.
\end{proposition}

\bigskip

In the case that one of the two (or both) are bounded, we get the usual commutation conditions.
The next result follows immediately from the previous one an the Fuglede-Putnam theorem.

\begin{result}
Consider a Hilbert \cst-module $E$ over a \cst-algebra $A$, let $S$ be a normal element in $\cL(E)$ and $T$ a normal element in $\cR(E)$. 
Then $S$ commutes with $T$ $\Leftrightarrow$
$S \, T \subseteq T \, S$.
\end{result}

\begin{corollary}
Consider a Hilbert \cst-module $E$ over a \cst-algebra $A$, let $S$,$T$ be normal elements in 
$\cL(E)$. Then $S$ commutes with $T$ $\Leftrightarrow$ $S \, T = T \, S$.
\end{corollary}

\bigskip

It is of course  not necessary to look at all the elements of $\text{C}_0(F)$. A certain well behaved subset will suffice.

\begin{result}
Consider a Hilbert \cst-module $E$ over a \cst-algebra $A$ and let $S$,$T$ be normal elements in $\cR(E)$ and let $F$ be a subset of  
$\,\C$ which is compatible with $S$.
Consider furthermore a bounded net $(e_i)_{i \in I}$ in $\text{K}(F)$ which converges strictly to $1$. 

Then $S$ and $T$ commute $\Leftrightarrow$ We have for every $i \in I$ that $(S \comp e_i(S)) \, T \subseteq T \, (S \comp e_i(S))$.
\end{result}
\begin{demo}
We already know that one implication is true. So we turn to the other one.

Therefore suppose that $(S \comp e_i(S)) \, T \subseteq T \, (S \comp e_i(S))$ for every $i \in I$.

Take $f \in \text{C}_0(\C)$. 

By the Fuglede-Putnam theorem, we get that $(S \comp e_i(S)) \, f(T) = f(T) \, (S \comp e_i(S))$ for every $i \in I$.

Choose $v \in D(S)$. Take $j \in I$.

We know  that $e_j(S) f(T) v \in D(S)$ and that $$S( e_j(S) f(T) v ) = 
(S \comp e_j(S)) \,  f(T) v = f(T) \, (S \comp e_j(S)) \, v = f(T)(e_j(S) S(v))$$

So we see that $(e_i(S) f(T) v)_{i \in I}$ converges to $f(T) v$ and that
$\bigl(S(e_i(S) f(T) v)\bigr)_{i \in I}$ converges to $f(T) S(v)$.

Therefore the closedness of $S$ implies that $f(T) v$ belongs to $D(S)$ and $S(f(T)v) = f(T)\,S(v)$.

\vspace{1mm}

So we have proven that $f(T) \, S \subseteq S \, f(T)$.
\end{demo}

\bigskip

Combining the first statement of result \ref{art9.res1} with proposition \ref{art3.prop2}, the next useful result becomes easy to prove.

\begin{corollary}
Consider a Hilbert \cst-module $E$ over a \cst-algebra $A$ and let $S$,$T$ be commuting  normal elements in $\cR(E)$. Let $F$,$G$ be subsets 
of  $\,\C$ such that $F$ is compatible with $S$ and $G$ is compatible with $T$. Then we have for every $f \in \text{C}(F)$ and 
$g \in \text{C}(G)$ that $f(S)$ and $g(T)$ commute.
\end{corollary}

\medskip

Because $z_T$ is a function of $T$ and equation 1.19 of \cite{Wor6}, we have of course the next corollary.

\begin{corollary}
Consider a Hilbert \cst-module $E$ over a \cst-algebra $A$ and let $S$,$T$ be normal elements in $\cR(E)$. Then $S$ commutes with $T$ 
$\Leftrightarrow$ $z_S \, z_T = z_T \, z_S$.
\end{corollary}

\medskip

The next result follows immediately from proposition \ref{art3.prop1}.

\begin{result}
Consider a Hilbert \cst-module $E$ over a \cst-algebra $A$. Let $B$ be a \cst-algebra $A$ and $\pi$ a non-degenerate $^*$-homomorphism 
from $B$ into $\cL(E)$. Suppose that $S$ and $T$ are two normal elements affiliated with $B$ such that $S$ and $T$ commute.
Then $\pi(S)$ and $\pi(T)$ commute.
\end{result}

\bigskip

As could be expected, the functional calculi of two commuting normal operators give a functional calculus on $\,\C \times \C$. This is 
essentially the content of proposition \ref{art9.prop2}.

\medskip

First, we state some unicity result.

\begin{result}
Consider a Hilbert \cst-module $E$ over a \cst-algebra $A$ and let $S$,$T$ be commuting normal elements in $\cR(E)$ and let $F$,$G$ be an
almost closed subsets of  $\,\C$.
Define $p$ as the projection of $F \times G$ on $F$ and $q$ as the projection of $F \times G$ on $G$. Consider non-degenerate 
$^*$-homomorphism $\pi$,$\th$ from $\text{C}_0(F \times G)$ into $\cL(E)$. 
If $\pi(p) = \th(p)$ and $\pi(q) = \th(q)$, then $\pi = \th$.
\end{result}
\begin{demo}
Choose $f \in \text{C}_0(F)$ and $g \in \text{C}_0(G)$. 

By lemma \ref{art3.lem1}, there exists $h \in \text{C}_0(\C)$ such that $f \subseteq h$. 
Using proposition \ref{prel3.prop1}, we see that $\pi(f \ot 1) = \pi(h \!\circ\! p) = \pi(h(p)) = h(\pi(p))$ and similarly,  
$\th(f \ot 1) = h(\th(p))$. Hence, $\pi(f \ot 1) = \th(f \ot 1)$.

One proves in the same way that $\pi(1 \ot g) = \th(1 \ot g)$. Consequently,
$\pi(f \ot g) = \th(f \ot g)$.
\end{demo}

\medskip

\begin{proposition} \label{art9.prop2}
Consider a Hilbert \cst-module $E$ over a \cst-algebra $A$ and let $S$,$T$ be commuting normal elements in $\cR(E)$ and let 
$F$,$G$ be subsets of  $\,\C$ such that $F$ is compatible with $S$ and $G$ is compatible with $G$.  Define $p$ as the projection 
of $F \times G$ on $F$ and $q$ as the projection of $F \times G$ on $G$.

Then there exists a unique non-degenerate $^*$-homomorphism $\pi$ from $\text{C}_0(F \times G)$ into $\cL(E)$ such that $\pi(p) = S$ 
and $\pi(q) = T$. We call $\pi$ the functional calculus of $S$,$T$ on $F \times G$. Furthermore,
\begin{itemize}
\item We have for every $f \in \text{C}(F)$ that $\pi(f \ot 1) = f(S)$
\item We have for every $g \in \text{C}(G)$ that $\pi(1 \ot g) = g(T)$
\end{itemize}
\end{proposition}
\begin{demo}
Call $\th$ the functional calculus of $S$ on $F$ and $\rho$ the functional calculus of
$T$ on $G$. Because $S$ and $T$ commute, we have by result \ref{art9.res1} that $\th$ and $\rho$ have commuting ranges. Hence, the universal 
property of the maximal tensor product implies the existence of a unique $^*$-homomorphism $\pi$ from $\text{C}_0(F \times G)$ such that
$\pi(f \ot g) = \th(f) \rho(g) = \rho(g) \th(f)$ for $f \in \text{C}_0(F)$ and 
$g \in \text{C}_0(G)$.

\medskip

Take $h \in \text{C}(F)$

Choose an approximate unit $(e_i)_{i \in I}$ for $\text{C}_0(F)$ in $\text{K}(F)$ and an approximate unit $(u_j)_{j \in J}$ for $\text{C}_0(G)$ in 
$\text{K}(G)$.
Then $(e_i \ot u_j)_{(i,j) \in I \times J}$ is an approximate unit for $\text{C}_0(F \times G)$ in $\text{K}(F \times G)$.

Hence, $(\pi(e_i \ot u_i))_{(i,j) \in I \times J}$ is a bounded net which converges strongly to $1$. 

\medskip

Fix $i \in I$, $j \in J$.

Because $e_i \ot u_j$ belongs to $\text{K}(F \ot G)$, lemma \ref{art1.lem1} implies that $\pi(e_i \ot u_j) \, E
\subseteq D(\pi(h \ot 1))$ and 
$$\pi(e_i \ot u_j) \, \pi(h \ot 1)  \subseteq \pi(h \ot 1) \, \pi(e_i \ot u_j) \hspace{1cm} (1) $$

We have also that $\pi(e_i \ot u_j) = e_i(S) \, u_j(T)$. This implies that $\pi(e_i \ot u_j) \, E \subseteq D(h(S))$. Because $S$ and $T$ commute, 
we have moreover that
$$\pi(e_i \ot u_j) \, h(S) \subseteq h(S)  \, \pi(e_i \ot u_j) \hspace{1cm} (2)$$

We see also that
\begin{eqnarray*}
& & \pi(h \ot 1) \, \pi(e_i \ot u_j) =  \pi((h \ot 1)(e_i \ot u_j)) = \pi((h e_i) \ot u_j) \\
& & \spat = (h e_i)(S) \, u_j(T) = h(S)\, e_i(S)\, u_j(T) = h(S) \, \pi(e_i \ot u_j) 
\hspace{1cm} (3)
\end{eqnarray*}

Combining these three results and using the closedness of both $\pi(h \ot 1)$ and $h(S)$, we get easily that $\pi(h \ot 1) = h(S)$.

\medskip

So we get in particular that $\pi(p) = \pi(\io_F \ot 1) = \io_F(S) = S$.

\medskip

The results concerning elements in $\text{C}(G)$ are proven in the same way.
\end{demo}

\medskip

\begin{notation}
Consider a Hilbert \cst-module $E$ over a \cst-algebra $A$ and let $S$,$T$ be commuting normal elements in $\cR(E)$ and let $F$,$G$ be subsets of  
$\,\C$ such that $F$ is compatible with $S$ and $G$ is compatible with $T$. Call $\pi$ the funtional calculus of $S$,$T$ on $F \times G$. For every 
$h \in \text{C}(F \times G)$, we define $h(S,T) = \pi(h)$, so $h(S,T)$ is a normal element in $\cR(E)$.
\end{notation}

\medskip

\begin{result}
Consider a Hilbert \cst-module $E$ over a \cst-algebra $A$ and let $S$,$T$ be commuting normal elements in $\cR(E)$ and let $F$,$G$ be subsets of  
$\,\C$ such that $F$ is compatible with $S$ and $G$ is compatible with $T$. Then we have the following properties :
\begin{itemize}
\item We have for every $f \in \text{C}(F)$ that $(f \ot 1)(S,T) = f(S)$.
\item We have for every $g \in \text{C}(G)$ that $(1 \ot g)(S,T) = g(T)$.
\end{itemize}
\end{result}

\medskip

\begin{corollary}
Consider a Hilbert \cst-module $E$ over a \cst-algebra $A$ and let $S$,$T$ be commuting normal elements in $\cR(E)$ and let $F$,$G$ be subsets of  
$\,\C$ such that $F$ is compatible with $S$ and $G$ is compatible with $T$. Then we have for every $f \in \text{C}_b(F)$ and
$g \in \text{C}_b(G)$ that $(f \ot g)(S,T) = f(S) g(T) = g(T) f(S)$.
\end{corollary}

\bigskip

Of course, restricting functions does not change the value.

\begin{result}
Consider a Hilbert \cst-module $E$ over a \cst-algebra $A$ and let $S$,$T$ be commuting normal elements in $\cR(E)$ and let $F$,$G$,$H$,$K$ be 
subsets of  $\,\C$ such that $F$,$H$ are compatible with $S$,  $G$,$K$ are compatible with $T$ such that $F \subseteq H$ and
$G \subseteq K$. Then we have for $h \in \text{C}(H \times K)$ that
$h(S,T) = h_{F \times G}(S,T)$.
\end{result}

The proof of this result is the same as the proof of lemma \ref{art3.lem2}.

\bigskip

As usual, the next result turns out to be very useful. It is an immediate consequence of proposition \ref{art3.prop1} and result \ref{art3.res4}.

\begin{result} \label{art9.res3}
Consider a Hilbert \cst-module $E$ over a \cst-algebra $A$ and let $S$,$T$ be commuting normal elements in $\cR(E)$ and let $F$,$G$ be subsets of  
$\,\C$ such that $F$ is compatible with $S$ and $G$ is compatible with $T$. Consider $k \in \text{C}(F \times G)$ and
consider an almost closed subset $H$ of  $\,\C$ such that $k(F \ot G) \subseteq H$. Then we have the following properties :
\begin{itemize}
\item The set $H$ is compatible with $k(S,T)$.
\item We have for every $h \in \text{C}(H)$ that $h(k(S,T)) = (h \!\circ\! k)(S,T)$.
\end{itemize}
\end{result}

\bigskip

This functional calculus for the pair $S$,$T$ allows us to define the product and the sum of two commuting normal operators. Call $p$ the projection 
of $\,\C \times \C$ on the first factor and $q$ the projection of $\,\C \times \C$ on the second factor. Then $p(S,T) = S$ and 
$q(S,T) = T$. Hence, looking at $(p+q)(S,T)$ and $(p q)(S,T)$ and refering to results
\ref{art1.res4} and \ref{art1.res2}, the following definition is justified.

\begin{definition}
Consider a Hilbert \cst-module $E$ over a \cst-algebra $A$ and let $S$,$T$ be commuting normal elements in $\cR(E)$.
\begin{trivlist}
\item[\ \, $\bullet$] The element $S \, T$ is closable and we define $S \comp T = \overline{S \, T}$. Then $S \comp T$ is a normal element in $\cR(E)$.
\item[\ \, $\bullet$] The element $S + T$ is closable and we define $S \dplus T = \overline{S + T}$. Then $S \dplus T$ is a normal element in $\cR(F)$.
\end{trivlist}
\end{definition}

\medskip

Refering to results \ref{art1.res4} and \ref{art1.res2} once more, we get immediately the following result.

\begin{corollary} \label{art9.cor1}
Consider a Hilbert \cst-module $E$ over a \cst-algebra $A$ and let $S$,$T$ be commuting normal elements in $\cR(E)$ and let $F$,$G$ be subsets of  
$\,\C$ such that $F$ is compatible with $S$ and $G$ is compatible with $T$. Consider $f \in \text{C}(F)$ and
$g \in \text{C}(G)$. Then
\begin{itemize}
\item $f(S) \comp g(T) = g(T) \comp f(S) = (f \ot g)(S,T)$
\item $f(S) \dplus g(T) = g(T) \dplus f(S) = (f \ot 1 + 1 \ot g)(S,T)$.
\end{itemize}
\end{corollary}

So we have in particular  that $S \dplus T = T \dplus S$ and $S \comp T = T \comp S$.

\bigskip

We will not pursue this matter any further except for the next two results which will be useful in the next section. 

\begin{result}
Consider a Hilbert \cst-module $E$ over a \cst-algebra $A$ and let $S$,$T$ be commuting normal elements in $\cR(E)$.
\begin{itemize}
\item If $S$ and $T$ are invertible, then $S \comp T$ is invertible.
\item If $S$ and $T$ are selfadjoint, then $S \comp T$ is selfadjoint.
\item If $S$ and $T$ are positive, then $S \comp T$ is positive.
\item If $S$ and $T$ are strictly positive, then $S \comp T$ is strictly positive.
\end{itemize}
\end{result}

\medskip

\begin{proposition} \label{art9.prop4}
Consider a Hilbert \cst-module $E$ over a \cst-algebra $A$ and let $S$,$T$ be commuting normal elements in $\cR(E)$. 
\begin{enumerate}
\item We have that $(S \comp T)^n = S^n \comp T^n$ for $n \in \N$.
\item If $S$ and $T$ are invertible, then $(S \comp T)^n = S^n \comp T^n$ for $n \in \Z$.
\item If $S$ and $T$ are positive, then $(S \comp T)^r = S^r \comp T^r$ for $r \in \R^+$.
\item If $S$ and $T$ are strictly positive, then $(S \comp T)^z = S^z \comp T^z$ for $z \in \C$.
\end{enumerate}
\end{proposition}

\medskip

\begin{remark} \rm 
As an example, we will prove the  first statement of the first result and the second statement of the second one. The proofs of the other 
ones are completely similar.

So take commuting invertible normal elements $S$,$T$ in $\cR(E)$.
Then $\C_0$ is compatible with $S$ and $T$. By corollary \ref{art9.cor1}, we get that
$$S \comp T = \io_{\Cvar_0}(S) \comp \io_{\Cvar_0}(T) = (\io_{\Cvar_0} \ot \io_{\Cvar_0})(S,T)$$
It is clear that $\io_{\Cvar_0} \ot \io_{\Cvar_0}$ is invertible in $\text{C}(\C_0 \times \C_0)$. This implies that $S \comp T$ is invertible.

\vspace{1mm}

Choose $n \in \Z$ and define the function $h$ from $\C_0$ into $\C_0$ such that
$h(c) = c^n$ for $c \in \C_0$. Then result \ref{art9.res3} implies that
\begin{eqnarray*}
(S \comp T)^n & = & h(S \comp T) = h\bigl((\io_{\Cvar_0} \ot \io_{\Cvar_0})(S,T)\bigr) 
= \bigl(h \!\circ\! (\io_{\Cvar_0} \ot \io_{\Cvar_0})\bigr)(S,T) \\
& = & (h \ot h)(S,T) = h(S) \comp h(T) = S^n \comp T^n
\end{eqnarray*}
\end{remark}

\bigskip\medskip

We end this section with the familiar result that the commuting relation involving strictly positive elements can be stated in terms of the 
unitary groups they generate.

\medskip

In \cite{Kus}, we prove the following Hilbert \cst-module version of a well known Hilbert space result (see e.g. \cite{Stra}).

\medskip

For any $z \in \C$, we use the notation $S(z) = \{ \, y \in \C \mid \text{Im}\,y \in [0,\text{Im}\,z] \, \}$ .

\begin{proposition} \label{art9.prop1}
Consider a Hilbert \cst-module $E$ over a \cst-algebra $A$ and let $S$ be a strictly positive element in $\cR(E)$. Let $z$ be a complex 
number and $v \in E$. 

Then we have that $v$ belongs to $D(S^{i z})$ $\Leftrightarrow$ There exists a function $f$ from $S(z)$ into $E$ such that $f$ is continuous 
on $S(z)$, $f$ is analytic on $S(z)^0$ and $f(t) = S^{i t} v$ for $t \in \R$.

If such a function $f$ exists, then $S^{iz} v = f(z)$.
\end{proposition}

\medskip

This allows us to prove the following commutation criteria.

\begin{result}
Consider a Hilbert \cst-module $E$ over a \cst-algebra $A$ and let $S$ be a strictly positive 
element in $\cR(E)$. Let $T$ a normal element in $\cR(E)$. 

Then $S$ commutes with $T$ $\Leftrightarrow$ We have for every $t \in \R$ and that $S^{it} \, T \subseteq T \, S^{it}$.
\end{result}
\begin{demo}
Suppose that $S^{it} \, T \subseteq T \, S^{it}$ for every $t \in \R$.

Choose $g \in \text{C}_0(G)$. Then Fuglede-Putnam implies that $S^{i t} \, g(T) = g(T) \, S^{it}$. The previous proposition implies 
(with a small effort) that $g(T) S \subseteq S g(T)$.
Hence, $S$ and $T$ commute by proposition \ref{art9.prop3} .
\end{demo}

\medskip

Applying the previous result, we get the following one.

\begin{corollary}
Consider a Hilbert \cst-module $E$ over a \cst-algebra $A$ and let $S$,$T$ be strictly positive 
element in $\cR(E)$. Then $S$ and $T$ commute $\Leftrightarrow$ We have for every $s,t \in \R$ that $S^{is} \, T^{it} = T^{it} \, S^{is}$.
\end{corollary}

\bigskip

\section{Tensor products of Hilbert \cst-modules}

In this section, we will look at some basic properties of regular operators  on the tensor products of Hilbert \cst-modules. 
We will always work with the minimal tensor product between \cst-algebras.

\medskip

In chapter 4 of \cite{Lan}, the tensor product of two Hilbert \cst-modules was defined in the following way.

\begin{definition}
Consider a Hilbert \cst-module $E$ over a \cst-algebra $A$ and a Hilbert \cst-module $F$ over a \cst-algebra $B$. Then $E \ot F$ 
is defined to be a Hilbert \cst-module over $A \ot B$ such that
\begin{itemize}
\item The set $E \od F$ is a dense subspace of $E \ot F$.
\item We have for every $v_1,v_2 \in E$ and $w_1,w_2 \in F$ that
$\lan v_1 \ot w_1 , v_2 \ot w_2 \ran = \lan v_1,v_2 \ran \ot \lan w_1 , w_2 \ran$.
\end{itemize}
\end{definition}

\medskip

In the same chaptet of \cite{Lan}, the tensor product of two adjointable operators was defined in the following way.

\begin{proposition}
Consider Hilbert \cst-modules $E$,$F$ over a \cst-algebra $A$ and Hilbert \cst-modules $G$,$H$ over a \cst-algebra $B$. Let $S \in \cL(E,F)$ 
and $T \in \cL(G,H)$. Then there exists a unique $R \in \cL(E \ot F, G \ot H)$ such that $R(v \ot w) = S(v) \ot T(w)$ for every $v \in E$ and 
$w \in G$. We have moreover that $\|R\| = \|S\| \, \|T\|$.
\end{proposition}

\medskip

\begin{lemma}
Consider Hilbert \cst-modules $E$,$F$ over a \cst-algebra $A$ and Hilbert \cst-modules $G$,$H$ over a \cst-algebra $B$. By the universality 
property of the algebraic tensor product, there exists a unique linear map $\pi$ from $\cL(E,F) \od \cL(G,H)$ into $\cL(E \ot G, F \ot H)$ 
such that $\pi(S \ot T) (v \ot w) = S(v) \ot T(w)$ for $S \in \cL(E,F)$, $T \in \cL(G,H)$, $v \in E$, $w \in G$. Then $\pi$ is injective.
\end{lemma}
\begin{demo}
Take $x \in \cL(E,F) \od \cL(G,H)$ such that $\pi(x) = 0$. Then there exists $S_1,\ldots\!,S_n \in \cL(E,F)$ and $T_1,\ldots\!,T_n \in \cL(G,H)$ 
such that $x = \sum_{i=1}^n S_i \ot T_i$
and such that $S_1,\ldots\!,S_n$ are linearly independent.

Choose $v \in G$, $w \in H$. Take $\om \in A^*$

Fix $p \in E$, $q \in F$. We have that
$$\sum_{i=1}^n \lan S_i(p), q \ran \ot \lan T_i(v) , w \ran
= \lan \pi(x) (p \ot v) , q \ot w \ran = 0 \ . $$

Applying $\om \ot \io$ to this equation gives us that
$$0 = \sum_{i=1}^n \lan S_i(p) , q \ran \,\, \om(\lan T_i(v) , w \ran)
= \lan [\, \sum_{i=1}^n \om(\lan T_i(v) , w \ran) \, S_i ] (p) , q \ran $$
Hence, we see that $\sum_{i=1}^n \om(\lan T_i(v) , w \ran) \, S_i = 0$. By the linear independence of $S_1,\ldots\!,S_n$, we get for every 
$i \in \{1,\ldots\!,n\}$ that
$\om(\lan T_i(v) , w \ran) = 0$.

\medskip

Consequently, $T_1 = \ldots = T_n = 0$. So $x = 0$.
\end{demo}

\medskip

From now on, we will forget about the mapping $\pi$ and consider $\cL(E,F) \od \cL(G,H)$ as a subspace of $\cL(E \ot G , F \ot H)$. So we 
get the following result.

\begin{result}
Consider Hilbert \cst-modules $E$,$F$ over a \cst-algebra $A$ and  Hilbert \cst-modules $G$,$H$ over a \cst-algebra $B$. Let $S \in \cL(E,F)$ 
and $T \in \cL(G,H)$.
Then $S \ot T$ is the unique element in $\cL(E \ot G, F \ot H)$ such that
$(S \ot T)(v \ot w) = S(v) \ot S(w)$ for $v \in E$ and $w \in G$. We have moreover that
$\|S \ot T\| = \|S\| \, \|T\|$.
\end{result}

\medskip

As in the case of Hilbert spaces, the minimal tensor product norm of an element in $\cL(E) \od \cL(F)$ is equal to the norm of it as an 
operator on $E \ot F$.

\begin{proposition}
Consider a Hilbert \cst-module $E$ over a \cst-algebra $A$ and  a Hilbert \cst-module $F$ over a \cst-algebra $B$. Then we have for every 
$x \in \cL(E) \od \cL(F)$ that
$\|x\|_{\min} = \|x\|$ (where the last norm is the norm of $x$ as an element in $\cL(E \ot F)$).
\end{proposition}
\begin{demo}
If we restrict the norm on $\cL(E \ot F)$ to $\cL(E) \od \cL(F)$, we get a \cst-norm on $\cL(E) \od \cL(F)$. This implies immediately that 
$\|x\|_{\min} \leq \|x\|$. We will turn to the other inequality.

\medskip

Take $\vep > 0$

Use the notations of section \ref{art8}. Define the set $K = \{ \, \rho \ot \th \mid \rho \in A^*_+ , \th \in B^*_+ \, \}$. Because we 
work with the minimal tensor product on $A \ot B$, the set $K$ is separating for $A \ot B$.  

Define the mapping $\eta$ from $\cL(E \ot F)$ into $\prod_{\om \in K} \cB((E \ot F)_\om)$ such that $(\eta(y))_\om = y_\om$ for every 
$y \in \cL(E \ot F)$ and $\om \in K$.
Then lemma \ref{art8.lem1} implies that $\eta$ is an injective $^*$-homomorphism and consequently isometric. 
Hence, $$\|x\| =  \|\eta(x)\| = \sup \, \{\, \|x_\om\| \mid \om \in K \, \}$$

This implies the existence of $\rho \in A^*_+$ and $\th \in B^*_+$ such that
$\|x_{\om \ot \th}\| \geq \|x\| - \vep$.

Define the mappings $\eta_\rho : \cL(E) \rightarrow \cB(E_\rho) : a \mapsto a_\rho$ and
$\eta_\th : \cL(F) \rightarrow \cB(F_\th) : b \mapsto b_\th$. Then $\eta_\rho$ and $\eta_\th$ are $^*$-homomorphisms. Hence we get the 
inequality  $\|x\|_{\min} \geq \|(\eta_\rho \ot \eta_\th)(x)\|$.

It is easy to check that $\lan \overline{v_1} \ot \overline{w_1} , \overline{v_2} \ot \overline{w_2} \ran = \lan \, \overline{v_1 \ot w_1} , 
\overline{v_2 \ot w_2} \, \ran$ for
$v_1,v_2 \in E$ and $w_1,w_2 \in F$. This implies the existence of a unitary operator $U$ from 
$E_\rho \ot F_\th$ to $(E \ot F)_{\rho \ot \th}$ such that $U(\overline{v} \ot \overline{w}) 
= \overline{v \ot w}$ for $v \in E$ and $w \in F$.

Now it is straightforward to check that $(\eta_\rho \ot \eta_\th)(x) = U^* \, x_{\om \ot \th} \, U$ which implies that 
$$\|x\|_{\min} \geq \|(\eta_\rho \ot \eta_\th)(x)\| = \|x_{\om \ot \th}\| \geq \|x\| - \vep $$

\medskip

From this all, we can conclude that  $\|x\|_{\min} \geq \|x\|$. Hence, 
$\|x\|_{\min} = \|x\|$.
\end{demo}

\medskip

So from now on, we can consider $\cL(E) \ot \cL(F)$ as the closure of $\cL(E) \od \cL(F)$ in
$\cL(E \ot F)$.

\bigskip\medskip

It is even possible to define the tensor product of two regular operators (see theorem 6.1
of \cite{Wor7} and chapter 10 of \cite{Lan}).

\begin{definition}  \label{art10.def1}
Consider Hilbert \cst-modules $E$,$F$ over a \cst-algebra $A$ and  Hilbert \cst-modules $G$,$H$ over a \cst-algebra $B$. Consider moreover 
elements $S \in \cR(E,F)$ and 
$T \in \cR(G,H)$. Then $S \od T$ is closable and we define $S \ot T$ to be the closure of $S \od T$.
\end{definition}

It is not difficult to check that $C \od D$ is a core for $S \ot T$ if $C$ is a core for
$S$ and $D$ is a core for $T$.

\begin{proposition}
Consider Hilbert \cst-modules $E$,$F$ over a \cst-algebra $A$ and  Hilbert \cst-modules $G$,$H$ over a \cst-algebra $B$. Consider moreover 
elements $S \in \cR(E,F)$ and 
$T \in \cR(G,H)$. Then $S \ot T$ belongs to $\cR(E \ot G, F \ot H)$  and $(S \ot T)^* = S^* \ot T^*$.
\end{proposition}

\medskip

The tensor product of normal (resp. selfadjoint, positive) elements will again be normal
(resp. selfadjoint, positive):

\begin{result}
Consider a Hilbert \cst-module $E$ over a \cst-algebra $A$ and  a Hilbert \cst-module $F$ over a \cst-algebra $B$. Consider moreover elements 
$S \in \cR(E)$ and 
$T \in \cR(F)$. Then we have the following properties.
\begin{itemize}
\item If $S$ and $T$ are normal, then $S \ot T$ is normal.
\item If $S$ and $T$ are selfadjoint, then $S \ot T$ is selfadjoint.
\item If $S$ and $T$ are positive, then $S \ot T$ is positive.
\item If $S$ and $T$ are strictly positive, then $S \ot T$ is strictly positive.
\end{itemize}
\end{result}
\begin{demo}
\begin{itemize}
\item The case for selfadjoint elements follows easily by the previous proposition.
\item Suppose that $S$ and $T$ are normal. By polarisation, we have that
  $\lan S(v) , S(w) \ran = \lan S^*(v) , S^*(w) \ran$ for every $v,w \in D(S)$.
A similar remark applies to $T$. Using these equalities, it is straightforward to check
that $D(S \od T) = D(S^* \od T^*)$ and that $$\lan (S \od T)(u) , (S \od T)(u) \ran
= \lan (S^* \od T^*)(u) ,  (S^* \od T^*)(u) \ran$$ for every $u \in D(S \od T)$.

This implies easily that $D(S \ot T) = D(S^* \ot T^*)$  and that $$\lan (S \ot T)(u) ,  (S \ot
T)(u) \ran = \lan (S^* \ot T^*)(u) , (S^* \ot T^*)(u) \ran$$ for every $u \in D(S \ot T)$.
So $S \ot T$ is normal by definition.

\item Suppose that $S$ and $T$ are positive. We know already that $S \ot T$ is selfadjoint.

Choose $v_1, \ldots\! , v_n \in D(S)$ and $w_1, \ldots\!, w_n \in D(T)$. Define $M \in
M_n(A)$ and $N \in M_n(B)$ such that $M_{ij} = \lan S(v_i) , v_j \ran$ and $N_{ij} =
\lan T(w_i) , w_j \ran$ for every $i,j \in \{1,\ldots\!,n\}$.

Using lemma 3.2 of \cite{Tak} and  the positivity of $S$, it follows readily that the
matrix $M$ is positive. Similarly we find that $N$ is positive. Using lemma 4.3 of
\cite{Lan}, we find that $\sum_{i,j=1}^n M_{ij} \ot N_{ij}$ is positive. But we have also
that
$$\sum_{i,j=1}^n M_{ij} \ot N_{ij} = \lan  (S \ot T) (\sum_{i=1}^n v_i \ot w_i) ,
(\sum_{j=1}^n v_j \ot w_j) \ran $$ which implies that the right hand side is positive.

So we have proven that $\lan (S \ot T)(u) , u \ran$ is positive for every $u \in D(S) \od
D(T)$. From this we get that $\lan (S \ot T)(u) , u \ran$ is positive for every $u \in D(S
\ot T)$. So $S \ot T$ is positive.
\item The result concerning strict positivity follows now easily.
\end{itemize}
\end{demo}

\medskip

\begin{result}
Consider Hilbert \cst-modules $E$,$F$ over a \cst-algebra $A$ and  Hilbert \cst-modules $G$,$H$ over a \cst-algebra $B$. 
Consider moreover elements $S \in \cR(E,F)$ and 
$T \in \cR(G,H)$.
\begin{itemize}
\item If $S$ and $T$ are invertible, then $S \ot T$ is invertible and $(S \ot T)^{-1} = S^{-1} \ot T^{-1}$.
\item If $S$ and $T$ are adjointable invertible, then $S \ot T$ is adjointable invertible.
\end{itemize}
\end{result}
\begin{demo}
Suppose that $S$ and $T$ are invertible. Using proposition \ref{art2.prop1}, it is easy to see that $S \ot T$ and $S^* \ot T^*$ 
are invertible.

We have moreover that $S^{-1} \od T^{-1} = (S \od T)^{-1} \subseteq (S \ot T)^{-1}$ which implies that $S^{-1} \ot T^{-1} 
\subseteq (S \ot T)^{-1}$.

In the same way, it follows that $(S^*)^{-1} \ot (T^*)^{-1} \subseteq (S^* \ot T^*)^{-1}$, which implies that $(S^{-1} \ot T^{-1})^* 
\subseteq \bigl((S \ot T)^{-1}\bigr)^*$.
Hence, $(S \ot T)^{-1} \subseteq S^{-1} \ot T^{-1}$. Consequently, $S^{-1} \ot T^{-1} = (S \ot T)^{-1}$.

\medskip

The second statement follows immediately from the first.
\end{demo}

\medskip

\begin{result}
Consider Hilbert \cst-modules $E$,$F$ over a \cst-algebra $A$ and  Hilbert \cst-modules $G$,$H$ over a \cst-algebra $B$. Consider 
moreover elements $S \in \cR(E,F)$ and  $T \in \cR(G,H)$. Then we have the equalities
$$(S \ot T)^* (S \ot T) = S^* S \ot T^* T \hspace{3cm} (S \ot T)(S \ot T)^* = S S^* \ot T^* T$$
\end{result}
\begin{demo}
We have that
$$S^* S \od T^* T \subseteq (S^* \od T^*)(S \od T) \subseteq (S^* \ot T^*)(S \ot T) = (S \ot T)^* (S \ot T) \ .$$
This implies that $S^* S \ot T^* T \subseteq (S \ot T)^* (S \ot T)$. Because these two
elements are selfadjoint, this inclusion is an equality.
\end{demo}

\medskip

Now we have our obligatory non-degenerate $^*$-homomorphism result.

\begin{proposition}
Consider a Hilbert \cst-module $E$ over a \cst-algebra $A$ and a Hilbert \cst-module $F$ over a \cst-algebra $B$. Let $C$ and $D$ be 
\cst-algebras, $\pi$ a non-degenerate $^*$-homomorphism from $C$ into $\cL(E)$ and $\th$ a non-degenerate $^*$-homomorphism from $D$ 
into $\cL(F)$. Then we have for every $S \in \cR(C)$ and $T \in \cR(D)$ that
$(\pi \ot \th)(S \ot T)$.
\end{proposition}
\begin{demo}
Because $D(S) \od D(T)$ is a core for $S \ot T$ and $E \od F$ is dense in $E \ot F$, the remark after theorem \ref{prel1.thm1} implies that 
$(\pi \ot \th)(D(S) \od D(T))(E \od F)$ is a core for $(\pi \ot \th)(S \ot T)$. Hence, $\pi(D(S)) E \od \th(D(T)) F$ is a core for 
$(\pi \ot \th)(S \ot T)$.

Because $\pi(D(S))E$ is a core for $\pi(S)$ and $\th(D(T))F$ is a core for $\th(T)$, the remark after definition \ref{art10.def1} gives us 
that $\pi(D(S)) E \od \th(D(T)) F$ is also a core for $\pi(S) \ot \th(T)$. 

It is easy to see that $$(\pi \ot \th)(S \ot T) \, (\pi(c) v \ot \th(d) w)
= (\pi(S) \ot \th(T)) \, (\pi(c) v \ot \th(d) w)$$
for $c \in D(S)$, $d \in D(T)$, $v \in E$ and $w \in F$. 
Consequently, $(\pi \ot \th)(S \ot T) =  \pi(S) \ot \th(T)$.
\end{demo}

\bigskip\medskip

\begin{lemma}
Consider a Hilbert \cst-module $E$ over a \cst-algebra $A$ and a Hilbert \cst-module $F$ over a \cst-algebra $B$. Define the injective 
non-degenerate $^*$-homomorphism $\pi$ from $\cK(E)$ into $\cL(E \ot F)$ such that $\pi(x) = x \ot 1$ for every $x \in \cK(E)$. 
Then we have for every $S \in \cR(\cK(E))$ that $\pi(S) = \tilde{S} \ot 1$.
\end{lemma}
\begin{demo}
We know that $D(S) E$ is a core for $\tilde{S}$ (result \ref{art7.res1}). This implies that $D(S) E \od F$ is a core for $\tilde{S} \ot 1$. 
On the other hand, we know also that $\pi(D(S))(E \od F)$ is a core for $\pi(S)$. This implies that $D(S) E \od F$ is also a core for $\pi(D(S))$.

We have moreover for every $x \in D(S)$ and $v \in E$, $w \in F$ that 
$$\pi(S)(x v \ot w) = \pi(S)(\pi(x)(v \ot w)) = \pi(S(x)) (v \ot w)
= S(x) v \ot w = \tilde{S}(x v) \ot w = (\tilde{S} \ot 1)( x v \ot w) $$

which implies that $\pi(S) = \tilde{S} \ot 1$.
\end{demo}

\medskip 

By corollary \ref{art7.cor1} and proposition \ref{art3.prop1} , this implies immediately the following result.

\begin{corollary} \label{art10.cor1}
Consider a Hilbert \cst-module $E$ over a \cst-algebra $A$ and a Hilbert \cst-module $F$ over a \cst-algebra $B$. Let $S$ be a normal element in 
$\cR(E)$ and $G$ a subset of $\,\C$ which is compatible with $S$. Then we have the following properties :
\begin{itemize}
\item The set $G$ is compatible with $S \ot 1$.
\item We have for every $f \in \text{C}(G)$ that $f(S \ot 1) = f(S) \ot 1$.
\end{itemize}
\end{corollary}

\medskip

We have of course similar results for elements of the form $1 \ot T$. 

\medskip

\begin{result} \label{art10.res1}
Consider a Hilbert \cst-module $E$ over a \cst-algebra $A$ and a Hilbert \cst-module $F$ over a \cst-algebra $B$. Let $S$ be a normal element in 
$\cR(E)$ and $T$ a normal element in $\cR(F)$. Then $S \ot 1$ and $1 \ot T$ commute and $(S \ot 1) \comp (1 \ot T) = S \ot T$.
\end{result}
\begin{demo}
The fact that $S \ot 1$ and $1 \ot T$ commute follows immediately from the previous result and the remark after it.
We have moreover that
$$S \od T  \subseteq (S \od 1)(1 \od T) \subseteq (S \ot 1)(1 \ot T) = (S \ot 1) \comp (1 \ot T)$$
which implies that $S \ot T \subseteq (S \ot 1) \comp (1 \ot T)$.
Because both are normal, this inclusion is an equality.
\end{demo}

\medskip

So we can also look at functional calculi for the pair $S$,$T$.

\medskip

Looking at the previous result,  corollary \ref{art10.cor1} and corollary \ref{art9.cor1}, we see that the following holds.

\begin{result}
Consider a Hilbert \cst-module $E$ over a \cst-algebra $A$ and a Hilbert \cst-module over a \cst-algebra $B$. Let $S$ be a normal element in 
$\cR(E)$ and $T$ a normal element in $\cR(F)$.
Consider subsets $F$,$G$ of $\,\C$ such that $F$ is compatible with $S$ and $G$ is compatible with $T$. Then we have for every 
$f \in \text{C}(F)$ and $g \in \text{C}(G)$ that $(f \ot g)(S \ot 1,1 \ot T) = f(S) \ot g(T)$.
\end{result}

\medskip

The next result is an easy consequence of result \ref{art10.res1}, corollary \ref{art10.cor1} and proposition \ref{art9.prop4}.

\begin{proposition}
Consider a Hilbert \cst-module $E$ over a \cst-algebra $A$ and a Hilbert \cst-module over a \cst-algebra $B$. Let $S$ be a normal element in 
$\cR(E)$ and $T$ a normal element in $\cR(F)$.
\begin{enumerate}
\item We have that $(S  \ot T)^n = S^n  \ot T^n$ for $n \in \N \cup \{0\}$.
\item If $S$ and $T$ are invertible, then $(S  \ot T)^n = S^n  \ot T^n$ for $n \in \Z$.
\item If $S$ and $T$ are positive, then $(S  \ot T)^r = S^r  \ot T^r$ for $r \in \R^+$.
\item If $S$ and $T$ are strictly positive, then $(S  \ot T)^z = S^z  \ot T^z$ for $z \in \C$.
\end{enumerate}
\end{proposition}

\medskip

\begin{corollary}
Consider Hilbert \cst-modules $E$,$F$ over a \cst-algebra $A$ and  Hilbert \cst-modules $G$,$H$ over a \cst-algebra $B$. Consider moreover 
elements $S \in \cR(E,F)$ and  $T \in \cR(G,H)$. Then we have that $|S \ot T| = |S| \ot |T|$.
\end{corollary}

\bigskip\medskip

We end this section with a problem connected with remark \ref{art4.rem1}. It is not clear to me under which conditions it is possible to cut 
out bigger pieces out of the spectrum and retain a functional calculus on this set. 

If you have a continuous function on a locally compact space, you can easily cut out closed sets which do not meet the image of this function.

\medskip

I wonder in particular if the following is true.

\begin{question}
Consider a Hilbert \cst-module $E$ over a \cst-algebra $A$ and let $T$ be a normal element in $\cR(E)$. Consider moreover a closed subset $K$ of 
$\si(T)$  such that the closure of 
$\,T \ot 1 - 1 \ot \io_K$ is invertible in $\cR(E \ot \text{C}_0(K))$.

Does there exists in this case a  non-degenerate $^*$-homomorphism $\pi$ from $\text{C}_0(\si(T) \setminus K)$ into $\cL(E)$ such that 
$\pi(\io_{\si(T) \setminus K} ) = T$.
\end{question}

It is not hard to check that the answer is affirmative in the commutative case and the case where $K$ is finite (by proposition \ref{art3.prop4}).


\bigskip



\begin{thebibliography}{VD}

\bibitem{Baa1} {\sc S. Baaj}, Multipicateurs non born\'{e}s. Th\`ese 3\`eme cycle, Universit\'{e} Paris 6 (1980).

\bibitem{Baa} {\sc S. Baaj \& P. Julg}, Th\'eorie bivariante de Kasparov et op\'erateurs non born\'es dans les \cst-modules hilbertiens.
{\it C.R Acad. Sci. Paris} {\bf 296} (1983), 875--878.


\bibitem{Bus1} {\sc  R.C. Busby},
Double centralizers and extensions of \cst-algebras.
{\it Trans. Amer. Math. Soc.} {\bf 132} (1968), 79--99.

\bibitem{Conw} {\sc J.B. Conway}, A Course in Functional Analysis. Springer-Verlag. (1985)

\bibitem{Kus} {\sc J. Kustermans}, One-parameter representations on \cst-algebras.
{\it Preprint Odense Universitet} (1997)

\bibitem{Kus1} {\sc J. Kustermans}, Regular \cst-valued weights on \cst-algebras.
{\it Preprint K.U. Leuven} (1997).


\bibitem{Lan} {\sc C. Lance},
Hilbert \cst-modules. Cambridge University Press (1995).

\bibitem{Wor7} {\sc K. Napi\'{o}rkowski \& S.L. Woronowicz} Operator theory in the \cst-algebra framework. To appear in {\it Reports Math. Phys.}


\bibitem{Sak}  {\sc   S. Sakai},
\cst-algebras and \wst-algebras.
Springer-Verlag, Berlin,  1971.


\bibitem{Stra} {\sc S. Stratila \& L. Zsid\'{o}},
Lectures on von Neumann algebras. {\it Abacus Press, Tunbridge Wells, England} (1979).

\bibitem{Tay}  {\sc D.C. Taylor}, The Strict Topology for Double Centralizer Algebras.
{\it Trans. Am. Math. Soc}{\bf 150} (1970), 633 -- 643

\bibitem{Tak}  {\sc  M. Takesaki},
Theory of Operator Algebras I.
{\it Springer-Verlag, New York} (1979).


\bibitem{Wor6}  {\sc  S.L. Woronowicz}, Unbounded elements affiliated with \cst-algebras
and non-compact quantum groups.
{\it Commun. Math. Phys.} {\bf  136} (1991),  399--432.

\bibitem{Wor5} {\sc S.L. Woronowicz},  \cst-algebras generated by unbounded elements.
{\it Preprint University of Warsaw}   



\bibitem{Zsido} {\sc I. Cior\v{a}nescu and L. Zsid\'{o}}, Analytic generators for one-parameter groups. {\it T\^{o}hoku Math. Journal} 
{\bf 28} (1976), 327--362.

\end{thebibliography}
\end{document}